\documentclass{article}
\usepackage{amsmath}
\usepackage{amsfonts}
\usepackage{color}      
\usepackage{epsfig}
\begin{document}

\newtheorem{theorem}{Theorem}
\newtheorem{lemma}{Lemma}
\newtheorem{corollary}{Corollary}
\newtheorem{definition}{Definition}
\newtheorem{conjecture}{Conjecture}
\newtheorem{identity}{Identity}
\newcommand{\Real}{{\bf R}}

\newcommand{\Ric}{\operatorname{Ric}}
\newcommand{\tr}{\operatorname{tr}}
\newcommand{\Hess}{\operatorname{Hess}}
\newcommand{\lift}{\operatorname{lift}}
\newcommand{\proj}{\operatorname{proj}}
\newcommand{\Comp}{\operatorname{Comp}}
\newcommand{\Perp}{\operatorname{perp}}
\newcommand{\data}{\operatorname{data}}
\newcommand{\graph}{\operatorname{graph}}
\newcommand{\II}{\text{I\hspace{-1 pt}I}}

\title{P.D.E.'s which Imply the Penrose Conjecture}

\author{Hubert L. Bray\thanks{Mathematics Department, Duke University,
Box 90320, Durham, NC 27708.  Supported in part by NSF grant
\#DMS-0706794.} \text{ }\text{ } \& \text{ } Marcus A.
Khuri\thanks{Mathematics Department, Stony Brook University, Stony
Brook, NY  11794. Supported in part by NSF grant \#DMS-0707086 and
a Sloan Fellowship.}}

\maketitle

\begin{abstract}

In this paper, we show how to reduce the Penrose conjecture to the
known Riemannian Penrose inequality case whenever certain
geometrically motivated systems of equations can be solved.
Whether or not these special systems of equations have general
existence theories is therefore an important open problem.  The
key tool in our method is the derivation of a new identity which
we call the generalized Schoen-Yau identity, which is of
independent interest. Using a generalized Jang equation, we
propose canonical embeddings of Cauchy data into corresponding
static spacetimes. In addition, our techniques suggest a more
general Penrose conjecture and generalized notions of apparent
horizons and trapped surfaces, which are also of independent
interest.

\end{abstract}

\section{Introduction}

In addition to their intrinsic geometric appeal, the Penrose
conjecture \cite{PENROSE2} and the positive mass theorem
\cite{SCHOENYAU3} are fundamental tests of general relativity as a
physical theory. In physical terms, the positive mass theorem
states that the total mass of a spacetime with nonnegative energy
density is also nonnegative. The Penrose conjecture, on the other
hand, conjectures that the total mass of a spacetime with
nonnegative energy density is at least the mass contributed by the
black holes in the spacetime.  In this section, we will explain
how these simple physical motivations translate into beautiful
geometric statements.

After special relativity, Einstein sought to explain gravity as a
consequence of the curvature of spacetime caused by matter.  In
contrast to Newtonian physics, gravity is not a force but instead
is simply an effect of this curvature.  As an analogy, consider a
heavy bowling ball placed on a bed which causes a significant
dimple in the bed.  Now roll a small golf ball off to one side of
the bowling ball.  Note that the path of the golf ball curves
around the bowling ball because of the curvature of the surface of
the bed.  In this analogy, the bowling ball represents the sun,
the golf ball represents the earth, and the surface of the bed
represents spacetime.  Whereas Newton explained the curvature of
the path of the smaller object by asserting an inverse square law
force of attraction between the two objects, Einstein declared
that the curvature of the smaller object's trajectory was due to
the curvature of spacetime itself, and that objects which did not
have forces (other than gravity) acting upon them followed
geodesics in the spacetime. That is, according to general
relativity, the sun and all of the planets are actually following
geodesics, curves with zero curvature, in the spacetime.

It should also be noted that general relativity is entirely
consistent with large scale experiments, whereas Newtonian physics
is not.  The most notable example may be the precession of the orbit
of Mercury around the Sun.  Whereas general relativity predicts the
rate at which the elliptical orbit precesses around the Sun to as
many digits as can be measured, Newtonian physics is off by almost
one percent, with all possible excuses for the discrepancy having
been eliminated.  The question, then, is how to turn the beautiful
and experimentally verified idea of matter causing curvature of
spacetime, which Einstein called his happiest thought, into a
precise mathematical theory.

First, assume that $(N^4,g_N)$ is a Lorentzian manifold, meaning
that the metric $g_N$ has signature $(-+++)$ at each point.  Note
that at each point, time-like vectors (vectors $v$ with $g_N(v,v) <
0$) are split into two connected components, one of which we will
call future directed time-like vectors, and the other of which we
will call past directed time-like vectors.

Next, define $T(v,w)$ to be the energy density going in the
direction of $v$ as measured by an observer going in the direction
of $w$, where $v,w$ are future-directed unit time-like vectors at
some point $p \in N$. In addition, suppose that $T$ is linear in
both slots so that T is a tensor.  Then the physical statement that
all observed energy densities are nonnegative translates into
\[T(v,w) \ge 0\]
for all future-directed (or both past-directed) time-like vectors
$v$ and $w$ at all points $p \in N$, known as the dominant energy
condition.

The goal, then, is to set $T$, which is called the stress-energy
tensor, equal to some curvature tensor.  A natural first idea is to
consider the Ricci curvature tensor since it is also a covariant
2-tensor.  In fact, this was Einstein's first idea.  However, the
second Bianchi identity on a manifold $N$ with metric tensor $g_N$
implies that
\[\mbox{div}(G) = 0,\]
where $G = Ric_N - \frac12 R_N g_N$, $Ric_N$ is the Ricci curvature
tensor, and $R_N$ is the scalar curvature.  This geometric identity
led Einstein to propose
\begin{equation}\label{Einsteinequation}
G = 8\pi T,
\end{equation}
known as the Einstein equation, since as an added bonus we
automatically get a conservation-type property for $T$, namely
$\mbox{div}(T) = 0$.  Naturally this is a very nice feature of the
theory since energy and momentum (the spatial components of the
energy vector) are conserved in every day experience.

The next step in pursuing this line of thought is to try to find
examples of spacetimes which satisfy the dominant energy condition,
the simplest case of which would be spacetimes with $G = 0$ which
are naturally called vacuum spacetimes.  Taking the trace implies
that such spacetimes (in 2+1 dimensions and higher) have zero scalar
curvature and therefore zero Ricci curvature as well.  The first
example (in 3+1 dimensions) is clearly Minkowski space
\[
   \left(\Real^4, -dt^2 + dx^2 + dy^2 + dz^2\right)
\]
which has zero Riemann curvature tensor. The second simplest example
of a spacetime with $G = 0$,
\begin{equation}\label{staticschwarz}
   \left(\Real \times (\Real^3 \setminus B_{m/2}(0)),
   -\left(\frac{1-\frac{m}{2r}}{1+\frac{m}{2r}}  \right)^2dt^2 +
   \left(1+\frac{m}{2r}\right)^4(dx^2 + dy^2 + dz^2)\right),
\end{equation}
where $r = \sqrt{x^2 + y^2 + z^2}$, is a one parameter family of
spacetimes called the Schwarzschild spacetimes.  When $m > 0$, these
spacetimes represent static black holes in a vacuum spacetime.

While the Schwarzschild spacetime can be covered by a single
coordinate chart (see Kruskal coordinates described in section
\ref{sec:coe}), the coordinate chart above only covers the exterior
region of the black hole and has a coordinate singularity (not an
actual metric singularity) on the coordinate cylinder $r = m/2$. For
our purposes, however, we will only be interested in the exterior
region of the Schwarzschild spacetime, which physically corresponds
to the region where observers have yet to pass into the event
horizon of the black hole, which is the point of no return from
which not even light can escape back out to infinity.

Spacetimes which may be expressed in the form
\[
   \left( \Real \times M, -\phi(x)^2 dt^2 + g\right),
\]
where $t \in \Real$, $x \in M$, and $g$ is a Riemannian (positive
definite) metric on $M$, are called static spacetimes.  This name is
appropriate since we see that the components of the spacetime metric
in this coordinate chart do not depend on $t$ but instead are
entirely functions of $x$.  Note also that static metrics are
defined not to have any time/spatial cross terms.  (Spacetimes which
allow time/spatial cross terms but where the metric components still
only depend on $x$ are called stationary spacetimes.)

An important result, first proved by Bunting and Masood-ul-Alam
\cite{BUNTINGALAM} using a very clever argument involving the
positive mass theorem, is that the only complete, asymptotically
flat \textit{static vacuum} spacetimes with black hole boundaries
(or no boundary) are the two spacetimes that we have listed so far,
Minkowski and Schwarzschild. This fact suggests that a thorough
understanding of these two spacetimes, including what makes them
special as compared to generic spacetimes, may be important for
understanding some of the most fundamental properties of general
relativity.

In fact, the Minkowski and Schwarzschild spacetimes are the extremal
spacetimes for the positive mass theorem and the Penrose conjecture,
respectively.  That is, the case of equality of the positive mass
theorem states that any space-like hypersurface of a spacetime
satisfying the hypotheses of the positive mass theorem which has $m
= 0$ can be isometrically embedded into the Minkowski spacetime.
Similarly, the case of equality of the Penrose conjecture (which,
while still a conjecture, has no known counter-examples in spite of
much examination) states that any space-like hypersurface of a
spacetime satisfying the hypotheses of the Penrose conjecture which
has $m = \sqrt{A/16\pi}$ (or to be more precise, the region outside
of the outermost minimal area enclosure of the apparent horizons)
can be isometrically embedded into the Schwarzschild spacetime.

Before we can state these theorems, though, we need to define a few
terms.  The basic object of interest in this paper is a space-like
hypersurface $M^3$ of a spacetime $N^4$, along with the induced
metric $g$ on $M^3$ and its second fundamental form $k$ in the
spacetime.

From this point on we will assume that $M^3$ has a global future
directed unit normal vector $n_{future}$ in the spacetime.  This
standard assumption is not necessary for stating the dominant
energy condition or, as we will see in section \ref{sec:genpen},
for defining generalized apparent horizons or stating the
generalized Penrose conjecture, but it is necessary for the
traditional definition for apparent horizons of black holes, as we
will see.  So for convenience, we will abuse terminology slightly
and also call
\begin{equation}\label{eqn:sffdef}
k(V,W) = -\langle \nabla_V W, n_{future}\rangle
\end{equation}
the second fundamental form of $M^3$, where $V,W$ are any vector
fields tangent to $M^3$ and $\nabla$ is the Levi-Civita connection
on the spacetime $N^4$. In this manner we are defining $k$ to be a
real-valued symmetric 2-tensor, where the true second fundamental
form, which takes values in the normal bundle to $M^3$, is $k \cdot
n_{future}$.

\begin{definition}
The triple $(M^3, g, k)$ is called the Cauchy data of $M^3$ for any
positive definite metric $g$ and any symmetric 2-tensor $k$.
\end{definition}

This name is appropriate because this is the data required to pose
initial value problems for p.d.e.'s such as the vacuum Einstein
equation $G=0$, or the Einstein equation coupled with equations
which describe how the matter evolves in the spacetime.  Note when
$M^3$ is flowed at unit speed orthogonally into the future that
\[
   \frac{d}{dt} g_{ij} = 2 k_{ij},
\]
so that $k$ is in fact the first derivative of $g$ in the time
direction (up to a factor).

Curiously, as we will see, the positive mass theorem and the Penrose
conjecture reduce to and are fundamentally statements about the
Cauchy data of space-like hypersurfaces of spacetimes, not the
spacetimes themselves.

\begin{definition}
At each point on $M^3$, define $\mu = T(n_{future},n_{future})$ to
be the energy density and the covector $J$ on $M^3$ to be the
momentum density, where $J(v) = T(n_{future},v)$, where $v$ is any
vector tangent to $M$.
\end{definition}

By the Einstein equation (equation \ref{Einsteinequation}) and the
Gauss-Codazzi identities \cite{ONEILL}, it follows that $\mu$ and
$J$ can be computed entirely in terms of the Cauchy data
$(M^3,g,k)$. In fact,
\begin{eqnarray}
 (8\pi) \; \mu \;=\; G(n_{future},n_{future})&=& (R + \mbox{tr}(k)^2 - \|k\|^2)/2
 \label{c1} \\
 (8\pi) \; J   \;=\; \hspace{.198in} G(n_{future},\cdot) \hspace{.198in}
 &=& \mbox{div}\left( k - \mbox{tr}(k) g\right),
 \label{c2}
\end{eqnarray}
where $R$ is the scalar curvature of $(M^3,g)$ at each point, and
the above traces, norms, and divergences are naturally taken with
respect to $g$ and the Levi-Civita connection of $g$.  Then the
dominant energy condition on $T$ implies that we must have
\begin{equation} \label{muj}
   \mu \ge |J|,
\end{equation}
which we will call the nonnegative energy density condition on
$(M^3, g, k)$, where again the norm is taken with respect to the
metric $g$ on $M^3$.

Equations \ref{c1} and \ref{c2} are called the constraint equations
because they impose constraints on the Cauchy data $(M^3,g,k)$ for
each initial value problem.  For example, we clearly need to impose
$\mu = 0$ and $J = 0$ on any Cauchy data which is meant to serve as
initial conditions for solving the vacuum Einstein equation $G = 0$.
However, for our purposes throughout the rest of this paper, we will
be interested in Cauchy data $(M^3,g,k)$ which only needs to satisfy
the nonnegative energy density condition in inequality \ref{muj}.
Since the assumption of nonnegative energy density everywhere is a
very common assumption, the theorems we prove will apply in a very
broad set of circumstances.

Next we turn our attention to the definition of the total mass of a
spacetime.  Looking back at the Schwarzschild spacetime, time-like
geodesics (which represent test particles) curve in the coordinate
chart as if they were accelerating towards the center of the
spacetime at a rate asymptotic to $m/r^2$ in the limit as $r$ goes
to infinity.  Hence, to be compatible with Newtonian physics (with
the universal gravitational constant set to 1) in the low field
limit, we must define $m$ to be the total mass of the Schwarzschild
spacetime.

More generally, consider any spacetime which is isometric to the
Schwarzschild spacetime with total mass $m$  for $r > r_0$ and which
is any smooth Lorentzian metric satisfying the dominant energy
condition on the interior region.  Of course the Schwarzschild
spacetime satisfies the dominant energy condition since it has $G =
0$.  Then the same argument as in the previous paragraph applies to
this spacetime, so its total mass must be $m$ as well.  This last
example inspires the following definition, which comes from
considering the $t=0$ slice of Schwarzschild spacetimes.

\begin{definition}
The Cauchy data $(M^3,g,k)$ will be said to be Schwarzschild at
infinity if $M^3$ can be written as the disjoint union of a compact
set $K$ and a finite number of regions $E_i$ (called \textit{ends}),
where $k=0$ on each end and each $(E_i,g)$ is isometric to
$\left(\Real^3 \setminus \bar{B}_{R_i}(0),
\left(1+\frac{m_i}{2r}\right)^4(dx^2 + dy^2 + dz^2)\right)$ for some
$m_i$ and some $R_i > \max(0,-m_i/2)$.  In addition, the mass of the
end $E_i$ will be defined to be $m_i$.
\end{definition}

We refer the reader to \cite{SCHOENYAU2} and \cite{WALD1} for more
general definitions of \textit{asymptotically flat} Cauchy data,
but for this paper the special case of being precisely
Schwarzschild at infinity is sufficiently interesting.

Typically we will be interested in Cauchy data with only one end.
However, sometimes it is convenient to allow for the possibility of
multiple ends.  Each end represents what we would normally think of
as a spatial slice of a universe, and the positive mass theorem and
the Penrose conjecture may be applied to each end independently.  In
fact, since ends can be compactified by adding a point at infinity
and then using a very large spherical metric on the end without
violating the nonnegative energy density assumption, without loss of
generality we may assume that any given Cauchy data has only one end
for the problems we will be considering.

\begin{theorem}\label{pmt}
(The Positive Mass Theorem, Schoen-Yau, 1981 \cite{SCHOENYAU2}; Witten, 1981 \cite{WITTEN}) \\
Suppose that the Cauchy data $(M^3,g,k)$ is complete, satisfies the
nonnegative energy density condition $\mu \ge |J|$, and is
Schwarzschild at infinity with total mass $m$.  Then
\[
   m \ge 0,
\]
and $m = 0$ if and only if $(M^3,g,k)$ is the pullback of the Cauchy
data induced on the image of a space-like embedding of $M^3$ into
the Minkowski spacetime.
\end{theorem}

The above theorem has an important special case when $k = 0$ which
is already extremely interesting.  Note that the nonnegative energy
condition reduces to simply requiring $(M^3,g)$ to have nonnegative
scalar curvature.

\begin{theorem}\label{rpmt}
(The Riemannian Positive Mass Theorem, Schoen-Yau, 1979 \cite{SCHOENYAU1}; Witten, 1981 \cite{WITTEN}) \\
Suppose that the Riemannian manifold $(M^3,g)$ is complete, has
nonnegative scalar curvature, and is Schwarzschild at infinity with
total mass $m$. Then
\[
   m \ge 0,
\]
and $m = 0$ if and only if $(M^3,g)$ is isometric to the flat metric
on $\Real^3$.
\end{theorem}

The adjective Riemannian was introduced by Huisken-Ilmanen in
\cite{HUISKENILMANEN} since the theorem is a statement about
Riemannian manifolds as opposed to Cauchy data in the more general
case.  We remind the reader that Cauchy data $(M^3,g,k)$ is still
required to have a positive definite metric $g$.

Notice that the Riemannian positive mass theorem is a beautiful
geometric statement about manifolds with nonnegative scalar
curvature.  In fact, Schoen-Yau were studying such manifolds
\cite{SCHOENYAU4} for purely geometric reasons when they first
realized that they could use minimal surface techniques to prove the
Riemannian positive mass theorem.  They then observed
\cite{SCHOENYAU2} that theorem \ref{pmt} (which is quite mysterious
from a geometric point of view without the physical motivation)
reduced to theorem \ref{rpmt} after solving a certain elliptic
p.d.e. on $(M^3,g,h)$ called the Jang equation, named after the
theoretical physicist who first introduced the equation in
\cite{JANG}.

Witten's proof of the positive mass theorem uses spinors and
proves both of the above statements by applying the
Lichnerowicz-Weitzenbock formula to a spinor which solves the
Dirac equation, and then integrating by parts.  This proof has a
strong appeal because it computes the total mass as an integral of
a nonnegative integrand.  However, so far it has not been clear
how to generalize this approach to achieve the Penrose conjecture,
although very interesting works in this direction include
\cite{HERZLICH} and \cite{KHURI2}.

Before we can state the Penrose conjecture, we need several more
definitions.  For convenience, we modify the topology of $M^3$ by
compactifying all of the ends of $M^3$ except for one chosen end by
adding the points $\{\infty_k\}$.  (However, the metric will still
not be defined on these new points.)

\begin{definition}\label{def:calS}
Define ${\bf\cal S}$ to be the collection of surfaces which are
smooth compact boundaries of open sets  $U$ in $M^3$, where $U$
contains the points $\{\infty_k\}$ and is bounded in the chosen end.
\end{definition}

All of the surfaces that we will be dealing with in this paper will
be in ${\cal S}$. Also, we see that all of the surfaces in ${\cal
S}$ divide $M^3$ into two regions, an inside (the open set) and an
outside (the complement of the open set). Thus, the notion of one
surface in ${\cal S}$ enclosing another surface in ${\cal S}$ is
well defined as meaning that the one open set contains the other.

\begin{definition}
Given any $\Sigma \in {\cal S}$, define $\tilde{\Sigma} \in {\cal
S}$ to be the outermost minimal area enclosure of $\Sigma$.
\end{definition}

That is, in the case that there is more than one minimal area
enclosure of the surface $\Sigma$, choose the outermost one which
encloses all of the others. The fact that an outermost minimal
area enclosure exists and is unique roughly follows from the
following: if $\partial A$ and $\partial B$ are both minimal area
enclosures of some surface, then so are $\partial(A \cup B)$ and
$\partial (A \cap B)$ since $|\partial(A \cup B)| + |\partial (A
\cap B)| = |\partial A| +|\partial B| = 2 A_{min}$ and both have
area at least $A_{min}$. A rigorous proof that the outermost
minimal area enclosure of a surface in an asymptotically flat
manifold exists and is unique is given in \cite{HUISKENILMANEN}.

\begin{definition}
Define $\Sigma \in {\cal S}$ in $(M^3,g,k)$ to be an apparent
horizon if it is one of the following three types of horizons,

\vspace{.1in} \noindent
a future apparent horizon if
\begin{equation}
H_\Sigma + \tr_\Sigma(k) = 0 \;\;\;\mbox{ on }\Sigma,
\label{futureah}
\end{equation}
a past apparent horizon if
\begin{equation}
H_\Sigma - \tr_\Sigma(k) = 0 \;\;\;\mbox{ on }\Sigma, \label{pastah}
\end{equation}
and a future and past apparent horizon if
\begin{equation} \label{futurepastah}
H_\Sigma = 0 \;\;\;\mbox{ and }\;\;\; \tr_{\Sigma}(k) = 0
\;\;\;\mbox{ on }\Sigma,
\end{equation}
where $H_\Sigma$ is the mean curvature of the surface $\Sigma$ in
$(M^3,g)$ (with the sign chosen to be positive for a round sphere in
flat $\Real^3$) and $\tr_{\Sigma}(k)$ is the trace of k restricted
to the surface $\Sigma$.
\end{definition}
Note that equation \ref{futurepastah} follows from assuming both
equations \ref{futureah} and \ref{pastah} everywhere on $\Sigma$.
Also note that $\Sigma$ is not required to be connected, although
from a physical point of view each component of $\Sigma$ is usually
thought of as the apparent horizon of a separate black hole.
Finally, observe that all three types of horizons are simply minimal
surfaces (surfaces with zero mean curvature) in the important
special case when $k=0$.

Physically, the only relevant apparent horizons for a spacecraft
flying around in a spacetime are future apparent horizons, because
spacecraft are only concerned about being trapped inside black holes
in the future. Mathematically, however, merely changing the choice
of global normal vector $n_{future}$ to $M^3$ in $N^4$ to
$-n_{future}$ changes the sign on $k$ which causes past apparent
horizons to become future apparent horizons, and vice versa.

Equations \ref{futureah}, \ref{pastah}, \ref{futurepastah} are
actually all conditions on the mean curvature vector of $\Sigma$ in
the spacetime. Note that at each point of $\Sigma^2$, the normal
bundle, of which the mean curvature vector is a section, is a
2-dimensional vector space with signature $(-+)$. Naturally, a basis
for this vector space is any outward future null vector along with
any outward past null vector. Since $\Sigma^2$ bounds a region in
$M^3$, outward is well-defined, and since there exists a global
normal vector $n_{future}$ to $M^3$, the future direction is
well-defined.

 Geometrically, if one flows a submanifold in the
normal directions $\vec{\eta}$, then the rate of change of the area
form of the submanifold is given by
\[
   \frac{d}{dt} dA = - \langle \vec\eta , \vec{H} \rangle dA
\]
where $\vec{H}$ is the mean curvature vector.  It turns out that the
mean curvature of a surface $\Sigma^2$ contained in a slice with
Cauchy data $(M^3,g,k)$ has coordinates $(\tr_{\Sigma}(k),-H)$,
where the first coordinate is in the unit future normal direction to
the slice and the second component is in the unit direction outward
perpendicular to the surface and tangent to the slice.  With this
convention, then the vector with components (1,1) is an outward
future null vector, and the vector with components (-1,1) is an
outward past null vector.

Hence, equation \ref{futureah} is equivalent to requiring that, at
each point of $\Sigma$, the dot product of the mean curvature vector
with any outward future null vector is zero (which implies that the
mean curvature vector is a real multiple of the outward future null
direction).  Similarly, equation \ref{pastah} is equivalent to
saying that the dot product of the mean curvature vector with any
outward past null vector is zero.  Thus, future apparent horizons
have the property that their areas do not change to first order when
flowed in outward future null directions.  The same is true for past
apparent horizons when flowed in outward past null directions.

We are now able to state the Penrose conjecture.  An excellent
survey of this conjecture is found in \cite{MARS}.

\begin{conjecture} \label{pconj}
(The Penrose Conjecture, 1973 \cite{PENROSE2} - Standard Version) \\
Suppose that the Cauchy data $(M^3,g,k)$ is complete, satisfies the
nonnegative energy density condition $\mu \ge |J|$, and is
Schwarzschild at infinity with total mass $m$ in a chosen end.  If
$\Sigma^2 \in {\cal S}$ is a future apparent horizon, then
\begin{equation}
   m \ge \sqrt{\frac{A}{16\pi}},   \label{pi}
\end{equation}
where $A$ is the area of the outermost minimal area enclosure
$\tilde{\Sigma}^2 = \partial U^3$ of $\Sigma^2$. Furthermore,
equality occurs if and only if $(M^3 \setminus U^3,g,k)$ is the
pullback of the Cauchy data induced on the image of a space-like
embedding of $M^3 \setminus U^3$ into the exterior region of a
Schwarzschild spacetime (which maps $\tilde\Sigma^2$ to a future
apparent horizon).
\end{conjecture}

Penrose's heuristic argument for a future apparent horizon in this
conjecture is described in more detail in \cite{BRAYCHRUSCIEL} and
\cite{MARS} but roughly goes as follows: If, as is generally
thought, asymptotically flat spacetimes eventually settle down to a
Kerr spacetime \cite{HEUSLER}, then in the distant future inequality
\ref{pi} will be satisfied since explicit calculation verifies this
fact for Kerr spacetimes, where $A$ is the area of the event
horizon.  Given that some energy may radiate out to infinity, the
total mass of these slices of Kerr may be less than the original
total mass. Also, by the Hawking area theorem \cite{HAWKINGELLIS}
(made more rigorous in \cite{CDGH}), and thus by the cosmic censor
conjecture \cite{PENROSE1} as well, the area of the event horizon is
nondecreasing in the spacetime evolution. Hence, this leads us to
conjecture inequality \ref{pi} in the initial Cauchy data slice, but
where $A$ is the total area of the event horizons of all of the
black holes.  The problem, though, is that, unlike apparent
horizons, event horizons are not determined by local geometry but
instead are defined in terms of which points in spacetime can
eventually escape out to infinity along future directed time-like
curves. Thus, in principle, there is no way to know which points
this includes without looking at the entire evolution of the
spacetime into the future.  However, in \cite{PENROSE2} Penrose
argued using the cosmic censor conjecture that future apparent
horizons, which are defined in terms of local geometry, must be
enclosed by event horizons. Thus, the area of $\tilde\Sigma$ serves
as a lower bound for the total area of the event horizons
\cite{JANGWALD}, \cite{HOROWITZ}, and the Penrose conjecture
follows. This same argument, but run in the opposite time direction,
yields the same conjecture for past apparent horizons as well. Thus,
in the conjecture one could replace ``future apparent horizon'' with
simply ``apparent horizon.''

It is also important to note, as Penrose did originally, that a
counterexample to the Penrose conjecture would be a very serious
issue for general relativity since it would imply that some part of
the above reasoning is false.  The consensus among many is that the
cosmic censor conjecture is the weakest link in the above argument.
If the cosmic censor conjecture turns out to be false, and naked
singularities (singularities not enclosed by the event horizons of
black holes) do develop in generic spacetimes, then this would
present a very interesting challenge to general relativity as a
physical theory.

However, like the positive mass theorem, setting $k = 0$ yields
another beautiful geometric statement about manifolds with
nonnegative scalar curvature, which is known to be true.

\begin{theorem} \label{rptheorem}
(The Riemannian Penrose Inequality, Bray, 2001 \cite{BRAY}) \\
Suppose that the Riemannian manifold $(M^3,g)$ is complete, has
nonnegative scalar curvature, and is Schwarzschild at infinity with
total mass $m$ in a chosen end. If $\Sigma^2 \in {\cal S}$ is a zero
mean curvature surface, then
\begin{equation}
   m \ge \sqrt{\frac{A}{16\pi}},   \label{rpi}
\end{equation}
where $A$ is the area of the outermost minimal area enclosure
$\tilde{\Sigma}^2 = \partial U^3$ of $\Sigma^2$. Furthermore,
equality occurs if and only if $(M^3 \setminus U^3,g)$ is isometric
to the Schwarzschild metric $\left(\Real^3 \setminus B_{m/2}(0),
\left(1+\frac{m}{2r}\right)^4(dx^2 + dy^2 + dz^2)\right)$.
\end{theorem}

In 1997, Huisken-Ilmanen proved a slightly weaker version of the
above result with the modification that $A$ is the area of the
\textit{largest connected component} of the outermost minimal area
enclosure of $\Sigma^2$ and with the additional assumption that
$H_2(M^3) = 0$. (This last topological condition can be replaced by
assuming that $\Sigma^2$ is already a connected component of the
outermost minimal area surface of $(M^3,g)$ by Meeks-Simon-Yau
\cite{MEEKSSIMONYAU}.)  Their method of proof, first proposed by the
theoretical physicists Geroch \cite{GEROCH} and Jang-Wald
\cite{JANGWALD}, uses a parabolic technique called inverse mean
curvature flow. Starting with a connected zero mean curvature
surface, Huisken-Ilmanen found a weak definition of inverse mean
curvature flow, where the surface is flowed out at each point in
$(M^3,g)$ with speed equal to the reciprocal of the mean curvature
of the surface at that point, for almost every surface in the flow.
Then they showed that the Hawking mass of the surface is
nondecreasing under this flow, equals the right hand side of the
Riemannian Penrose inequality initially, and limits to the left hand
side of the Riemannian Penrose inequality as the surface flows out
to large round spheres going to infinity. Both the physicists'
insight into proposing this idea and the mathematicians' cleverness
at generalizing the argument to something which could be made
rigorous are remarkably beautiful.

Bray's proof also involves a flow, but of the Riemannian manifold
$(M^3,g)$.  The flow of metrics stays inside the conformal class
of the original metric and eventually flows to a Schwarzschild
metric (shown as the case of equality metric).  The conformal flow
of metrics is chosen so as to keep the area of the outermost
minimal area enclosure of $\Sigma$ constant.  Also, the total mass
of the Riemannian manifold is nonincreasing by a clever argument
(first used by Bunting and Masood-ul-Alam in \cite{BUNTINGALAM})
using the positive mass theorem after a reflection of the manifold
along a zero mean curvature surface and a conformal
compactification of one of the resulting two ends.  Then since the
Schwarzschild metric gives equality in inequality \ref{rpi}, the
inequality follows for the original Riemannian manifold $(M^3,g)$.

All three systems of equations discussed in this paper which imply
the Penrose conjecture are based on a new geometric identity which
we call the generalized Schoen-Yau identity.  The identity is proved
in section \ref{sec:gsyi}, but with the lengthy computations
relegated to the appendices for readability.  This new idenity is a
generalization of equation 2.25 in Schoen-Yau's paper
\cite{SCHOENYAU2}.  The original Schoen-Yau identity was used to
reduce the positive mass theorem to the Riemannian positive mass
theorem by solving a p.d.e. called the Jang equation.  For all three
systems, our technique will involve a generalization of the Jang
equation to solving a system of two equations, the first of which is
a generalized Jang equation in all three cases. Rather than spending
time explaining the Jang equation, we will go straight to our proof
since the Jang equation appears as a special case of our method
(which for future reference is the case $\phi = 1$).

As a final comment, the Penrose conjecture can be generalized to a
statement about Cauchy data on $n$-manifolds motivated by
considering $(n+1)$-dimensional spacetimes, where $n \ge 3$.  In
fact, the positive mass theorem was proved by Schoen-Yau in
dimensions $n \le 7$ and by Witten in any number of dimensions, but
with the additional assumption that $M^n$ is spin.  The Riemannian
Penrose inequality was proved by Bray \cite{BRAY} in dimension $3$
using a proof which that author and Dan Lee \cite{BRAYLEE} have
generalized to manifolds in dimensions $n \le 7$, and in a slightly
weaker form by Huisken-Ilmanen \cite{HUISKENILMANEN} in dimension 3.
Since we will be reducing the general case of the Penrose conjecture
to the Riemannian Penrose inequality, the techniques presented here
have, at a minimum, the potential to address the Penrose conjecture
for manifolds with dimensions $n \le 7$. However, we will focus on
$n=3$ for simplicity.

\section{The Case of Equality}\label{sec:coe}

In this section we carefully study the case of equality of the
Penrose conjecture for the obvious reason that all of our estimates
used to prove the conjecture must give equality in these cases.
Also, we want to make sure that our techniques apply to all of the
case of equality examples as a check that we are not making
unjustified assumptions.

We refer the reader to \cite{ONEILL} for a discussion of the
Schwarzschild spacetime in Kruskal coordinates and follow those
conventions (except for the names of the two functions $\alpha$ and
$\beta$ defined in a moment).  Understanding the Schwarzschild
spacetime in Kruskal coordinates is essential since this is the
simplest global coordinate chart for the spacetime.  In Kruskal
coordinates, the entire Schwarzschild spacetime is expressed as the
subset $uv > -2m/e$ of $\Real^2 \times S^2$ with coordinates
$(u,v,\sigma \in S^2)$ and line element
\begin{equation}\label{kruskal}
   2 \beta(r) du dv + r^2 d\sigma^2,
\end{equation}
where $d\sigma^2$ is the standard round unit sphere metric on $S^2$,
$r>0$ is a function of $u,v$ determined by
\[
   uv = \alpha(r) = (r - 2m)e^{(r/2m) - 1},
\]
and
\[
   \beta(r) = (8m^2/r) e^{1 - (r/2m)}.
\]
The first quadrant region described by $u,v > 0$ is defined to be an
exterior region and is isometric to
\begin{equation}\label{standschwarz}
   \left(\Real \times (\Real^3 \setminus B_{2m}(0)) ,
   -\left(1-\frac{2m}{r} \right)dt^2 + \left(1-\frac{2m}{r}
   \right)^{-1}dr^2 + r^2 d\sigma^2
   \right)
\end{equation}
under the isometry
\[
u = \sqrt{\alpha(r)}e^{-t/4m}, \hspace{.5in} v
=\sqrt{\alpha(r)}e^{t/4m},
\]
which we leave as an exercise for the interested reader to check.
Note that we have now defined three different coordinate chart
representations for the Schwarzschild spacetime, the two above in
equations \ref{kruskal} and \ref{standschwarz}, and our original one
in equation \ref{staticschwarz}.

A key point is that two of these coordinate chart representations of
the exterior region of the Schwarzschild spacetime are written in
the form of a static spacetime.  For example, using the coordinates
in equation \ref{standschwarz}, the exterior region can be expressed
as
\begin{equation} \label{staticspacetimeorig}
   \left(\Real \times M^3, -\phi^2 dt^2 + g \right),
\end{equation}
where
\[
   \phi^2 = 1 - \frac{2m}{r},
\]
which of course gives us
\[
   r = \frac{2m}{1-\phi^2} \;\;\;\mbox{ and }\;\;\; r-2m =
   \frac{2m\phi^2}{1-\phi^2}.
\]
Hence, if we think of a slice of the static spacetime expressed as
the graph of $t = f(x)$ in the static spacetime, $x \in M^3$, we
see that
\begin{equation}\label{fff}
   f = 2m \log(v/u),
\end{equation}
\begin{equation}\label{ppp}
   \frac{2m\phi^2}{1-\phi^2} \; \exp\left({\frac{\phi^2}{1-\phi^2}}\right) =
   \alpha(r) = uv.
\end{equation}

The reason that these last two equations are important is that it
allows us to understand the behavior of $f$ and $\phi$ as they
approach the boundary of the exterior region $\{ x \;|\; u>0 ,
v>0\}$ of the Schwarzschild spacetime.  Our slice (intersected with
the exterior region of the Schwarzschild spacetime) has a future
apparent horizon boundary if $u=0$ everywhere on the boundary, a
past apparent horizon boundary if $v=0$ everywhere on the boundary,
and a future and past apparent horizon boundary if $u=v=0$
everywhere on the boundary.

The mixed case where $u=0$ on part of the boundary and $v=0$ on the
rest of the boundary does not represent a traditional apparent
horizon boundary.  However, we note that whenever this boundary is
area-outerminimizing, we are in fact in a case of equality of the
Penrose conjecture.  This observation helps motivate the definition
of a generalized apparent horizon in the next section.

Also, while the $u=0$ level set and the $v=0$ level set on $M^3$ are
both smooth (since the gradients of $u$ and $v$ on $M^3$ are never
zero since $M^3$ is space-like), the boundary of $\{ x \in M^3 \;|\;
u(x) > 0, v(x) > 0 \}$ in $M^3$ need not be smooth since the zero
level sets of $u$ and $v$ do not need to intersect smoothly.
However, when the boundary has corners it is never area
outerminimizing and thus not a case of equality of the Penrose
conjecture.

 We also note that apparent horizons
outside of the exterior region, say with $u=0$ but with $v < 0$ on
part of the apparent horizon, are not area-outerminimizing since
they have negative mean curvature at some points.  Consequently,
these last apparent horizons are enclosed by surfaces with less area
and are therefore not cases of equality of the Penrose conjecture
either.

Kruskal coordinates reveals that the Schwarzschild spacetime is
smooth on the boundary of the exterior region and certainly does not
have any singularities there. However, static coordinate
representations of the Schwarzschild spacetime have coordinate chart
singularities there (which do not represent anything geometric or
physical).  Hence, while Kruskal coordinates $u$ and $v$ are smooth
on any slice, even up to the apparent horizon boundary, $f$ and
$\phi$ are not necessarily.

In fact, we see that $f$ goes to $\pm \infty$ logarithmically at the
apparent horizon boundary typically (when $u$ or $v$ goes to zero
and the other stays positive).  Also, $\phi^2$ vanishes on the
apparent horizon boundary only linearly if either $u$ or $v$ is
strictly positive, which means that the derivative of $\phi$ is
going to $\infty$.  However, in the future and past apparent horizon
boundary case where $u,v$ both go to zero, then $\phi^2$ vanishes
quadratically and $\phi$ is smooth up to the boundary.  It is also
true that $f$ is smooth up to the boundary in this case by
L'Hopital's rule since $u$ and $v$, which equal zero on the future
and past apparent horizon, have nonzero derivatives there (since the
hypersurface is space-like). These observations are helpful since we
will be dealing with slices of the exterior region of the
Schwarzschild spacetime viewed in static coordinates for the rest of
this paper.

\section{Generalized Apparent Horizons and the Generalized Penrose
Conjecture}\label{sec:genpen}

In this section we describe the most general version of the Penrose
conjecture that we believe to be true.  Naturally this is an
important question to consider since proofs of a conjecture may be
more easily found when the most natural version of the conjecture is
understood.

\begin{definition}
Define the smooth surface $\Sigma^2 \in {\cal S}$ in $(M^3,g,k)$ to
be a generalized apparent horizon if
\begin{equation}
   H_{\Sigma} = |\tr_{\Sigma}(k)|
\end{equation}
and a generalized trapped surface if
\begin{equation}
   H_{\Sigma} \le |\tr_{\Sigma}(k)|.
\end{equation}
\end{definition}
In terms of the mean curvature vector $\vec{H}$ of $\Sigma^2$ in
the spacetime, a generalized trapped surface is one where
$\vec{H}$ is not strictly inward space-like anywhere on $\Sigma$.
Also note that this definition of a generalized apparent horizon
does not need a globally defined future directed unit normal to
$M^3$ since the definition is unaffected by a change of sign of
the second fundamental form $k$.  A related class of surfaces,
referred to as ``$\ast$-surfaces", appears in a different context
in \cite{SENOVILLA}.

Referring to the previous section, note that any smooth slice $M^3$
of the Schwarzschild spacetime which smoothly intersects (which is
often not the case) with the boundary of the first quadrant $\{u \ge
0 \;,\; v \ge 0\}$ of Kruskal coordinates intersects in a
generalized apparent horizon. These generalized apparent horizons
also give equality in the Penrose conjecture, so it is natural to
include them in the statement of a generalized Penrose conjecture.
Also note that traditional apparent horizons, if they are not
already generalized apparent horizons, are at least always
generalized trapped surfaces.

Another consideration which leads to this definition of generalized
apparent horizons is the case when a surface with multiple connected
components is a future apparent horizon on some connected components
and a past apparent horizon on the others.  While Penrose's original
heuristic argument does not apply to this surface, the techniques
that we develop in this paper seem to apply perfectly well.  Thus,
we would like a generalized Penrose conjecture which includes this
case as well.

After a talk on generalized apparent horizons by the first author
at the Niels Bohr International Academy's program ``Mathematical
Aspects of General Relativity" in April 2008, Robert Wald posed
the following insightful question: Is it possible for generalized
trapped surfaces to exist as boundaries of space-like slices of
Minkowski space? (A similar query was posed by Mars and Senovilla
in \cite{MARSSENOVILLA}.)  This question raises the issue of
whether or not generalized trapped surfaces always yield a
positive contribution to the ADM mass, which of course is a
prerequisite for a generalized version of the Penrose conjecture.
If one could find a generalized trapped surface which was the
boundary of a space-like slice of Minkowski space, then the total
mass of the slice would be zero, making a Penrose-type inequality
for the surface impossible.

In response to this question, the second author of this paper showed
that no such generalized trapped surface in Minkowski space exists
\cite{KHURI}. Furthermore, he showed that Witten's proof of the
positive mass theorem also works for asymptotically flat manifolds
with generalized trapped surface boundary and gives a positive lower
bound on the total mass. This result suggests that generalized
trapped surfaces and generalized apparent horizons have some
physical significance in that such surfaces, along with nonnegative
energy density $\mu \ge |J|$ everywhere in the spacetime, always
imply that the total mass is positive.  Finding the best possible
lower bound on the total mass motivates conjecturing a generalized
Penrose inequality.

In addition, a discussion between the first author and Tom Ilmanen
led to two more conjectures about generalized apparent horizons and
generalized trapped surfaces, which are known to be true in the
special case $k=0$.  We are pleased that Michael Eichmair
\cite{EICHMAIR} has announced proofs of these two conjectures
(except for the topological part of conjecture \ref{TI-MSY}) using
elliptic techniques (whereas Ilmanen's original ideas used parabolic
techniques).  We omit the $n=2$ case in these next two conjectures
because they are less relevant for our present purposes, but we
understand that Eichmair's results apply there as well.

\begin{conjecture} (Tom Ilmanen, 2006)\label{TI} \\
Given complete, asymptotically flat Cauchy data $(M^n,g,k)$, $3\le
n\le 7$, with a generalized trapped surface $\Sigma^{n-1}$, then
there exists a unique outermost generalized trapped surface
$\bar\Sigma$ which is a generalized apparent horizon.
\end{conjecture}

\begin{conjecture}(Tom Ilmanen, 2006)\label{TI-MSY}\\
Furthermore, $\bar\Sigma$ is strictly area outerminimizing (every
other surface which encloses it has larger area), and for $n = 3$,
the region exterior to $\bar\Sigma$ is diffeomorphic to $\Real^3$
minus a finite number of disjoint closed balls.
\end{conjecture}

The above conjecture is a generalization of Meeks-Simon-Yau
\cite{MEEKSSIMONYAU}, which is the case when $k=0$.  The topological
conclusions of this last conjecture, like the original
Meeks-Simon-Yau result, make this conjecture particularly
interesting for its own sake as well as important for the Jang-IMCF
system of equations we will describe later in the paper. We also
encourage the reader to study the related theorems of Andersson and
Metzger \cite{ANDERSSONMETZGER} on future and past apparent
horizons, which are relevant for this discussion.

All together, these considerations lead us to make the following
generalized Penrose conjecture.  Since traditional apparent horizons
are always generalized trapped surfaces, this conjecture implies the
original Penrose conjecture.

\begin{conjecture} \label{gpconj} (The Generalized Penrose
Conjecture) \\
Suppose that the Cauchy data $(M^3,g,k)$ is complete, satisfies the
nonnegative energy density condition $\mu \ge |J|$, and is
Schwarzschild at infinity with total mass $m$ in a chosen end. If
$\Sigma^2 \in {\cal S}$ is a generalized trapped surface, then
\begin{equation}
   m \ge \sqrt{\frac{A}{16\pi}},
\end{equation}
where $A$ is the area of the outermost minimal area enclosure
$\tilde{\Sigma}^2 =
\partial U^3$ of $\Sigma^2$. Furthermore, equality occurs if and
only if $(M^3 \setminus U^3,g,k)$ is the pullback of the Cauchy data
induced on the image of a space-like embedding of $M^3 \setminus
U^3$ into the exterior region of a Schwarzschild spacetime (which
maps $\tilde\Sigma^2$ to a generalized apparent horizon).
\end{conjecture}

We note that this conjecture is true when $k=0$ by \cite{BRAY}.
In this case, $\Sigma$ has nonpositive mean curvature and acts as
a barrier to imply the existence of an outermost minimal area
enclosure of $\Sigma$ which is minimal.

It is important to note that conjectures \ref{TI} and \ref{TI-MSY}
imply that the generalized Penrose conjecture (and hence the
original Penrose conjecture) follows from the following important
case of the generalized Penrose conjecture.

\begin{conjecture} \label{gpconj1} (The Generalized Penrose
Conjecture - Outermost Case) \\
Suppose that the Cauchy data $(M^3,g,k)$ is complete, satisfies the
nonnegative energy density condition $\mu \ge |J|$, and is
Schwarzschild at infinity with total mass $m$ in a chosen end.
Suppose also that $\Sigma^2 = \partial U^3 \in {\cal S}$ is a
strictly area outerminimizing generalized apparent horizon, that no
other generalized trapped surfaces enclose it, and that $M^3
\setminus \overline{U^3}$ is diffeomorphic to $\Real^3$ minus a
finite number of disjoint closed balls. Then
\begin{equation}
   m \ge \sqrt{\frac{A}{16\pi}},
\end{equation}
where $A$ is the area of $\Sigma^2$. Furthermore, equality occurs if
and only if $(M^3 \setminus U^3,g,k)$ is the pullback of the Cauchy
data induced on the image of a space-like embedding of $M^3
\setminus U^3$ into the exterior region of a Schwarzschild spacetime
(which maps $\Sigma^2$ to a generalized apparent horizon).
\end{conjecture}

Most of the remainder of this paper will focus on proving the above
conjecture. Naturally, when attempting a difficult conjecture, it
makes sense to consider the simplest case which still captures the
essential subtleties of the problem, and the above conjecture is
arguably that case.  We will focus on the three dimensional case in
this paper, but the above conjecture is the same in higher
dimensions up to and including seven, but without any hypothesis on
the topology of $M^n \setminus \overline{U^n}$, where the
conjectured inequality becomes $m \ge c_n A^{(n-2)/(n-1)}$ - see
\cite{BRAYLEE}.

Furthermore, the condition of not having any generalized trapped
surfaces outside of $\Sigma$ may turn out to be very important.
There is reason to believe that the generalized Jang equation, which
is a p.d.e.~we will define later in the paper which is central to
all of our approaches, may blow up on surfaces with $|H_\Sigma| =
|\tr_\Sigma(k)|$ (which are clearly generalized trapped surfaces).
Since we are leaving the existence theory of our proposed p.d.e.'s
which imply the Penrose conjecture open, those who are considering
studying these existence theories will have to understand this
possible behavior carefully.

\section{Proof of the Penrose Conjecture in a Special Case \label{specialcasesection}}

In this section we will prove the Penrose conjecture, conjecture
\ref{pconj}, with two extra assumptions, and show how the conjecture
follows from the Riemannian Penrose inequality, theorem
\ref{rptheorem}.  This special case, where a correct approach is
quite clear, will help us motivate the general case which is not so
obvious.

A major hint in the statement of the Penrose conjecture is the case
of equality.  Since the Penrose conjecture is an equality for any
slice (space-like hypersurface) of the exterior region of the
Schwarzschild spacetime with an apparent horizon boundary, we know
that all of our techniques must preserve this equality in every
estimate we derive.

On the other hand, if we are given some Cauchy data $(M^3,g,k)$
which comes from a slice of the Schwarzschild spacetime, it may be
difficult to recognize it as such.  However, our techniques must
absolutely be able to recognize these Cauchy data as the instances
where we get equality in all of our inequalities.

More generally, suppose $(M^3,g,k)$ comes from a slice of the static
spacetime
\begin{equation} \label{staticspacetime}
   \left(\Real \times M^3, -\phi^2 dt^2 + \bar{g} \right),
\end{equation}
where $\phi$ is a real-valued function on $M$ and $\bar{g}$ is some
other Riemannian metric on $M$.  Notice that the Schwarzschild
spacetime can be expressed in the form of equation
\ref{staticspacetime}.  However, while the Schwarzschild spacetime
is vacuum (meaning it has zero Einstein curvature and consequently
zero Ricci curvature), we are making no such requirement on
$\left(\Real \times M^3, -\phi^2 dt^2 + \bar{g} \right)$.

Given a real-valued function $f$ on $M$, define the graph map
\begin{equation}\label{graphmap}
   F: M \mapsto \Real \times M
\end{equation}
where $F(x) = (f(x), x)$. Then a short calculation reveals that the
pullback of the induced metric on the image of $F$ in a coordinate
chart is $\bar{g}_{ij} - \phi^2 f_i f_j$, so setting
\[
   \bar{g}_{ij} = g_{ij} + \phi^2 f_i f_j,
\]
guarantees that the pullback of the induced metric on the image of
the graph map $F$ is precisely $g$. A similar type of calculation
(but which is much longer and so is carried out in the appendices)
yields that the pullback of the second fundamental form of the image
of the graph map $F$ in the static spacetime to $(M^3,g)$ is
\begin{equation}\label{a0}
   h_{ij} = \frac{\phi Hess_{ij}f + \phi_i f_j + f_i \phi_j}
            {\left(1 + \phi^2|df|^2_g \right)^{1/2}},
\end{equation}
where subscripts on $f$ and $\phi$ represent coordinate chart
partial derivatives and the Hessian of $f$ is taken with respect to
the metric $g$ (or the Levi-Civita connection of $g$ if one
prefers).  These considerations lead us to the following special
case of the Penrose conjecture which has an elegant and relatively
short proof using the Gauss-Codazzi identities and the Riemannian
Penrose inequality.

\begin{theorem}\label{specialcasetheorem}
The Penrose conjecture as stated in conjecture \ref{pconj} follows
for
\\ $(M^3,g,k)$ if there exist two smooth functions $f$ and $\phi$
on $M^3$ such that
\begin{equation}\label{a1}
   k_{ij} = h_{ij} = \frac{\phi Hess_{ij}f + \phi_i f_j + f_i \phi_j}
            {\left(1 + \phi^2|df|^2_g \right)^{1/2}}
            \;\;\mbox{ outside of }\; \Sigma
\end{equation}
and
\begin{equation}\label{a2}
  \phi = 0 \mbox{ on }\Sigma,
\end{equation}
where $\phi > 0$ outside of $\Sigma$ and $f$ has compact support.
\end{theorem}

\textit{Proof:}  We will reduce the Penrose conjecture on
$(M^3,g,k)$ to the Riemannian Penrose inequality on $(M^3,\bar{g})$.
To do this we need to show that \vspace{.1in}

$\bullet\;\;\;$ the scalar curvature $\bar{R}$ of $\bar{g}$ is
nonnegative and that

$\circ\;\;\;$ $\Sigma$ has zero mean curvature $\bar{H}$ in
$(M^3,\bar{g})$.  \vspace{.1in}

\noindent Then the fact that $\bar{g}$ measures areas to be at least
as large as $g$ does implies that the area of any surface in
$(M^3,\bar{g})$ is at least as large as the area of that same
surface in $(M^3,g)$.  Thus,
\[
   \bar{A} := |\tilde{\Sigma}_{\bar{g}}|_{\bar{g}} \ge |\tilde{\Sigma}_{\bar{g}}|_g
      \ge |\tilde{\Sigma}_g|_g =: A,
\]
where $\tilde{\Sigma}_{\bar{g}}$ and $\tilde{\Sigma}_g$ are the
outermost minimal area enclosures of $\Sigma$ in $(M^3,\bar{g})$ and
$(M^3,g)$, respectively.  Since $f$ has compact support, the masses
of the two manifolds are the same. Then by the Riemannian Penrose
inequality
\[
   m = \bar{m} \ge \sqrt{\frac{\bar{A}}{16\pi}} \ge
   \sqrt{\frac{A}{16\pi}},
\]
which proves that Penrose conjecture on $(M^3,g,k)$.  Thus, all that
is left to prove are the two bullet points ($\bullet$) and
($\circ$).  \vspace{.1in}

\textit{Proof of ($\circ$):} Since $\bar{g}_{ij} = g_{ij} + \phi^2
f_i f_j$ and $\phi = 0$ on $\Sigma$ and $\phi$ and $f$ are smooth,
the two metrics are the same up to first order on $\Sigma$. But the
mean curvature of a surface, which is the main term in the first
variation of area formula, only depends on the metric and the first
derivatives of the metric.  Hence, $\bar{H} = H$.

Since $\Sigma$ is an apparent horizon, $H = \pm
\mbox{tr}_\Sigma(k)$. But since $\phi = 0$ on $\Sigma$, derivatives
along $\Sigma$ of $\phi$ are zero as well, so our assumption on the
special form of $k$ in equation \ref{a1} implies that
$\mbox{tr}_\Sigma(k)$ = 0. Hence,
\[
   \bar{H} = H = \mbox{tr}_\Sigma(k) = 0.
\]
 \vspace{.1in}

\textit{Proof of ($\bullet$):}  Working inside of the static
spacetime in equation \ref{staticspacetime}, let $n$ be the future
pointing normal vector to the image of $M^3$ under the graph map $F$
from equation \ref{graphmap} and let $\bar{n}$ be the future
pointing normal vector to $M^3$ viewed as the $t = 0$ slice of the
static spacetime.  Then these two vector fields on hypersurfaces can
be extended to the entire spacetime by requiring that these extended
vector fields are invariant under translation in the time coordinate
(which is an isometry of the spacetime).

The trick is to compute $G(n,\bar{n})$ using the Gauss-Codazzi
identities, but in two different ways.  We are given the nonnegative
energy density condition on $(M^3,g,k)$ that $\mu \ge |J|$.  Since
we are in the very special case that $k$ actually equals the second
fundamental form $h$ of the graph, $(M^3,g,h)$ has $\mu \ge |J|$
too.  This is equivalent to saying that $G(n,w) \ge 0$ for all
future time-like vectors $w$ in the spacetime.  Letting $w =
\bar{n}$ thus implies that
\begin{equation}\label{ggg}
   G(n,\bar{n}) \ge 0.
\end{equation}
On the other hand, applying the Gauss-Codazzi identities to the
$t=0$ slice of the static spacetime gives us
\begin{eqnarray*}
 (8\pi) \; \bar\mu \;=\; G(\bar{n},\bar{n})  &=& (\bar{R} + \mbox{tr}(\bar{p})^2 - \|\bar{p}\|^2)/2  \\
 (8\pi) \; \bar{J} \;=\; \,G(\bar{n},\cdot)\, &=& \, \mbox{div}\left( \bar{p} - \mbox{tr}(\bar{p}) \bar{g}\right)
\end{eqnarray*}
where $\bar{p}$ is the second fundamental form of the $t=0$ slice,
which of course is zero by the time symmetry of the spacetime.
Hence, $(8\pi)\bar\mu = \bar{R}/2$ and $\bar{J} = 0$.  Thus, if we
let
\[
   n = \alpha \bar{n} + (\mbox{vector tangent to $t=0$ slice}),
\]
where $\alpha$ is a positive function on $M$, we have that
\begin{equation}\label{ggg2}
   G(\bar{n},n) = \alpha G(\bar{n}, \bar{n}) = \alpha\bar{R}/2.
\end{equation}
But $G$ is symmetric, so by inequality \ref{ggg}, $\bar{R} \ge 0$,
which completes the proof of ($\bullet$) and the proof of the
Penrose inequality in this special case.

The case of equality of the above theorem would follow from
conjecture \ref{coeconj} in section \ref{sec:exist1}.  We refer the
reader to that section for discussion on the case of equality since
the main purpose of this section was to motivate the identities
computed in the next section.

\section{The Generalized Schoen-Yau Identity}\label{sec:gsyi}

The proof of the Penrose conjecture in the special case presented
above suggests how the Gauss-Codazzi identities can be used to
compute a formula for the scalar curvature $\bar{R}$ of $\bar{g} = g
+ \phi^2 df^2$ in terms of the scalar curvature $R$ of $g$, the
graph function $f$, and the warping factor $\phi$.  In this section
we will derive this formula and then show how this formula leads to
an identity central to our approach to the Penrose conjecture.

From this point on we will abuse terminology slightly and always
refer to the image of the graph map $F(M)$ simply as $M$ and the
$t=0$ slice of the constructed spacetime as $\bar{M}$.  This
notation is convenient since then $(M,g)$ and $(\bar{M},\bar{g})$
are space-like hypersurfaces of the spacetime $\left(\Real \times
M^3, -\phi^2 dt^2 + \bar{g} \right)$.  Let $\pi: M \mapsto
\bar{M}$ be the projection map $\pi(f(x),x) = (0,x)$ to the $t=0$
slice of the spacetime.

\begin{figure}[htbp]
\begin{center}

\input{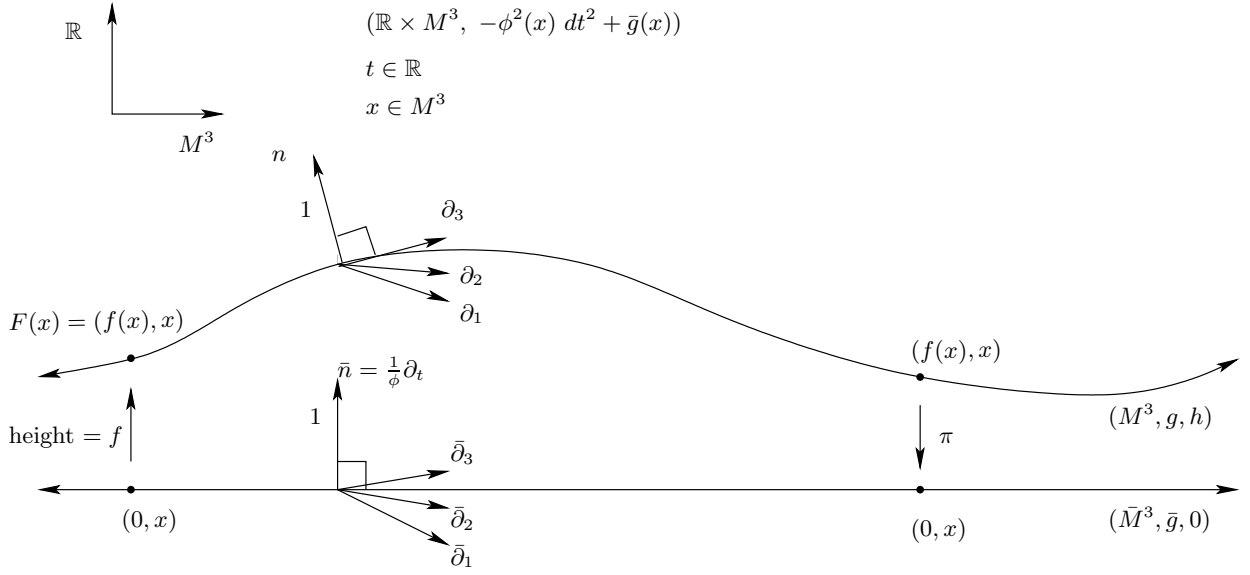}

\caption{Schematic diagram of the constructed static spacetime}
\label{figure1}
\end{center}
\end{figure}

Establishing some notation, let $\bar{\partial}_0 = \partial_t$ and
$\{\bar{\partial}_i\}$ be coordinate vectors tangent to $\bar{M}$.
Define
\begin{equation}
   \partial_i = \bar\partial_i + f_i\bar\partial_0
\end{equation}
to be the corresponding coordinate vectors tangent to $M$ so that
$\pi_*(\partial_i) = \bar\partial_i$. Then in this coordinate chart,
we have that
\begin{equation}\label{metric}
   g_{ij} = \bar{g}_{ij} - \phi^2 f_i f_j.
\end{equation}
It is convenient to write
\begin{equation}\label{inversemetric}
   g^{ij} = \bar{g}^{ij} + v^i v^j ,
\end{equation}
where
\begin{equation}\label{vdef}
   v^i = \frac{\phi f^{\bar{i}}}{(1-\phi^2|df|^2_{\bar{g}})^{1/2}} = \frac{\phi
   f^i}{(1+\phi^2|df|^2_g)^{1/2}} .
\end{equation}
We also define
\begin{equation}\label{vvecdef}
   \bar{v} = v^i \bar\partial_i \;\;\;\;\;\;\mbox{ and }\;\;\;\;\;\; v = v^i \partial_i
\end{equation}
so that $\pi_*(v) = \bar{v}$, and observe the useful identity
\begin{equation}
(1-\phi^2|df|^2_{\bar{g}}) \cdot (1+\phi^2|df|^2_g) = 1,
\end{equation}
which is evident by looking at the ratios of the volume forms.  See
appendix \ref{app:secondfundamentalform} for more discussion on
these calculations.

In this paper we use the convention that a barred index (as in
$f^{\bar{i}}$ above) denotes an index raised (or lowered) by
$\bar{g}$ as opposed to $g$.  That is, $f^{\bar{i}} =
\bar{g}^{ij}f_j$, where as usual $f_j = \partial f / \partial x_j$
in the coordinate chart. In general, barred quantities will be
associated with the $t=0$ slice $(\bar{M},\bar{g})$ and unbarred
quantities will be associated with the graph slice $(M,g)$.

In appendix \ref{app:identities} we compute that the second
fundamental form of the graph slice $(M,g)$ in our constructed
static spacetime is
\begin{eqnarray}
   h_{ij} &=& \frac{\phi \overline{\mbox{Hess}}_{ij}f + (f_i \phi_j + \phi_i f_j)
              - \phi^2 \langle df , d\phi \rangle_{\bar{g}} f_i f_j}
              { (1-\phi^2|df|^2_{\bar{g}})^{1/2}} \\
          &=& \frac{\phi \mbox{Hess}_{ij}f + (f_i \phi_j + \phi_i f_j)}
              { (1+\phi^2|df|^2_g)^{1/2}},
\end{eqnarray}
which we list now for future reference.

Finally, we extend $h$ and $k$ trivially in our constructed static
spacetime so that $h(\partial_t, \cdot) = 0 = k(\partial_t,\cdot)$
and such that these extended 2-tensors equal the original 2-tensors
when restricted to $M$.  Note that this gives
$h(\partial_i,\partial_j) = h(\bar\partial_i, \bar\partial_j)$, so
we can call this term $h_{ij}$ without ambiguity.  The same is true
for $k_{ij}$ and components of 1-forms like $f_i$ and $\phi_i$.
However, we remind the reader that the Hessian of a function, which
is the covariant derivative of the differential of a function,
depends on the connection and hence the metric since we will always
be using the respective Levi-Civita connections on $(M,g)$ and
$(\bar{M},\bar{g})$.

Now we are ready to proceed to compute a formula for $\bar{R}$. It
is a short calculation to verify that, in the constructed spacetime,
\begin{equation}
   \langle n , \bar{n} \rangle = - (1-\phi^2|df|^2_{\bar{g}})^{-1/2}
                               = - (1+\phi^2|df|^2_g)^{1/2}
\end{equation}
Thus,
\begin{equation}
   \bar{n} = (1+\phi^2|df|^2_g)^{1/2} n +
   \mbox{tan}_{\mbox{graph}}(\bar{n})
\end{equation}
where another short calculation reveals that
\begin{equation}
   \mbox{tan}_{\mbox{graph}}(\bar{n}) = -\phi f^j \partial_j
                                      = -\phi \nabla f.
\end{equation}

As in the previous section, the trick is to compute $G(n,\bar{n})$
two different ways using the Gauss-Codazzi identities.  As before,
applying these identities to the $t=0$ slice $(\bar{M},\bar{g})$
of the constructed spacetime gives us
\begin{eqnarray*}
   G(n,\bar{n}) &=& (1+\phi^2|df|^2_g)^{1/2} G(\bar{n},\bar{n}) \\
                &=& (1+\phi^2|df|^2_g)^{1/2} \cdot \bar{R}/2
\end{eqnarray*} since the $t=0$ slice has zero second fundamental
form.  On the other hand, applying the Gauss-Codazzi identities to
the graph slice $(M,g)$ yields
\begin{eqnarray*}
   G(n,\bar{n}) &=& (1+\phi^2|df|^2_g)^{1/2} G(n,n) +
                     G(n,\mbox{tan}_{\mbox{graph}}(\bar{n})) \\
&=& (1+\phi^2|df|^2_g)^{1/2}[R + (\mbox{tr}_g h)^2 - \|h\|^2_g]/2 \\
& &  + \mbox{div}(h - (\mbox{tr}_g h)g) (-\phi\nabla f).
\end{eqnarray*}
Combining the two previous equations, we get our first desired
result
\begin{equation}
   \bar{R} = R + (\mbox{tr}_g h)^2 - \|h\|_g^2
               + 2(d(\mbox{tr}_g h) - \mbox{div}(h))(v).
\end{equation}

Of course, what we are given in the hypotheses of the Penrose
conjecture is that $\mu \ge |J|_g$, where
\begin{eqnarray*}
 (8\pi) \; \mu \;=\; G(n,n)     &=& (R + \mbox{tr}_g(k)^2 - \|k\|^2)/2  \\
 (8\pi) \; J   \;=\; \,G(n,\cdot)\, &=& \mbox{div}(k) -
 d(\mbox{tr}_g(k)),
\end{eqnarray*}
for some symmetric 2-tensor $k$.   Hence,
\begin{eqnarray}\label{eqn:startofappendix}
\bar{R} &=& 16\pi(\mu - J(v))
        + (\mbox{tr}_g h)^2 -  (\mbox{tr}_g k)^2
        - \|h\|_g^2 + \|k\|_g^2   \notag \\ &&
        + 2 \, v(\mbox{tr}_g h) - 2 \, v(\mbox{tr}_g k)
        - 2 \, \mbox{div}(h)(v) + 2 \, \mbox{div}(k)(v).
\end{eqnarray}

Note that $\mu - J(v) \ge 0$ since $|v|_g \le 1$.  Hence, as we saw
in the previous section, if we can choose a $\phi$ and an $f$ so
that $h = k$, then we immediately get that $\bar{R} \ge 0$. However,
we are interested in investigating if a more general relationship
between $h$ and $k$ can give a similar result.

Our procedure is to convert our formula for $\bar{R}$ to an
expression in terms of the $\bar{g}$ metric.  Arguably $\bar{g}$ is
more natural than $g$ since it is the metric induced on the $t=0$
slice of the static spacetime.  To perform the conversion, we need
several identities for arbitrary symmetric 2-tensors $k$ which are
proven in appendix \ref{app:identities} and which we list here.

\begin{identity}\label{i1}
\[
   (\mbox{tr}_g(k))^2 - \|k\|_g^2 = (\mbox{tr}_{\bar{g}}k)^2 -
   \|k\|_{\bar{g}}^2 + 2 k(\bar{v},\bar{v}) \mbox{tr}_{\bar{g}}k -
   2|k(\bar{v},\cdot)|_{\bar{g}}^2
\]
\end{identity}

\begin{identity}\label{i2}
\[
   v(\mbox{tr}_g k) = \bar{v}(\mbox{tr}_{\bar{g}}k + k(\bar{v},\bar{v}))
\]
\end{identity}

\begin{identity}\label{i3}
\[
   \overline\Gamma_{ij}^k - \Gamma_{ij}^k = h_{ij}v^k - \phi f_i f_j
   \phi^{\bar{k}}
\]
\end{identity}

\begin{identity}\label{i4}
\begin{eqnarray*}
   \mbox{div}(k)(v) &=& \overline{\mbox{div}}(k)(\bar{v})
   + (\overline\nabla_{\bar{v}}k)(\bar{v},\bar{v})
   - 2 |\bar{v}|_{\bar{g}}^2\; k\left(\bar{v},\frac{\overline\nabla \phi}{\phi}\right)
   \\ &&
   + \langle h(\bar{v},\cdot) , k(\bar{v},\cdot) \rangle_{\bar{g}}
   + 2 h(\bar{v},\bar{v}) k(\bar{v},\bar{v})
   + (\mbox{tr}_{\bar{g}}h) k(\bar{v},\bar{v})
\end{eqnarray*}
\end{identity}

\begin{identity}\label{i5}
\[
  v_{\bar{\imath}\bar{;}j} = h_{ij} + v_{\bar{\imath}} h(\bar{v},\cdot)_j
                 - \frac{\phi_i v_{\bar{\jmath}}}{\phi}
\]
\end{identity}

\begin{identity}\label{i6}
\[
   \overline{\mbox{div}}(k)(\bar{v}) = \overline{\mbox{div}}(k(\bar{v},\cdot))
        - \langle h , k \rangle_{\bar{g}}
        - \langle h(\bar{v},\cdot) , k(\bar{v},\cdot) \rangle_{\bar{g}}
        + k\left(\bar{v},\frac{\overline\nabla \phi}{\phi}\right)
\]
\end{identity}

\begin{identity}\label{i7}
\[
   (\overline\nabla_{\bar{v}} k)(\bar{v},\bar{v}) = \bar{v}(k(\bar{v},\bar{v}))
      - 2 \langle h(\bar{v},\cdot) , k(\bar{v},\cdot) \rangle_{\bar{g}}
      - 2 h(\bar{v},\bar{v}) k(\bar{v},\bar{v})
      + 2 |\bar{v}|_{\bar{g}}^2 \;k\left(\bar{v},\frac{\overline\nabla\phi}{\phi}\right)
\]
\end{identity}

\begin{identity}\label{i8}
\begin{eqnarray*}
   \mbox{div}(k)(v) &=& \overline{\mbox{div}}(k(\bar{v},\cdot))
      + \bar{v}(k(\bar{v},\bar{v}))
      + k\left(\bar{v}, \frac{\overline\nabla\phi}{\phi}\right) \\&&
      - \langle h , k \rangle_{\bar{g}}
      - 2\langle h(\bar{v},\cdot) , k(\bar{v},\cdot) \rangle_{\bar{g}}
      + (\mbox{tr}_{\bar{g}}h) k(\bar{v},\bar{v})
\end{eqnarray*}
\end{identity}

Identities \ref{i1} and \ref{i2} are short calculations.  Identity
\ref{i3} is used in the proof of identity \ref{i4}. Plugging
identities \ref{i6} and \ref{i7} (which are proved using identity
\ref{i5}) into identity \ref{i4} results in identity \ref{i8}.
Finally, plugging identities \ref{i1}, \ref{i2}, and \ref{i8} into
our formula for $\bar{R}$ results in the main identity of this
paper.

\begin{identity}\label{mainidentity}
(\textbf{The Generalized Schoen-Yau Identity})
\begin{eqnarray*}
\bar{R} &=& 16\pi(\mu - J(v)) + \|h-k\|_{\bar{g}}^2+2|q|_{\bar{g}}^2
        - \frac{2}{\phi}\overline{\mbox{div}}(\phi q) \\&&
        + (\mbox{tr}_{\bar{g}}h)^2 - (\mbox{tr}_{\bar{g}}k)^2
        + 2 \bar{v}(\mbox{tr}_{\bar{g}}h - \mbox{tr}_{\bar{g}}k)
        + 2 k(\bar{v},\bar{v})(\mbox{tr}_{\bar{g}}h - \mbox{tr}_{\bar{g}}k)
\end{eqnarray*}
where
\[
   q = h(\bar{v}, \cdot) - k(\bar{v}, \cdot) = h(v, \cdot) - k(v, \cdot) \; .
\]
\end{identity}
Note that the two definitions of $q$ exist on the entire constructed
static spacetime and are equal since both $h$ and $k$ are extended
trivially in the constructed static spacetime. We also observe that
\[
   \frac{1}{\phi}\overline{\mbox{div}}(\phi q) = \mbox{div}_{ST}(q),
\]
where $\mbox{div}_{ST}$ is the divergence operator in the
constructed static spacetime.

In the special case that $\phi = 1$, the above identity was derived
by a different method by Schoen-Yau as equation 2.25 of
\cite{SCHOENYAU2} (in fact the procedure in \cite{SCHOENYAU2} may
also be used to obtain identity \ref{mainidentity} and will be
presented in a future paper). In that paper, they used the Jang
equation,
\[
0 = \mbox{tr}_{\bar{g}}(h - k)
\]
to reduce the positive mass theorem to the Riemannian positive
mass theorem.  While imposing the Jang equation in the special
case that $\phi = 1$ does not imply that $\bar{R} \ge 0$ as would
be most desirable, $\bar{R} \ge 2|q|_{\bar{g}}^2 - 2 \,
\overline{\mbox{div}}(q)$ implies that there exists a conformal
factor on $\bar{g}$ such that the conformal metric has nonnegative
scalar curvature and total mass less than or equal to that of
$\bar{g}$ and $g$.  Then the Riemannian positive mass theorem
applied to the metric conformal to $\bar{g}$ implies the positive
mass theorem on $(M,g)$.  This approach does not quite work for
the Penrose conjecture because the conformal factor needed to
achieve nonnegative scalar curvature changes the area of the
horizon in a way which is difficult to control.

\section{The Generalized Jang Equation}

All of the approaches to the Penrose conjecture that we consider in
this paper use the generalized Schoen-Yau identity.  This identity
plays a central role in the remainder of our discussions because it
directly relates the nonnegative energy condition on $(M^3,g,k)$
(which implies that $\mu \ge J(v)$ since $|v|_g < 1$) to the scalar
curvature of $(M^3,\bar{g})$.

Furthermore, this generalized Schoen-Yau identity strongly motivates
the generalized Jang equation,
\begin{equation}
   0 = \mbox{tr}_{\bar{g}}(h - k),
\end{equation}
which on the original manifold $(M^3,g)$ with Cauchy data
$(M^3,g,k)$ is the equation
\begin{equation}\label{gjecoords}
0 = \left(g^{ij} - \frac{\phi^2f^if^j}{1+\phi^2|df|_g^2}     \right)
 \left( \frac{\phi Hess_{ij}f + \phi_i f_j + f_i \phi_j} {\left(1
+\phi^2|df|^2_g \right)^{1/2}} - k_{ij} \right)
\end{equation}
when one substitutes the formulas for $h$ and $\bar{g}^{ij}$ in a
coordinate chart.  (In this paper we adopt Einstein's convention
that whenever there are both raised and lowered indices, summation
is implied, so the above formula is a summation over $i,j$ both
ranging from $1$ to $3$.)

Of course the original Jang equation, which again is the special
case $\phi(x) = 1$, only had one free function, $f$, whereas the
generalized Jang equation has two free functions, $f$ and $\phi$.
Hence, to get a determined system of equations, we need to specify
one more equation.  Later in the paper we will propose various
choices for this second equation, but our choice for the first
equation will always be the generalized Jang equation above.

Once the generalized Jang equation is specified, the generalized
Schoen-Yau identity simplifies greatly to
\begin{equation}\label{ssyi}
\bar{R} = 16\pi(\mu - J(v)) + \|h-k\|_{\bar{g}}^2 + 2|q|_{\bar{g}}^2
        - \frac{2}{\phi}\overline{\mbox{div}}(\phi q).
\end{equation}
It is important to note that the first three terms of the right hand
side of the above equation are all nonnegative since $\mu \ge |J|_g$
and $|v|_g < 1$.

\subsection{Boundary Conditions}\label{sec:boundaryconditions}

Examining the case of equality slices of the Schwarzschild spacetime
described in section \ref{sec:coe} leads us to propose the following
boundary conditions on generalized apparent horizons.  At a minimum,
these boundary conditions are satisfied almost everywhere for slices
of the exterior region of the Schwarzschild spacetime with
generalized apparent horizon boundaries.

\vspace{.1in} \hspace{.15in} {\bf Boundary Conditions on Generalized
Apparent Horizons}

\textit{ \\
Given a generalized apparent horizon $\Sigma$ with mean curvature
$H_\Sigma = |\tr_\Sigma^g(k)|$ and outward unit normal $\nu$ in
$(M^3,g)$, we require that $\phi = 0$ and
\begin{equation}
   \langle \nu,v \rangle_g = \mbox{sign}(\tr_\Sigma^g(k))
\end{equation}
on $\Sigma$, where as usual
\[
   v = \frac{\phi \nabla f}{(1 + \phi^2 |df|_g^2)^{1/2}}
\]
and $v$ is extended to the boundary $\Sigma$ by continuity.}

\noindent
Note that these boundary conditions are consistent with
$f$ blowing up to $+\infty$ where $\tr_\Sigma^g(k) < 0$, blowing
down to $-\infty$ where $\tr_\Sigma^g(k) > 0$, and $f$ staying
bounded where $\tr_\Sigma^g(k) = 0$ on $\Sigma$.

The hope is that these or similar boundary conditions imply that
$\Sigma$, which was a generalized apparent horizon in $(M^3,g)$,
becomes a minimal surface with
\begin{equation}\notag
   \bar{H}_\Sigma = 0
\end{equation}
in $(M^3,\bar{g})$.  We discuss the general case of this question in
appendix \ref{sec:transformationofthemeancurvature}.  For now, we
observe two important special cases.

The first important special case is when $\Sigma$ is a traditional
apparent horizon, either future or past, and $H_\Sigma > 0$.  If we
also assume that $f$ goes to $\pm \infty$ on each connected
component of $\Sigma$ in a reasonable fashion, then the level sets
of $f$ converge to $\Sigma$.  The formula for the mean curvature of
the level sets of $f$ in the new metric $\bar{g}$ is
\begin{eqnarray*}
   \bar{H} &=& (1 + \phi^2|df|^2)^{-1/2} H  \\ &=& (1-|v|_g^2)^{1/2} H
\end{eqnarray*}
since $\bar{g} = g + \phi^2 df^2$ does not change the metric on
the level sets of $f$, stretches lengths perpendicular to the
level sets of $f$ by a factor of $(1 + \phi^2|df|^2)^{1/2}$, and
by the first variation formula for area.  Then if we assume that
$f$ and $\phi$ behave similarly to the case of equality slices of
Schwarzschild, we get the following lemma.

\begin{lemma}
Suppose that $(M^3,g,k)$ has a smooth interior boundary $\Sigma$
which is a future [past] apparent horizon with $H_\Sigma>0$.  Then
if $f$ blows up [blows down] logarithmically, $|df|$ blows up
asymptotic to $1/s$, and $\phi^2$ goes to zero asymptotic to $s$
(where $s$ is the distance to $\Sigma$ in $(M^3,g)$), then the limit
of the mean curvatures $\bar{H}$ of the level sets of $f$ in
$(M^3,\bar{g})$ is zero.
\end{lemma}

The second important special case is the case where the boundary is
a future and past apparent horizon.  In this case, based on the case
of equality slices of Schwarzschild, we expect $f$ to stay bounded
and smooth and $\phi$ to stay smooth as well.

\begin{lemma}
Suppose that $(M^3,g,k)$ has a smooth interior boundary $\Sigma$
which is a future and past apparent horizon (which by definition has
$H_\Sigma = 0$ and $\tr_\Sigma^g(k) = 0$).  Then if $f$ is bounded
and smooth and $\phi$ is smooth and equals zero on $\Sigma$, then
$\bar{H}_\Sigma = H_\Sigma = 0$.
\end{lemma}
The proof of this lemma appeared in this paper already in section
\ref{specialcasesection}.  The point is that since $\bar{g} = g +
\phi^2 df^2$, both metrics $\bar{g}$ and $g$ are the same up to
first order on $\Sigma$ since $f$ and $\phi$ are smooth and $\phi =
0$ on $\Sigma$. Then since the mean curvature is only a function of
the metric and first derivatives of the metric, the two mean
curvatures are equal, and since $H_\Sigma=0$, both are zero.

\subsection{Blowups, Blowdowns, and Outermost Horizons}

One phenomenon of the original Jang equation ($\phi = 1$) is that
$f$ can blowup to $\infty$ on future apparent horizons or blowdown
to $-\infty$ on past apparent horizons, and this feature is still
present in the generalized Jang equation (given plausible
assumptions about the behavior of $f$ and $\phi$). More
importantly, according to \cite{SCHOENYAU2}, blowups and blowdowns
of $f$ with the original Jang equation can \textit{only} occur on
apparent horizons.

An important question, then, is to understand when blowups can
occur with the generalized Jang equation.  Certainly blowups and
blowdowns can still occur on traditional apparent horizons.  A
reasonable conjecture is that the blowup properties of the
generalized Jang equation are the same as the original Jang
equation as long as $\phi$ is smooth and strictly positive.
However, if $\phi$ is allowed to go to zero, then we have already
seen that $f$ can have a mixture of blowup, blowdown, and bounded
behavior on generalized apparent horizons in case of equality
slices of the Schwarzschild spacetime.  Those who study the
existence theories of the systems of equations proposed in this
paper will need to understand these issues.

A reasonable hope, however, is that as long as our boundary $\Sigma$
is already an outermost generalized apparent horizon, so that no
other generalized trapped surfaces enclose $\Sigma$, then $f$ stays
bounded away from $\Sigma$ as long as $\phi$ stays strictly
positive. A relevant calculation which is useful for studying the
question of when $f$ can blowup or blowdown is the following.

Recall the standard identity
\[
   \Delta f = \mbox{Hess} f (\nu, \nu) + H_\Sigma \cdot \nu(f) + \Delta_\Sigma f
\]
for the Laplacian of a function in terms of the Laplacian of that
function restricted to a hypersurface with mean curvature $H$ and
outward unit normal $\nu$.  If we let $\Sigma$ be any level set of
$f$, then we get that
\[
   \tr^g_\Sigma(\mbox{Hess} f) = \mp |df|_g H_\Sigma
\]
for blowup and blowdown respectively.  Thus, the generalized Jang
equation implies that
\begin{eqnarray*}
   0 &=& \tr_{\bar{g}}(h - k) \\
     &=& \bar{g}^{ij} (h_{ij} - k_{ij}) \\
     &=& \left[(g^{ij} -\nu^i\nu^j) + (\nu^i\nu^j - v^i
     v^j)\right](h_{ij} - k_{ij}) \\
     &=& \tr^g_\Sigma(h-k) +
     \frac{(h-k)(\nu,\nu)}{1+\phi^2|df|_g^2}\\
     &=& \frac{\phi \;\tr^g_\Sigma(\mbox{Hess}
     f)}{(1+\phi^2|df|_g^2)^{1/2}} - \tr^g_\Sigma(k) +
     \frac{(h-k)(\nu,\nu)}{1+\phi^2|df|_g^2}\\
     &=& \mp \frac{\phi |df|_g}{(1+\phi^2|df|_g^2)^{1/2}} H_\Sigma - \tr^g_\Sigma(k) +
     \frac{(h-k)(\nu,\nu)}{1+\phi^2|df|_g^2}\\
     &=&     \langle \nu,v \rangle_g    H_\Sigma - \tr^g_\Sigma(k) +
     \frac{(h-k)(\nu,\nu)}{1+\phi^2|df|_g^2}\\
\end{eqnarray*}
on level sets of $f$, where we have used the facts that $v$ and
$\nu$ are collinear,
\[
   |v|_g^2 = 1 - \frac{1}{1+\phi^2|df|_g^2},
\]
and the formulas for $\bar{g}^{ij}$ and $h_{ij}$ from section
\ref{sec:gsyi}.

\begin{lemma}\label{lemma:hbounded}
When $f$ is blowing up or blowing down, the term $h(\nu,\nu)$ is
bounded if $\phi^2 df$ is assumed to be smooth and nonzero in the
limit up to the boundary and $\phi = 0$ on the boundary, which is
true in case of equality slices of Schwarzschild.
\end{lemma}

\textit{Proof:} Referring back to section \ref{sec:coe}, in a smooth
slice of Schwarzschild $\phi^2 df$ can be expressed in terms of the
smooth Kruskal coordinate variables $u,v$.  Since by equation
\ref{fff}
\[
   df = 2m\left(\frac{dv}{v} - \frac{du}{u}\right)
\]
and by equation \ref{ppp}
\[
   \phi^2 = uv \gamma(uv)
\]
for some smooth function $\gamma \ne 0$ for $\phi^2 < 1$, we have
that
\[
   \phi^2 df = 2m \gamma(uv) \left(u dv - v du \right).
\]
The fact that our case of equality slices of Schwarzschild are
spacelike implies that $|du|_g \ne 0 \ne |dv|_g$, and the fact that
we are assuming blowup or blowdown implies that exactly one of $u,v$
is going to zero on the boundary.  Hence, not only is $\phi^2 df$
smooth up to the boundary, it is also nonzero in the limit up to the
boundary.

Then the fact that $h(\nu,\nu)$ is bounded up to the boundary
assuming this smoothness follows from the short calculation that
\begin{equation}
   h(\nu,\nu) = \frac{\nabla_\nu(\phi^2 df)(\nu)}{(\phi^2 + |\phi^2
   df|_g^2)^{1/2}},
\end{equation}
which completes the proof of the lemma.

Thus, referring back to our calculation before the lemma, given
blowup or blowdown with behavior on $f$ and $\phi$ as seen in the
case of equality slices of Schwarzschild in section \ref{sec:coe},
$\phi^2$ goes to zero linearly, $|df|_g$ goes to infinity like $1/s$
so that $\phi^2|df|_g^2$ goes to infinity like $1/s$, and
$h(\nu,\nu)$ stays bounded.  Then since $k$ is given to be smooth
and therefore bounded, we conclude that the generalized Jang
equation implies that
\[
   0 = \mp H_\Sigma - \tr^g_\Sigma(k)
\]
on surfaces with this type of blowup or blowdown of $f$, which of
course are the equations for future and past apparent horizons,
respectively.

\section{The Jang - Zero Divergence Equations}
\label{sec:exist1}

Looking at equation \ref{ssyi}, the most direct way to get $\bar{R}
\ge 0$ is to set $\overline{\mbox{div}}(\phi q) = 0$.  We will call
the resulting system of equations, equations \ref{direct:eqn1} and
\ref{direct:eqn2}, the Jang - zero divergence equations.  The
following existence conjecture for these equations implies the
outermost case of the generalized Penrose conjecture, conjecture
\ref{gpconj1}, using the Riemannian Penrose inequality, and is
therefore an important open problem.

\begin{conjecture}\label{exist1}
Given asymptotically flat Cauchy data $(M,g,k)$ with an outermost
generalized apparent horizon boundary $\Sigma$, there exists a
solution $(f, \phi)$ to the system of equations
\begin{eqnarray}
   0 &=& \mbox{tr}_{\bar{g}}(h - k)       \label{direct:eqn1}\\
   0 &=& \overline{\mbox{div}}(\phi q),   \label{direct:eqn2}
\end{eqnarray}
with $\lim_{x \rightarrow \infty} f(x) = 0$, $\phi^2 |\nabla f|^2 =
o(r^{-1})$, $\nabla(\phi^2 |\nabla f|^2) = o(r^{-2})$, \\
and $\lim_{x \rightarrow \infty} \phi(x) = 1$, where $\bar{g} = g +
\phi^2 df^2$,
\begin{equation}
 h = \frac{\phi \mbox{Hess}f + (df \otimes d\phi + d\phi \otimes df)}
              { (1+\phi^2|df|^2_g)^{1/2}} \; ,
\end{equation}
$q = h(v, \cdot) - k(v, \cdot)$, and $v = \phi \nabla f /
(1+\phi^2|df|^2_g)^{1/2}$, such that $\Sigma$ has zero mean
curvature in the $\bar{g}$ metric.
\end{conjecture}
The boundary conditions on $f$ which lead to $\Sigma$ having zero
mean curvature in $(M^3,\bar{g})$ are discussed in the previous
section and in appendix \ref{sec:transformationofthemeancurvature}.

Also, we comment that while equation \ref{direct:eqn2} is third
order in $f$, subtracting derivatives of equation \ref{direct:eqn1}
can remove the third order terms of $f$ in favor of Ricci curvature
terms.  While the resulting system has quadratic second order terms
in $f$ in the second equation, the system is degenerate elliptic.

The above system may also be reduced to a system of 1st order
equations by introducing new variables.  If we let $\alpha = df$,
then the above system has a solution whenever the first order system
with variables $\phi$ (a 0-form), $\alpha$ (a 1-form), and $\beta$
(a 2-form)
\begin{eqnarray}
   0 &=& d\alpha \\
   0 &=& d\beta \\
   0 &=& \mbox{tr}_{\bar{g}}(h - k) \\
   \phi q &=& d^{\bar{*}} \beta
\end{eqnarray}
has a solution, where $d^{\bar{*}}$ is the d star operator with
respect to the $\bar{g}$ metric which sends 2-forms to 1-forms. In
these variables, $\bar{g} = g + \phi^2 \alpha^2$,
\begin{equation}
h = \frac{\phi \nabla\alpha + (\alpha \otimes d\phi + d\phi \otimes
\alpha)} { (1+\phi^2|\alpha|^2_g)^{1/2}} \; ,
\end{equation}
$q = h(v, \cdot) - k(v, \cdot)$, and $v = \phi \vec\alpha /
(1+\phi^2|\alpha|^2_g)^{1/2}$, where $\vec\alpha$ is the dual vector
to $\alpha$ with respect to $g$.

\begin{theorem}
Conjectures \ref{exist1} and \ref{coeconj} (defined below to handle
the case of equality) imply conjecture \ref{gpconj1}, the outermost
case of the generalized Penrose conjecture.
\end{theorem}

{\it Proof:}  As was discussed in the previous section, the point of
requiring $\Sigma$ to be outermost in the above conjecture is so
that $f$ does not blowup or blowdown on the interior of $M$. Then
the method of proof assuming conjecture \ref{exist1} is basically
the same as the proof of the Penrose conjecture in the special case
in section \ref{specialcasesection}. The total mass of $(M^3,g)$ is
the same as the total mass of $(M^3,\bar{g})$ since the total mass
is defined in terms of the $1/r$ rate of decay of the metrics which
are equal since $||\bar{g} - g||_g = \phi^2 |df|_g^2$. Also, since
$\bar{g}$ measures lengths, areas, etc. to be greater than or equal
to that measured by $g$,
\begin{equation}
   \bar{A} = |\tilde{\Sigma}_{\bar{g}}|_{\bar{g}} \ge |\tilde{\Sigma}_{\bar{g}}|_g
      \ge |\tilde{\Sigma}_g|_g = A.
\end{equation}
Hence, the Penrose conjecture on $(M^3,g,k)$ follows from the
Riemannian Penrose inequality on $(M^3,\bar{g})$.

Thus, all that remains is to show that the Riemannian Penrose
inequality can be applied to $(M^3,g)$.  The existence theorem
already gives us that $\bar{H} = 0$, so the last thing to check is
that $\bar{R} \ge 0$, which follows directly from the generalized
Schoen-Yau identity and equation \ref{ssyi}.  This proves the
inequality part of the Penrose conjecture on $(M^3,g,k)$.

In the case of equality of the Penrose conjecture, clearly we must
have equality in all of our inequalities.  Since the case of
equality of the Riemannian Penrose inequality is solely when
$\bar{g}$ is the Schwarzschild metric which has zero scalar
curvature $\bar{R}$, equation \ref{ssyi} gives us
\begin{equation}\label{coejang}
   0 = 16\pi(\mu - J(v)) + \|h-k\|_{\bar{g}}^2 + 2|q|_{\bar{g}}^2.
\end{equation}
Since each of these three terms is nonnegative, each must be zero.
Hence, $k = h$.  If we could argue that $\phi = \phi_0$, where
$\phi_0$ is the warping factor from the Schwarzschild spacetime,
then we would have that $k = h$ is the second fundamental form  and
$g = \bar{g} - \phi^2 df^2$ is the induced metric of a slice of a
Schwarzschild spacetime, as desired.

However, there is a delicate point here.  In fact, $\phi$ does not
have to equal $\phi_0$ for $(M^3,g,k)$ to be the Cauchy data from a
slice of a Schwarzschild spacetime.  If $\phi = c \phi_0$ for some
constant $c
> 0$, then defining $df_0 = c \, df$ (which can be integrated to recover $f_0$)
implies that $(f_0, \phi_0)$ and $(f,\phi)$ produce the same metrics
and second fundamental forms. This may seem like a minor point at
first, but in fact this statement is still true if (and only if
which we leave as an exercise) $dc = 0$ on the open region $D$ where
$df \ne 0$.  Thus, $c$ may be different constants on each connected
component of $D$.  Again, $df_0$, which is still closed, may be
integrated to recover $f_0$ since the $t=0$ slice of Schwarzschild
is simply connected.  Thus, we have the following lemma.

\begin{lemma}\label{dk=0}
If $\phi = c \phi_0$, where $dc = 0$ on $\{x\;|\; df \ne 0\}$, then
$(M^3,g,k)$ comes from a slice of the Schwarzschild spacetime.
\end{lemma}

To prove the case of equality of the Penrose conjecture then, we
need to prove the hypotheses of the above lemma.  Looking back at
equation \ref{coejang}, we see that we must also have
\[
   0 = \mu - J(v) = (\mu - |J|) + |J|(1 - |v|) + (|J||v| - J(v)),
\]
where all norms are with respect to $g$. Again, since each of the
three grouped terms is nonnegative, all must be zero.  Since $|v|
< 1$, the second term equaling zero implies that $|J| = 0$ so that
the first term equalling zero implies that $\mu = 0$.

In the appendices we compute that in the static spacetime $\bar{g} -
\phi^2dt^2$,
\[
   n = (1-\phi^2|df|_{\bar{g}}^2)^{1/2} (\bar{n} + \phi \overline\nabla f)
\]
and that if $\bar{R} = 0$,
\[
   J = G(n, \cdot) =
   (1-\phi^2|df|_{\bar{g}}^2)^{1/2} \left[
   \overline{\mbox{Ric}} - \frac{\overline{\mbox{Hess}}\phi}{\phi} +
   \frac{\bar\Delta \phi}{\phi} \bar{g}
   \right](\phi \bar{\nabla}f , \cdot)
\]
where $\cdot$ is a tangent vector to the graph slice $(M^3,g)$ in
the first instance and its component tangent to the $t=0$ slice in
the second. Then since the Schwarzschild spacetime has $G = 0$,
$\bar{R} = 0$, $\bar\Delta \phi_0 = 0$, and
\[
\overline{\mbox{Ric}} = \frac{\overline{\mbox{Hess}}\phi_0}{\phi_0},
\]
$J=0$ in the case of equality implies the overdetermined equation
(when $df \ne 0$) for $\phi$ that
\begin{equation}\label{overdeterminedequation}
   0 =
   \left[\frac{\overline{\mbox{Hess}}\phi_0}{\phi_0} - \frac{\overline{\mbox{Hess}}\phi}{\phi} +
   \left(\frac{\bar\Delta \phi}{\phi} - \frac{\bar\Delta \phi_0}{\phi_0}\right) \bar{g}
   \right](\bar{\nabla}f , \cdot).
\end{equation}

\begin{conjecture} \label{coeconj}
Equation \ref{overdeterminedequation} implies the hypotheses of
lemma \ref{dk=0}.
\end{conjecture}

Clearly the hypotheses of lemma \ref{dk=0} imply equation
\ref{overdeterminedequation}, but we need the converse to be true as
well.  Assuming conjecture \ref{exist1} is true, a proof of
conjecture \ref{coeconj} would finish the case of equality part of
the outermost case of the generalized Penrose conjecture.

\section{Einstein-Hilbert Action Methods}\label{sec:EHM}

Equation \ref{ssyi} is a remarkable equation which deserves very
careful consideration.  Since we need a lower bound on the scalar
curvature $\bar{R}$ of $(M^3,\bar{g})$, the only troublesome term in
that equation is the last one, the divergence term. In the previous
section, we dealt with this last term by setting it equal to zero.
In this section, we make the natural observation that divergence
terms can also be dealt with by integrating them.

\begin{theorem}\label{thm:EHM}
If $\;\bar{g} = g + \phi^2df^2\;$ on $M^3$ with boundary $\Sigma^2$
and Cauchy data $(M^3,g,k)$ satisfying the nonnegative energy
condition $\mu \ge |J|$, the generalized Jang equation
\[
   0 = \tr_{\bar{g}}(h - k)
\]
is satisfied, and $f$ and $\phi$ behave at infinity and on the
boundary $\Sigma$ such that equations \ref{condition1} and
\ref{condition2} are satisfied as expected, then
\begin{equation}
   \int_M \bar{R} \phi \;\overline{dV} \ge 0,
\end{equation}
where $\bar{R}$ is the scalar curvature of $\bar{g}$ and
$\overline{dV}$ is the volume form of $\bar{g}$.
\end{theorem}

In other words, no matter what $\phi(x)$ is (as long as certain
boundary conditions are satisfied), the generalized Jang equation by
itself already gives a lower bound on the integral of the scalar
curvature of $\bar{g}$, weighted by $\phi$.  Of course the choice of
$\phi$ affects $f$ since $\phi$ appears in the generalized Jang
equation.

In the next couple of sections we discuss two different inequalities
of the form
\begin{equation}\label{form}
   \bar{m} - \sqrt{\frac{\bar{A}}{16\pi}} \ge \int_M Q(x)\bar{R}(x)\overline{dV},
\end{equation}
for some $Q(x) \ge 0$, where each inequality is based on one of the
two proofs of the Riemannian Penrose inequality. The expression for
$Q(x)$ differs in the two cases and will be described later.
However, if we then choose $\phi(x) = Q(x)$ to be our second
equation to be coupled with the generalized Jang equation, then
existence of such a system implies
\begin{equation}
   {m} - \sqrt{\frac{{A}}{16\pi}} \;\ge\;
   \bar{m} - \sqrt{\frac{\bar{A}}{16\pi}} \;\ge\; \int_M
   Q(x)\bar{R}(x)\overline{dV} \;=\;
    \int_M \bar{R} \phi \;\overline{dV} \;\ge\; 0,
\end{equation}
proving the corresponding form of the Penrose conjecture for the
original Cauchy data $(M^3,g,k)$.  We will call any method of proof
as above an Einstein-Hilbert action method.  So far we know of only
two methods of this form, the Jang-IMCF equations presented in the
next section, and the Jang-CFM equations discussed in the section
after that.

Note that any inequality of the form of inequality \ref{form} proves
the Riemannian Penrose inequality for $(M^3,\bar{g})$ as a special
case since then $\bar{R} \ge 0$ by hypothesis.  Thus, unless one
finds a new proof of the Riemannian Penrose inequality, the only way
to hope to prove an inequality like \ref{form} is either to adapt
currently known proofs of the Riemannian Penrose inequality, as we
are about to do in this paper, or to use the Riemannian Penrose
inequality itself. This last idea deserves additional consideration.

\textit{Proof of theorem \ref{thm:EHM}:}  Applying the divergence
theorem to equation \ref{ssyi} gives us that
\begin{eqnarray*}
    \int_M \bar{R} \phi \;\overline{dV} &\ge& 2\int_{\Sigma -
    S_\infty} \phi q(\bar\nu) \overline{dA}\\
    &=& 2 \int_{\Sigma} \phi (h - k)(\bar{v},\bar\nu) \overline{dA},\\
\end{eqnarray*}
where $k$ is assumed to converge to zero at infinity (or have
compact support), and $df$ is assumed to decay at least as fast as
$1/r^2$ at infinity (with reasonable bounds on $\mbox{Hess} f$ as
well) so that
\begin{equation}\label{condition1}
   0 = \lim_{r\rightarrow\infty} \int_{S_r} \phi (h - k)(\bar{v},\bar\nu)
   \overline{dA},
\end{equation}
where $\bar\nu$ is the unit outward normal vector to $\Sigma$  and
the sphere at infinity in $(
{M}^3,\bar{g})$.

In appendix \ref{sec:transformationofthemeancurvature}, we observe
that
\[
   \bar{\nu} = \left(\frac{1+\phi^2|df|_g^2}{1+\phi^2|(df|_\Sigma)|_g^2}\right)^{1/2}
   (\nu - \langle \nu,v \rangle v)
\]
and
\[
   \overline{dA} = (1+\phi^2|(df|_\Sigma)|_g^2)^{1/2} dA.
\]
Then since $1+\phi^2|df|_g^2 = 1/(1-|v|_g^2)$, we get that
\begin{eqnarray*}
    \int_M \bar{R} \phi \;\overline{dV}
    &\ge& 2 \int_{\Sigma} \phi (h - k)(\bar{v},\bar\nu) \overline{dA}\\
    &=& 2 \int_{\Sigma} \frac{\phi (h - k)(v,\nu - \langle \nu,v \rangle
    v)}{(1-|v|_g^2)^{1/2}} dA \\ &=& 0
\end{eqnarray*}
if we assume that
\begin{equation}\label{condition2}
   0 = \lim_{\Sigma_\epsilon \rightarrow \Sigma} \int_{\Sigma_\epsilon} \frac{\phi (h - k)(v,\nu - \langle \nu,v \rangle
    v)}{(1-|v|_g^2)^{1/2}} dA
\end{equation}
for some smooth family of surfaces $\Sigma_\epsilon$ converging to
$\Sigma$, proving the theorem.

In the case that $f$ is blowing up (or down) everywhere on $\Sigma$,
then choosing $\Sigma_\epsilon$ to be the level sets of $f$
simplifies things even more since then $v = \pm |v|\nu$. Then the
integrand becomes $\phi (1-|v|_g^2)^{1/2} (h-k)(v,\nu)$ and equals
zero with the usual boundary behavior since $\phi = 0$ on the
boundary, $|v|_g$ is going to one on the boundary, and both
$h(\nu,\nu)$ and $k(\nu,\nu)$ are bounded by lemma
\ref{lemma:hbounded}.  Equation \ref{condition2} is also clearly
satisfied in the case of a future and past apparent horizon where we
assume that $f$ and $\phi$ stay smooth and bounded, since all of the
terms in the integrand will be bounded, and $\phi = 0$ on the
boundary $\Sigma$.

As a final comment on theorem \ref{thm:EHM} before moving on, we
fully admit that a better understanding of the boundary behavior of
$f$ and $\phi$ is needed.  This better understanding should be able
to be achieved when the existence theories for the equations we are
proposing are discovered.

Before we get into applications of this theorem, it is worth noting
that
\[
  E(\bar{g},\phi) = \int_M \bar{R} \phi \;\overline{dV}
\]
is, up to a boundary term, the Einstein-Hilbert action of the
quotiented  static spacetime
\[
   \left(S^1 \times M, -\phi^2 dt^2 + \bar{g}\right),
\]
where we have turned the usual $\Real$ time coordinate into an $S^1$
of length one to get a finite integral.  The Einstein-Hilbert action
is defined to be the total integral of the scalar curvature $R^{ST}$
of the spacetime.  In the appendix we observe that
\[
   R^{ST} = \bar{R} - 2 \frac{\bar{\Delta}\phi}{\phi},
\]
and since $dV^{ST} = \phi \;\overline{dV}$, the Einstein Hilbert
action of the quotiented spacetime is
\begin{eqnarray*}
   \int_{S^1 \times M} R^{ST} dV^{ST}
   &=& \int_M (\bar{R}\phi - 2 \bar{\Delta}\phi) \;\overline{dV} \\
   &=& \int_M \bar{R}\phi \;\overline{dV} - 2
   \int_{\partial M} \langle\overline\nabla\phi,\bar\nu\rangle_{\bar{g}} \overline{dA}. \\
\end{eqnarray*}
We further observe that the boundary term vanishes when $\phi$ is
harmonic on $(M,\bar{g})$, as is the case in the Schwarzschild
spacetime.

Finally, the vacuum Einstein equation $G = 0$ is the
Euler-Lagrange equation which results from requiring a spacetime
to be a critical point of the Einstein-Hilbert action.  Since the
Minkowski and Schwarzschild spacetimes are the only vacuum static
spacetimes (with no boundary or black hole boundary)
\cite{BUNTINGALAM}, it follows that they are the only two static
spacetimes which are critical points of the Einstein-Hilbert
action, or equivalently $E(g,\phi)$, since boundary terms are
irrelevant for variations away from the boundary.

\section{The Jang-IMCF Equations}\label{sec:JangIMCF}

In this section we show how inverse mean curvature flow in
$(M^3,\bar{g})$ can be used to determine a warping factor $\phi$ for
the generalized Jang equation to get a system of equations which,
when there are solutions, implies the Penrose conjecture for a
single black hole when $H_2(M^3) = 0$.  Alternatively, the method
presented here has the potential to address the outermost case of
the Penrose conjecture as stated in conjecture \ref{gpconj1}, with
the additional assumption that the outermost generalized apparent
horizon is connected.  We will call the system of equations we are
proposing in this section the Jang-IMCF equations. An important open
problem is to find an existence theory for these equations.

Before we state the Jang-IMCF equations, we need to review inverse
mean curvature flow.  As introduced by Geroch \cite{GEROCH} and
Jang-Wald \cite{JANGWALD}, a smooth family of surfaces $\Sigma(t)$
in $(M^3,\bar{g})$ is said to satisfy inverse mean curvature flow
if the speed in the outward normal direction of the family of
surfaces as $t$ increases at each point is equal to $1/\bar{H}$,
where $\bar{H}>0$ is the mean curvature of the surface at that
point.  This flow has the important and surprising property that
the Hawking mass of $\Sigma(t)$ is nondecreasing when
$(M^3,\bar{g})$ has nonnegative scalar curvature $\bar{R}$.

To be more precise, define the Hawking mass of a surface $\Sigma$ in
$(M^3,\bar{g})$ to be
\[
   m_H(\Sigma) = \sqrt{\frac{|\Sigma|_{\bar{g}}}{16\pi}}\left(1 -
   \frac{1}{16\pi} \int_\Sigma \bar{H}^2 \overline{dA}
   \right),
\]
where all quantities are computed in $(M^3,\bar{g})$.  Then we can
compute the rate of change of the Hawking mass of a surface when
flowed out orthogonally with speed $\eta = 1/\bar{H}$ in
$(M^3,\bar{g})$ by using the first variation formula
\[
   \frac{d}{dt}(\overline{dA}) = (\eta \bar{H}) \overline{dA} = \overline{dA},
\]
the second variation formula
\[
   \frac{d}{dt} \bar{H} = -\Delta \eta - \|\overline{\II}\|_{\bar{g}}^2 \eta -
   \overline{\mbox{Ric}}(\bar\nu, \bar\nu) \eta,
\]
and the Gauss equation
\[
   \overline{\mbox{Ric}}(\bar\nu, \bar\nu) = \frac{1}{2} \bar{R} -
   \bar{K} + \frac{1}{2}\bar{H}^2 - \frac{1}{2} \|\overline{\II}\|_{\bar{g}}^2,
\]
where $\overline{\II}$ is the second fundamental form of $\Sigma$ in
$(M^3,\bar{g})$, $\overline{\mbox{Ric}}$ is the Ricci curvature of
$(M^3,\bar{g})$, and $\bar{K}$ is the Gauss curvature of $\Sigma$,
to get
\[
   \frac{d}{dt} \left(m_H(\Sigma(t))\right) = \sqrt{\frac{|\Sigma(t)|_{\bar{g}}}{16\pi}}
   \left[ \frac12 + \frac{1}{16\pi}\int_{\Sigma(t)} \frac{2|\overline\nabla \bar{H}|_{\bar{g}}^2}{\bar{H}^2}
   + \bar{R} - 2\bar{K} + \|\overline{\II}\|_{\bar{g}}^2 - \frac12 \bar{H}^2
   \right].
\]
The essential assumption that $\Sigma(t)$ is connected is used to
conclude that
\[
   \int_{\Sigma(t)} \bar{K} \overline{dA} = 2\pi \chi(\Sigma(t)) \le
   4\pi
\]
by the Gauss-Bonnet formula, which, along with
\[
   \|\overline{\II}\|_{\bar{g}}^2 \ge \frac12 \tr(\overline{\II})^2 =
   \frac12 \bar{H}^2
\]
allows us to conclude that
\[
   \frac{d}{dt} \left(m_H(\Sigma(t))\right) \ge \sqrt{\frac{|\Sigma(t)|_{\bar{g}}}{16\pi}}
   \int_{\Sigma(t)}  \frac{\bar{R}}{16\pi} \overline{dA}.
\]

If $(M^3,\bar{g})$ has nonnegative scalar curvature, then the above
equation implies that the Hawking mass of the smooth family of
surfaces determined by inverse mean curvature flow is nondecreasing.
However, we get a more general result if we integrate the above
equation in $t$ and use the co-area formula
\[
    \overline{dA}dt = \frac{1}{\eta}\overline{dV} = \bar{H} \overline{dV}
\]
to conclude that for a smooth family of surfaces satisfying inverse
mean curvature flow which foliates $M^3$, that
\begin{equation*}
   m_H(\Sigma(\infty)) - m_H(\Sigma(0)) = \int_{M} \bar{H} \sqrt{\frac{|\Sigma(t)|_{\bar{g}}}{16\pi}}
   \left(\frac{\bar{R}}{16\pi}\right) \overline{dV}
\end{equation*}
where at each point $x \in M^3$, $\bar{H}$ is the mean curvature
of the surface $\Sigma(t)$ through the point $x$.  Then assuming
that $\Sigma(0)$ has $\bar{H} = 0$ and is area outerminimizing and
that $(M^3,\bar{g})$ is sufficiently asymptotically flat, we
conclude our main result that
\begin{equation}
   \bar{m} - \sqrt{\frac{\bar{A}}{16\pi}} \ge \int_M Q(x)
   \left(\frac{\bar{R}(x)}{16\pi}\right)
   \overline{dV}(x)
\end{equation}
where $\bar{A}$ is the area of $\Sigma = \Sigma(0)$, the area
outerminimizing minimal boundary of $M^3$, and
\[
   Q(x) = \bar{H} \sqrt{\frac{|\Sigma(t)|_{\bar{g}}}{16\pi}}.
\]

More generally, Huisken-Ilmanen \cite{HUISKENILMANEN} observed
that there exists a weak notion of inverse mean curvature flow in
which the surfaces $\Sigma(t)$ jump outward to their outermost
minimal area enclosures whenever they are not already that
surface.  A key step in their approach is to represent the family
of surfaces $\Sigma(t)$ as the levels sets of a real-valued
function $u(x)$ on $M^3$ called the level set function.  Then if
\[
   \Sigma(t) = \partial \{ x \;|\; u(x) \le t \},
\]
it follows that $\eta = 1/|\overline\nabla u|_{\bar{g}}$ and
\[
   \bar{H} = \overline{\mbox{div}} \left(\frac{\overline\nabla u}{|\overline\nabla u|_{\bar{g}}}
   \right),
\]
so that inverse mean curvature flow on the level sets of $u(x)$ is
equivalent to
\begin{equation}\label{weakIMCF}
   \overline{\mbox{div}} \left(\frac{\overline\nabla u}{|\overline\nabla u|_{\bar{g}}}
   \right) = |\overline\nabla u|_{\bar{g}}.
\end{equation}
Huisken-Ilmanen then proceed to define a notion of weak solutions
to the above level set equation using an energy minimization
technique. These solutions have ``jump regions" where
$\overline\nabla u = 0$ corresponding to where the family of
surfaces $\Sigma(t)$ is not continuously varying but instead
``jumps" over these regions. Furthermore, Huisken-Ilmanen, using
elliptic regularization, proved that weak solutions of their
inverse mean curvature flow always exist.  We refer the reader to
their beautiful work \cite{HUISKENILMANEN}. However, using their
generalized inverse mean curvature flow, we achieve the following
theorem.

\begin{theorem}\label{thm:IMCFineq}
Given an asymptotically flat $(M^3,\bar{g})$ with $H_2(M^3) = 0$ and
a minimal connected boundary $\Sigma$ which bounds an interior
region, then
\begin{equation}\label{qqq}
   \bar{m} - \sqrt{\frac{\bar{A}}{16\pi}} \ge  \int_M Q(x)
   \left(\frac{\bar{R}(x)}{16\pi}\right)
   \overline{dV}(x)
\end{equation}
where $\bar{A}$ is the area of the outermost minimal area enclosure
$\tilde{\Sigma} = \partial U^3$ of $\Sigma$, $\bar{R}$ is the scalar
curvature and $\overline{dV}$ is the volume form of $\bar{g}$, and
\begin{equation}
   Q = |\overline\nabla u|_{\bar{g}}\sqrt{\frac{\bar{A}e^u}{16\pi}},
\end{equation}
where $u(x)$ is a weak solution to Huisken-Ilmanen inverse mean
curvature flow equalling zero on $\Sigma$.
\end{theorem}

\textit{Proof:}  As described in Huisken-Ilmanen's paper, if
$\Sigma$ is not already its own outermost minimal area enclosure, it
immediately jumps to it.  During this initial jump, but only on this
first jump, the area of the surface may decrease.  Hence, $\bar{A}$
must be defined to be the area of the outermost minimal area
enclosure $\tilde\Sigma$ of $\Sigma$, which also has zero mean
curvature by the maximum principle using $\Sigma$ as a barrier.
Also, since $H_2(M^3) = 0$, it follows that $\tilde{\Sigma}$ is
connected.  Since each component of $\tilde{\Sigma}$ bounds a
region, it follows that if $\tilde{\Sigma}$ did have more than one
connected component, all of the components except for one could be
removed (by either filling in holes or removing disconnected
regions), thereby decreasing the area. Then starting the flow at
$\tilde\Sigma$, our previous calculations generalize. The condition
that $H_2(M^3) = 0$ is also used to guarantee that each $\Sigma(t)$
is connected after each jump and therefore has Euler characteristic
$\le 2$ as needed in the computation of the rate of change of the
Hawking masses of $\Sigma(t)$.

By the first variation formula mentioned earlier in this section,
inverse mean curvature flow grows the area form exponentially.
Hence,
\[
   |\Sigma(t)|_{\bar{g}} = \bar{A} e^t.
\]
Thus, we have that
\begin{equation*}
   Q(x) = \overline{\mbox{div}} \left(\frac{\overline\nabla u}{|\overline\nabla u|_{\bar{g}}}
   \right)\sqrt{\frac{\bar{A}e^u}{16\pi}},
\end{equation*}
which equals the desired result by equation \ref{weakIMCF}.


Theorem \ref{thm:IMCFineq} deserves careful consideration.  In the
case that $\bar{R} \ge 0$, we recover a Riemannian Penrose
inequality for a single black hole. More generally, however, since
$\bar{R}/16\pi$ is energy density, we see that we have a kind of
integral of energy density on the right hand side of equation
\ref{qqq}, modified by the factor $Q(x)$.  On the flat metric on
$\Real^3$ and starting inverse mean curvature flow at a point, $Q(x)
= 1$, and on the Schwarzschild metric, $Q(x)$ is the harmonic
function going to one at infinity and equally zero on the minimal
neck.  This last fact, which can be verified by direct calculation,
will turn out to be important since this harmonic function also
equals the warping factor $\phi(x)$ in the Schwarzschild metric.

Theorem \ref{thm:IMCFineq} and theorem \ref{thm:EHM} together
motivate the system of equations,
\begin{eqnarray}
 0 &=& \mbox{tr}_{\bar{g}}(h - k)      \\
   |\overline\nabla u|_{\bar{g}} &=&
   \overline{\mbox{div}} \left(\frac{\overline\nabla u}{|\overline\nabla u|_{\bar{g}}}
   \right) \\
  Q &=& |\overline\nabla u|_{\bar{g}}\sqrt{\frac{\bar{A}e^u}{16\pi}} \\
   c \,\phi &=& Q,
\end{eqnarray}
where for our later convenience we choose $c =
\sqrt{\frac{\bar{A}}{16\pi}}$. The first equation is the generalized
Jang equation again. The second equation is the level set
formulation of inverse mean curvature flow on $(M^3,\bar{g})$. The
third equation is the definition of $Q(x)$ in terms of the inverse
mean curvature flow level set function $u(x)$. The new equation,
then, is the fourth equation, which sets $\phi(x)$ equal to $Q(x)$,
up to a constant. Then by theorems \ref{thm:IMCFineq} and
\ref{thm:EHM}, we conclude that
\begin{eqnarray*}
   \bar{m} - \sqrt{\frac{\bar{A}}{16\pi}} &\ge&  \int_M Q(x)
   \left(\frac{\bar{R}(x)}{16\pi}\right)
   \overline{dV}(x) \\
   &=& \frac{c}{16\pi} \int_M \bar{R} \phi
   \overline{dV} \;\;\ge\;\; 0.
\end{eqnarray*}

Then recalling that $\bar{g}$ measures areas at least as large as
$g$ does, we have that
\begin{equation*}
   \bar{A} = |\tilde{\Sigma}_{\bar{g}}|_{\bar{g}} \ge |\tilde{\Sigma}_{\bar{g}}|_g
      \ge |\tilde{\Sigma}_g|_g = A,
\end{equation*}
where again $\tilde{\Sigma}_{\bar{g}}$ is the outermost minimal area
enclosure of $\Sigma$ in $(M^3,\bar{g})$ and $\tilde\Sigma_g$ is the
outermost minimal area enclosure of $\Sigma$ in $(M^3,g)$.  Recall
also that since $H_2(M^3) = 0$, $\Sigma$ connected (and bounding a
region) implies that both $\tilde\Sigma_{\bar{g}}$ and
$\tilde\Sigma_g$ are also connected (and bound a region).  Hence, if
we can solve the above system with boundary conditions so that
$\Sigma$ has zero mean curvature in $(M^3,\bar{g})$ and so that the
total masses of $(M^3,g)$ and $(M^3,\bar{g})$ are the same, then we
would be able to conclude that
\begin{equation*}
   m = \bar{m} \ge \sqrt{\frac{\bar{A}}{16\pi}} \ge
   \sqrt{\frac{{A}}{16\pi}},
\end{equation*}
which would prove the Penrose conjecture for a single black hole in
the case that $H_2(M^3) = 0$.

In the case of equality in the above inequalities, $(M^3,\bar{g})$
has to be a time symmetric slice of the Schwarzschild spacetime by
the original Huisken-Ilmanen result.  Thus, inverse mean curvature
flow yields precisely the spherically symmetric spheres of
Schwarzschild, so $u$ is easy to compute.  Direct computation then
reveals that $Q(x)$ is the harmonic function in $(M^3,\bar{g})$
which equals zero on $\Sigma(0)$ and goes to one at infinity.  Since
$\phi$ equals $Q$ (up to a multiplicative constant, which is
irrelevant after a constant rescaling of the time coordinate in what
follows), we get that $\left(\Real \times M^3, -\phi^2 dt^2 +
\bar{g} \right)$ is isometric to a Schwarzschild spacetime.  Hence,
$g = \bar{g} - \phi^2 df^2$ is the induced metric on a slice of
Schwarzschild with graph function $f(x)$.  Finally, examining the
case of equality of theorem \ref{thm:EHM} (and that theorem's use of
the generalized Schoen-Yau identity) forces $\|h-k\|_{\bar{g}}^2 =
0$, which of course implies that $k_{ij} = h_{ij}$.  Hence, the
original Cauchy data $(M^3,g,k)$ is the induced Cauchy data on a
slice of a Schwarzschild spacetime with graph function $f(x)$.

Thus, understanding this system of equations, and whatever existence
theory might be associated with it, is a very interesting and
important open problem.  A first step is to observe that $Q$ does
not need to be defined in the system.  Hence, our system is
equivalent to

\vspace{.1in} \noindent {\bf The Jang - Inverse Mean Curvature Flow
Equations}
\begin{eqnarray}
 0 &=& \mbox{tr}_{\bar{g}}(h - k)      \\
   |\overline\nabla u|_{\bar{g}} &=&
   \overline{\mbox{div}} \left(\frac{\overline\nabla u}{|\overline\nabla u|_{\bar{g}}}
   \right) \\
   \phi &=& |\overline\nabla u|_{\bar{g}} \; e^{u/2}, \label{thirdeqn}
\end{eqnarray}
where we recall that
\begin{equation*}
\bar{g} = g + \phi^2 df^2
\end{equation*}
 and
\begin{equation*}
 h = \frac{\phi \mbox{Hess}f + (df \otimes d\phi + d\phi \otimes df)}
              { (1+\phi^2|df|^2_g)^{1/2}} \; ,
\end{equation*}
which can be thought of as three equations and three free functions
$f$, $u$, and $\phi$ on the original Cauchy data $(M^3,g,k)$.

In fact, the third equation, equation \ref{thirdeqn}, can be used to
solve for $\phi$ in terms of $u$, $du$, and $df$.  The purpose of
this is to recognize that the Jang-IMCF equations may also be
thought of as two equations and two free functions $f$ and $u$ once
we substitute for $\phi$.  Since only first derivatives of $f$ and
$u$ appear in the expression for $\phi$ below, the resulting
equivalent system is two second order equations in $f$ and $\phi$.

Unfortunately, the expression for $\phi$ in terms of $f$ and $u$ on
$(M^3,g)$ is a bit messy, but at least it is explicit.  From
equation \ref{thirdeqn}, we get
\begin{equation} \label{phiphi}
   \phi = |du|_{\bar{g}} \; e^{u/2}
\end{equation}
which is simple enough except that $\phi$ also appears in the
expression for $\bar{g}$.  Next we note that
\begin{eqnarray*}
   |du|_{\bar{g}}^2 &=& \bar{g}^{ij} u_i u_j \\
   &=& \left(g^{ij} - \frac{\phi^2 f^i f^j}{1 + \phi^2 |df|_g^2}\right) u_i
   u_j \\
   &=& |du|_g^2 - \frac{\phi^2 \langle df , du \rangle_g^2}{1 +
   \phi^2|df|_g^2},
\end{eqnarray*}
which, when combined with equation \ref{phiphi} gives us
\begin{equation*}
   \phi^2 = e^u \left( |du|_g^2 - \frac{\phi^2 \langle df , du \rangle_g^2}{1 +
   \phi^2|df|_g^2} \right).
\end{equation*}
It follows that $\phi^2$ solves the quadratic equation,
\begin{equation*}
   |df|_g^2 \cdot \phi^4 + B \cdot \phi^2 - e^u|du|_g^2 = 0,
\end{equation*}
where $B = 1 + e^u \left(  \langle df , du \rangle_g^2 - |df|_g^2
|du|_g^2 \right)$.  Thus,
\begin{equation}
   \phi^2 = \frac{-B + \sqrt{B^2 +
   4e^u|df|_g^2|du|_g^2}}{2|df|_g^2},
\end{equation}
which is clearly always nonnegative (and where we disregard the
negative square root in the quadratic formula since that solution is
nonpositive). Thus, an equivalent formulation of the Jang-IMCF
equations is
\begin{eqnarray}
 0 &=& \mbox{tr}_{\bar{g}}(h - k)   \label{gggjjjeee}   \\
   |\overline\nabla u|_{\bar{g}} &=&
   \overline{\mbox{div}} \left(\frac{\overline\nabla u}{|\overline\nabla u|_{\bar{g}}}
   \right)
\end{eqnarray}
where
\begin{equation}
\bar{g} = g + \phi^2 df^2,
\end{equation}
\begin{equation}
 h = \frac{\phi \mbox{Hess}f + (df \otimes d\phi + d\phi \otimes df)}
              { (1+\phi^2|df|^2_g)^{1/2}} \; ,
\end{equation}
\begin{equation}\label{phiformula}
   \phi = \frac{\sqrt{-\frac{B}{2} + \sqrt{\frac{B^2}{4} +
   e^u|df|_g^2|du|_g^2}}}{|df|_g},
\end{equation}
and
\begin{equation}
   B = 1 + e^u \left(  \langle df , du \rangle_g^2 - |df|_g^2 |du|_g^2
   \right),
\end{equation}
which can be thought of as two equations and two free functions $f$
and $u$ on the original Cauchy data $(M^3,g,k)$.  (In equation
\ref{phiformula}, $\phi = |du|_g e^{u/2}$ when $|df|_g = 0$ by
equation \ref{phiphi}).

We end this section with a general discussion of the some of the
challenges involved in finding an existence theory for the
Jang-IMCF equations. First, note that these equations reduce to
the Huisken-Ilmanen IMCF equation on $(M^3,g)$ when $\tr_g(k)=0$
since then we can choose $f=0$ (which implies $\bar{g} = g$) to
satisfy the generalized Jang equation (equation \ref{gggjjjeee}).
Thus, clearly a notion of a weak solution to this system of
equations is required. Furthermore, the notion of ``jumps" must
still be involved when there are regions in which $du = 0$.  Note
that when $du = 0$, then $\phi = 0$.  Thus, if $f$ stays smooth
and bounded, it would follow that $h = 0$, which means that the
generalized Jang equation cannot be solved unless $\tr_g(k) = 0$
in this region as well.  If $\tr_g(k) \ne 0$ in this region, then
this would suggest that $f$ needs to be unbounded or undefined in
this region.  Clearly this is an important issue to understand.

Given these and other considerations, one might be tempted to be
pessimistic about finding a general existence theory for the
Jang-IMCF equations.  In fact, it was once the case that most were
pessimistic about the original inverse mean curvature flow
proposed by Geroch \cite{GEROCH}, right up until Huisken-Ilmanen
\cite{HUISKENILMANEN} found an amazingly beautiful and natural
existence theory for a generalized version of inverse mean
curvature flow. Thus, there is also precedent for optimism.

\section{The Jang-CFM Equations}

In this section we comment that there is at least one other
Einstein-Hilbert action method in addition to the Jang-IMCF
equations.  So far we have seen how the Penrose conjecture would
follow from a general existence theory for the Jang-Zero Divergence
equations presented in section \ref{sec:exist1} or, for a single
black hole in dimension three, from a general existence theory for
the Jang-IMCF equations presented in section \ref{sec:JangIMCF}. In
this section, we briefly discuss a third system of equations whose
existence theory would also imply the Penrose conjecture.  The
precise statement of this third system is a bit laborious and so we
do not state it here, but only describe it and the additional
considerations it involves.

In section \ref{sec:EHM}, we explained how any inequality of the
form
\begin{equation}\label{form2}
   \bar{m} - \sqrt{\frac{\bar{A}}{16\pi}} \ge \int_M Q(x)\bar{R}(x)\overline{dV},
\end{equation}
for some $Q(x) \ge 0$, leads to a system of equations which implies
the Penrose conjecture.  The first equation in the system is the
generalized Jang equation and the second equation in the system is
simply $\phi(x) = Q(x)$ (times a constant if one likes). In section
\ref{sec:JangIMCF} we pursued this approach in detail for the
Huisken-Ilmanen inverse mean curvature flow.

Bray's proof \cite{BRAY} of the Riemannian Penrose inequality,
when revisited, also yields an inequality of the form of equation
\ref{form2}.  This proof of the Riemannian Penrose inequality
involves a conformal flow of metrics (CFM) which flows an initial
asymptotically flat metric with nonnegative scalar curvature to a
Schwarzschild metric in the limit as the flow parameter goes to
infinity.  Furthermore, the area of the horizon stays constant,
and (by the Riemannian positive mass theorem it turns out that)
the total mass is nonincreasing during the flow.

To generalize the conformal flow of metrics (CFM) proof to get an
inequality as in equation \ref{form2}, we first need to generalize
the positive mass theorem to get an inequality of the form
\begin{equation}\label{form3}
   \tilde{m} \ge \int_M Q(x)\tilde{R}(x)\tilde{dV},
\end{equation}
for some $Q(x) \ge 0$ on some $(M,\tilde{g})$.  Witten's spinor
proof \cite{WITTEN} of the Riemannian positive mass theorem provides
such a result, for example, at least whenever a spinor solution to
the Dirac equation exists (since we are not assuming $\tilde{R} \ge
0$ anymore, there is an issue now). Also, a result of this type can
be found by multiplying the metric $\tilde{g}$ by a conformal factor
to achieve zero scalar curvature globally (when such a factor
exists), and then measuring how much the mass changes.  This last
idea is made precise by Jauregui in \cite{JAUREGUI}.  Finally, one
can also use inverse mean curvature flow starting from any point to
prove an inequality of the above form in equation \ref{form3}.  This
third approach currently has the advantage over the first two in
that it is known to work in all cases in dimension three by the
previous section and the work of Huisken and Ilmanen
\cite{HUISKENILMANEN} and Streets \cite{STREETS}.

In the conformal flow of metrics $(M,g_t)$ with total masses $m(t)$,
\begin{equation}
   {m}'(t) = -\frac{1}{2} \tilde{m}(t).
\end{equation}
Hence,
\begin{equation}
   {m} - \sqrt{\frac{{A}}{16\pi}} \ge \int_0^\infty \frac12
   \tilde{m}(t) dt
\end{equation}
since the areas $A(t)$ of the horizons of $(M,g_t)$ stay constant
and the flow of metrics converges to Schwarzschild where ${m} -
\sqrt{\frac{{A}}{16\pi}} = 0$.  Then plugging equation \ref{form3}
into the above equation and accounting how the scalar curvature
transforms conformally gives a result of the desired form in
equation \ref{form2}. Hence, modulo the existence questions needed
to get an equality of the form of equation \ref{form3}, we get a
generalization of the Riemannian Penrose inequality.

One difference between the Jang-CFM equations and the Jang-IMCF
equations, however, is that the Jang-CFM equations are not local.
That is, $Q(x)$ in this case does not satisfy a local p.d.e.~at each
point and instead has a more complicated expression.  Hence, for the
Jang-CFM equations to have an existence theory, the theory would
have to work for a wide range of possible $Q$.  On the other hand,
the $Q(x)$ from the Jang-CFM equations has the potential to have
better regularity than the $Q(x)$ from the Jang-IMCF equations.  For
the Jang-IMCF equations, $c \phi = Q$ is not necessarily continuous
or even positive, and in fact equals zero in jump regions of the
inverse mean curvature flow on $(M^3,\bar{g})$, which as discussed
at the end of the previous section, introduces additional analytical
challenges.

\section{Open Problems}

The two most interesting and important open problems discussed in
this paper are finding an existence theory for the Jang-Zero
Divergence Equations (which would prove the Penrose conjecture) and
finding an existence theory for the Jang-IMCF equations (which would
prove the Penrose conjecture for a single black hole when $H_2(M^3)
= 0$).  Another interesting problem is to find a general existence
theory for any Einstein-Hilbert action method as long as the
associated $Q(x)$ has certain properties.  If the $Q(x)$ from the
Jang-CFM equations qualified for such a theory, this would also
prove the Penrose conjecture.

Another interesting problem is to find additional Einstein-Hilbert
action methods by finding new inequalities of the form of equation
\ref{form2}.  Since the special case of $\bar{R} = 0$ implies the
Riemannian Penrose inequality, one would either have to find a new
proof of the Riemannian Penrose inequality or use the Riemannian
Penrose inequality itself to find such a generalization.  There may
be reasonable ideas to try in this latter approach.

There is also the question of the physical interpretation of
inequalities of the form of equation \ref{form2}.  Rewriting the
inequality gives us
\begin{equation}
   \bar{m} \ge \sqrt{\frac{\bar{A}}{16\pi}} + \int_M Q(x)\bar{R}(x)\overline{dV},
\end{equation}
which could be interpreted as saying that the total mass of a
time-symmetric slice of a spacetime (not necessarily with
nonnegative energy density) is at least equal to the mass
contributed by the black holes (the first term) plus a weighted
integral of the energy density (the second term), since energy
density at each point can be interpreted as $\bar\mu =
\bar{R}/16\pi$. The purpose of $Q$ can be interpreted as the need to
account for potential energy.  Also, $Q$ should go to zero (and does
in the IMCF and CFM cases) at and inside the horizons of the black
holes since matter inside black holes should not affect the total
mass.

We also believe that the generalized Schoen-Yau identity and the
generalized Jang equation have much potential for many possible
applications in the study of general relativity.  One point of view,
for example, is that the Jang-Zero Divergence equations give a
canonical way of embedding Cauchy data $(M^3,g,k)$ into a static
spacetime.  If one is interested in understanding how the initial
Cauchy data evolves under the vacuum Einstein equations, or some
other equation coupled with the Einstein equation, then one could
compute how the canonical static metrics associated with the
evolving Cauchy data slices evolve.  One nice property of this
approach is that if the initial Cauchy data is a slice of the
Schwarzschild spacetime, and we are solving the vacuum Einstein
equations for example, then while the Cauchy data is evolving in
what might appear to be complicated ways, the associated canonical
static spacetime remains the Schwarzschild spacetime.  Also, since
the generalized Jang equation blows up on horizons, this method
could only be used to study the exterior region of spacetimes
outside the apparent horizons of black holes.  There may be some
advantages to this restriction if this becomes a natural way to
avoid spacetime singularities.

\appendix

\section{Introduction to the Appendices}

The target audience of these appendices are graduate students and
other researchers who are interested in entering geometric
relativity as a field to study.  As such, we have included more
detail in these calculations than is typical.  We justify this
choice in part with the fact that there are so many computations,
many people would have a hard time duplicating all of these
computations in a reasonable amount of time, even with well chosen
hints.  We also hope that these appendices will be useful to
students and researchers who are interested in practicing their
computational skills.  We recommend the book ``Semi-Riemannian
Geometry with Applications to Relativity'' by Barrett O'Neill as an
excellent introduction to the differential geometry of general
relativity, and we mostly follow that book's notation here. Readers
should go through the calculations in these appendices in order
since notational conventions which are established in one appendix
apply to the appendices which follow as well.

The authors would like to thank Alan Parry for helping with the
TeXing and Jeff Jauregui for helpful comments improving the
readability of these appendices.

\section{Curvature of Static Spacetimes}\label{app:curv}

In this section we compute the Einstein curvature, Ricci curvature,
and scalar curvature of the general static spacetime metric
\[\tilde{g} = -\phi(x)^2 dt^2 + \bar{g}\] on $\Real \times M^3$, where
$t \in \Real$, $x \in M^3$ and $\bar{g}$ is a positive definite
metric on $M^3$.

First, choose a coordinate chart on $(M^3,\bar{g})$ with coordinates
$(x^1, x^2, x^3)$ and let $x^0 = t$.  Then $\{\bar\partial_\alpha =
\frac{\partial}{\partial x^\alpha}\}_{\alpha=0}^3$ is a basis of the
tangent plane at each point of the spacetime.  Let
$\{\bar\omega^\alpha\}_{\alpha=0}^3$ be the corresponding dual basis
of one forms at each point of the spacetime so that
$\bar\omega^\alpha(\bar\partial_\beta) = \delta^\alpha_\beta$. (We
use bars over these bases instead of tildes to be consistent with
section \ref{sec:gsyi} and subsequent appendices). Finally, we
define the components of $\tilde{g}$ (which are the same as the
components of $\bar{g}$ for tangent vectors to $M^3$) to be
\[
\tilde{g}_{\alpha\beta} = \langle
{\bar\partial}_\alpha,{\bar\partial}_\beta \rangle_{\tilde{g}}
\;\;\;\;\;( = \tilde{g}({\bar\partial}_\alpha, {\bar\partial}_\beta)
\mbox{ by convention})
\]
so that $\tilde{g}_{00} = -\phi(x)^2$, $\tilde{g}_{0i} =
\tilde{g}_{i0} = 0$, and $\tilde{g}_{ij} = \langle
{\bar\partial}_i,{\bar\partial}_j \rangle_{\bar{g}}$ for $1 \le i,j
\le 3$.

\vspace{.1in} \noindent \textbf{(Notation:  We adopt the convention
that Greek indices always range from $0$ to $3$ and Latin indices
always range from $1$ to $3$.  Also, we adopt Einstein's convention
that any time an index is both an upper and lower index in an
expression, summation over that index is implied.)}\vspace{.1in}

\noindent (Recall also that $\tilde{g}^{\alpha\beta}$ are the
components of the inverse matrix of $\tilde{g}$ expressed as a
matrix at each point of the coordinate chart, and that indices of a
tensor may be raised or lowered by contracting with
$\tilde{g}^{\alpha\beta}$ or $\tilde{g}_{\alpha\beta}$, respectively
\cite{ONEILL}).

 By the Koszul formula \cite{ONEILL}, the Levi-Civita connection
 $\tilde\nabla$ of $\tilde{g}$ can be expressed in terms of its components as
\[\tilde\nabla_{{\bar\partial}_\alpha}{\bar\partial}_\beta =
 {{\tilde\Gamma}_{\alpha \beta}}^{\ \ \gamma} {\bar\partial}_\gamma,\]
where
\begin{equation}\label{eqn:Christoffel}
   {{\tilde\Gamma}_{\alpha\beta}}^{\ \ \gamma} = \frac12 \tilde{g}^{\gamma\theta} \left(
   \tilde{g}_{\alpha\theta,\beta} + \tilde{g}_{\beta\theta,\alpha} - \tilde{g}_{\alpha\beta,\theta}   \right)
\end{equation}
are called the Christoffel symbols of ${\tilde\nabla}$.

\vspace{.1in} \noindent \textbf{(Notation:  Commas denote
differentiation with respect to the coordinate chart so that
$\tilde{g}_{\alpha\beta,\theta} = \frac{\partial
\tilde{g}_{\alpha\beta}}{\partial x^\theta}$). }  \vspace{.1in}

Plugging in our expressions for $\tilde{g}_{\alpha\beta}$, short
calculations reveal that
\begin{equation*}
   0 = {{\tilde\Gamma}_{00}}^{\ \ 0} = {{\tilde\Gamma}_{0i}}^{\ \ j} = {{\tilde\Gamma}_{i0}}^{\ \ j} =
   {{\tilde\Gamma}_{ij}}^{\ \ 0} ,
\end{equation*}
\begin{equation*}
   {{\tilde\Gamma}_{0i}}^{\ \ 0} = {{\tilde\Gamma}_{i0}}^{\ \ 0} = \frac{\phi_i}{\phi},
   \;\;\;\;\;\mbox{ and that }\;\;\;\;\;
   {{\tilde\Gamma}_{00}}^{\ \ i} = \phi\cdot\phi^{\bar\imath},
\end{equation*}
where $\phi_i = \phi_{,i} = \frac{\partial\phi(x)}{\partial x^i}$
and, as previously stated, $\phi^{\bar\imath} =\bar{g}^{ij}\phi_j
 = \tilde{g}^{ij}\phi_j = \tilde{g}^{i\alpha}\phi_\alpha = \phi^{\tilde\imath}$ since $\phi$
is only a function of $x$ and does not depend on $t$.  Note that we
are using the convention that a raised index with a bar or tilde
over it denotes raising the index with $\bar{g}$ or $\tilde{g}$,
respectively.

Since the Lie bracket of coordinate vector fields is zero, it
follows from the definition of the Riemann curvature tensor that
\begin{align*}
  R_{ijk}^{\ \ \ l} &=
  {\bar\omega}^l({\tilde\nabla}_{{\bar\partial}_i}{\tilde\nabla}_{{\bar\partial}_j}{\bar\partial}_k-{\tilde\nabla}_{{\bar\partial}_j}{\tilde\nabla}_{{\bar\partial}_i}{\bar\partial}_k).
\intertext{Hence,}
 R_{ijk}^{\ \ \ l} &= {\bar\omega}^{l}
\left({\tilde\nabla}_{{\bar\partial}_{i}}
\left({\tilde\Gamma}_{jk}^{\ \
\alpha}{\bar\partial}_{\alpha}\right)-{\tilde\nabla}_{{\bar\partial}_{j}}
  \left({\tilde\Gamma}_{ik}^{\ \ \alpha}{\bar\partial}_{\alpha}\right)\right) \\
  &= {\bar\omega}^{l} \left(\left({\tilde\Gamma}_{jk}^{\ \ \alpha}\right)_{,i} {\bar\partial}_{\alpha}+{\tilde\Gamma}_{jk}^{\ \ \alpha}{\tilde\Gamma}_{i\alpha}^{\ \ \, \beta}{\bar\partial}_{\beta}- \left( {\tilde\Gamma}_{ik}^{\ \ \alpha} \right)_{,j} {\bar\partial}_{\alpha}- {\tilde\Gamma}_{ik}^{\ \ \alpha}{\tilde\Gamma}_{j\alpha}^{\ \ \, \beta} \bar\partial_{\beta}\right) \\
  &=\left({\tilde\Gamma}_{jk}^{\ \ l} \right)_{,i}-\left({\tilde\Gamma}_{ik}^{\ \ l} \right)_{,j}+\sum_{\alpha}\left({\tilde\Gamma}_{jk}^{\ \ \alpha}{\tilde\Gamma}_{i\alpha}^{\ \ l}- {\tilde\Gamma}_{ik}^{\ \ \alpha} {\tilde\Gamma}_{j\alpha}^{\ \ \, l}\right)
\end{align*}
For the beginner, we note that about half of the text books define
the Riemann curvature tensor to be the negative of what we used
above.  However, all texts eventually end up with the same
definition of the Ricci curvature (defined in a moment) which is
agreed to be a positive multiple of the metric on the standard
sphere.

Plugging in our formulas for the Christoffel symbols, we thus
compute that
\begin{align*}
  \tilde{R}_{0jk}^{\ \ \ 0} &= \left( {\tilde\Gamma}_{jk}^{\ \ 0} \right)_{,0}- \left( {\tilde\Gamma}_{0k}^{\ \ 0}\right)_{,j}+ \sum_{\alpha} \left({\tilde\Gamma}_{jk}^{\ \ \alpha} {\tilde\Gamma}_{0\alpha}^{\ \ 0}- {\tilde\Gamma}_{0k}^{\ \ \alpha} {\tilde\Gamma}_{j\alpha}^{\ \ 0} \right) \\
  &= \left( -\frac{\phi_k}{\phi} \right)_{,j} + \sum_{m} {\tilde\Gamma}_{jk}^{\ \ m} \cdot \frac{\phi_m}{\phi} - \frac{\phi_k}{\phi} \cdot \frac{\phi_j}{\phi} \\
  &= -\frac{{\overline\Hess}_{jk}\phi}{\phi} \\
\intertext{from which it follow that}
  \tilde{R}_{j00}^{\ \ \ j} &= \phi {\overline\Delta}\phi
\end{align*}
where $\overline\Hess$ is the Hessian and $\overline\Delta$ is the
Laplacian on $(M^3,\bar{g})$.  The second computation follows from
the first by first lowering the raised 0 index (introducing a factor
of $-\phi^2$), using the antisymmetry of the Riemann curvature
tensor to switch indices, and then taking the trace of the Hessian.
We remind the reader that the Latin letters $j$, $k$, $m$ range from
1 to 3. The beginning student should review the definition of the
Hessian, the Laplacian, and the use of normal coordinates.  In this
case, we note that we may choose normal coordinates on $M^3$ such
that ${{\tilde\Gamma}_{jk}}^{\ \ m} = 0$ at a single point.
Similarly, it is straightforward to verify that
\begin{equation*}
\tilde{R}_{0j0}^{\ \ \ 0} = 0 = \tilde{R}_{kj0}^{\ \ \ \, k}.
\end{equation*}

It turns out that the above components of the Riemann curvature
tensor are all that we need to compute the Ricci curvature of the
spacetime. For example,
\begin{equation*}
  \tilde\Ric_{jk}=\tilde{R}_{\alpha jk}^{\ \ \ \, \alpha}=\tilde{R}_{a jk}^{\ \ \ \, a} +  \tilde{R}_{0jk}^{\ \ \ 0}
  = \bar{R}_{a jk}^{\ \ \ \, a} +  \tilde{R}_{0jk}^{\ \ \ 0}
  = \overline{\Ric}_{jk} + \tilde{R}_{0jk}^{\ \ \ 0}.
\end{equation*}
The first equality is the definition of the Ricci curvature as the
trace of the Riemann curvature tensor.  For the second equality
recall our summation convention for Latin and Greek indices stated
above.  The third equality is a consequence of the Gauss equation
for submanifolds since the $t=0$ slice of our spacetime has zero
second fundamental form by symmetry.  The fourth equality is simply
the definition of the Ricci curvature $\overline\Ric$ of
$(M^3,\bar{g})$. Also,
\begin{equation*}
   \tilde\Ric_{00}=\tilde{R}_{\alpha00}^{\ \ \ \, \alpha}=\tilde{R}_{000}^{\ \ \ 0}+R_{j00}^{\ \ \ j}=\tilde{R}_{j00}^{\ \ \ j}
\end{equation*}
by antisymmetry of the Riemann curvature tensor.  Finally,
\begin{equation*}
   \tilde\Ric_{j0}=\tilde{R}_{\alpha j0}^{\ \ \ \, \alpha}=\tilde{R}_{0j0}^{\ \ \ 0}+\tilde{R}_{kj0}^{\ \ \ \, k}=0.
\end{equation*}

Thus, putting it all together, we have formulas for the components
of the Ricci curvature of the static spacetime metric $\tilde{g}=
-\phi(x)^2 dt^2 + \bar{g}$ on $\Real \times M^3$  , namely
\begin{align*}
  \tilde\Ric_{00} &= \phi \overline{\Delta}\phi \\
  \tilde\Ric_{jk} &= \overline{\Ric}_{jk}-\frac{\overline{\Hess}_{jk}\phi}{\phi} \\
  \tilde\Ric_{j0} &= \tilde\Ric_{0j}=0
\end{align*}
in terms of the Ricci curvature $\overline\Ric$ of $(M^3,\bar{g})$
and the Hessian and Laplacian of $\phi$ on $(M^3,\bar{g})$.

Next we can compute the scalar curvature of the spacetime by taking
the trace of the Ricci curvature,
\begin{equation*}
   \tilde{R} = \tilde{g}^{jk} \tilde\Ric_{jk}
   = \bar{R} - 2 \frac{\overline\Delta\phi}{\phi}.
\end{equation*}
Hence, since $\tilde{G} = \tilde{\Ric} - \frac12 \tilde{R}
\tilde{g}$, the components of the Einstein curvature tensor of the
static spacetime metric are
\begin{align*}
  \tilde{G}_{00} &= \frac12 \bar{R}\phi^2  \\
  \tilde{G}_{jk} &= \overline{\Ric}_{jk}-\frac{\overline{\Hess}_{jk}\phi}{\phi} +
  \left(\frac{\overline\Delta \phi}{\phi} - \frac{\bar{R}}{2}\right) \bar{g} \\
  \tilde{G}_{j0} &= \tilde\Ric_{0j}=0
\end{align*}
as desired.

\section{The Second Fundamental Form of the
Graph}\label{app:secondfundamentalform}

In this section we will compute the second fundamental form of a
space-like slice of the static spacetime $\left(\Real \times M^3,
\tilde{g} \right)$, where

\begin{equation} \label{staticspacetime2}
   \tilde{g} = -\phi^2 dt^2 + \bar{g},
\end{equation}
$\phi$ is a real-valued function on $M$,
 and $\bar{g}$ is a
Riemannian metric on $M$. Given a real-valued function $f$ on $M$,
define the graph map
\begin{equation}\label{graphmap2}
   F: M \mapsto \Real \times M
\end{equation}
where $F(x) = (f(x), x)$.

As we established in section \ref{sec:gsyi}, we will abuse
terminology slightly and always refer to the image of the graph map
$F(M)$ simply as $M$ and the $t=0$ slice of the constructed
spacetime as $\bar{M}$.  This notation is convenient since then
$(M,g)$ and $(\bar{M},\bar{g})$ are space-like hypersurfaces of the
spacetime $\left(\Real \times M^3, -\phi^2 dt^2 + \bar{g} \right)$
(given appropriate bounds on the gradient of $f$). Let $\pi: M
\mapsto \bar{M}$ be the projection map $\pi(f(x),x) = (0,x)$ to the
$t=0$ slice of the spacetime.

\begin{figure}[htbp]
\begin{center}

\input{bray_figure1.pstex_t}

\caption{Schematic diagram of the constructed static spacetime}
\label{figure2}
\end{center}
\end{figure}

We repeat some notation and definitions from section \ref{sec:gsyi}
for clarity.  Let $\bar{\partial}_0 = \partial_t$ and
$\{\bar{\partial}_i\}$ be coordinate vectors tangent to $\bar{M}$.
Define
\begin{equation}
   \partial_i = \bar\partial_i + f_i\bar\partial_0
\end{equation}
to be the corresponding coordinate vectors tangent to $M$ so that
$\pi_*(\partial_i) = \bar\partial_i$.  Then in this coordinate
chart, the components of the metrics $g$ and $\bar{g}$ induced from
the spacetime are
\[
g_{ij} = \langle \partial_i, \partial_j \rangle \;\;\mbox{   and
}\;\; \bar{g}_{ij} = \langle \bar\partial_i, \bar\partial_j \rangle
\]
where the angle brackets refer to the spacetime metric.  Then it
follows immediately that
\begin{equation}\label{eqn:metric2}
   g_{ij} = \bar{g}_{ij} - \phi^2 f_i f_j.
\end{equation}
We comment here that the reader should think of $g$, $\phi$, and
$f$ as the variables that we get to choose which determine
$\bar{g}$. The metric $g$ comes from the initial Cauchy data
$(M,g,k)$ and $\phi$ and $f$ are functions which will satisfy a
system of equations of our choosing.

The inverse of $\{g_{ij}\}$ turns out to be
\begin{equation}\label{eqn:inversemetric2}
   g^{ij} = \bar{g}^{ij} + v^i v^j ,
\end{equation}
where
\begin{equation}\label{vdef2}
   v^i = \frac{\phi f^{\bar{i}}}{(1-\phi^2|df|^2_{\bar{g}})^{1/2}} = \frac{\phi
   f^i}{(1+\phi^2|df|^2_g)^{1/2}} .
\end{equation}
The above computation is most easily verified at each point in
normal coordinates of $\bar{g}$ at that point, where the gradient of
$f$ is assumed to lie in the first coordinate direction.  The second
part of equation \ref{vdef2} can be computed in the same manner, but
where we consider that $\bar{g}_{ij} = {g}_{ij} + \phi^2 f_i f_j$
and then use normal coordinates as before, but this time for the
metric $g$.  (For the beginner, the use of normal coordinates is
exemplified in more detail in a moment.)

We also define
\begin{equation}\label{vvecdef2}
   \bar{v} = v^i \bar\partial_i \;\;\;\;\;\;\mbox{ and }\;\;\;\;\;\; v = v^i \partial_i
\end{equation}
so that $\pi_*(v) = \bar{v}$, and observe the useful identity
\begin{equation}\label{eqn:usefulidentity2}
(1-\phi^2|df|^2_{\bar{g}}) \cdot (1+\phi^2|df|^2_g) = 1,
\end{equation}
which follows directly from computing the ratio of the volume forms
of $g$ and $\bar{g}$ two different ways, namely with respect to $g$
and then $\bar{g}$.

As established in section \ref{sec:gsyi}, we use the convention that
a barred index (as in $f^{\bar{i}}$ above) denotes an index raised
(or lowered) by $\bar{g}$ as opposed to $g$.  That is, $f^{\bar{i}}
= \bar{g}^{ij}f_j$, where as usual $f_j = \partial f / \partial x^j$
in the coordinate chart. In general, barred quantities will be
associated with the $t=0$ slice $(\bar{M},\bar{g})$ and unbarred
quantities will be associated with the graph slice $(M,g)$.

Our next step is to compute the unit normal vector ${n}$ to the
graph slice $(M,g)$ defined by $f$.  It is straightforward to verify
that
\begin{equation}
  {n} = \frac{\bar{\partial}_{0} + \phi^{2}f^{\bar{k}}\bar{\partial}_{k}}{\phi \left(1-\phi^{2}\left\vert df \right\vert^{2}_{\bar{g}} \right)^{1/2}}
\end{equation}
has unit length in the spacetime metric, is perpendicular to the
tangent vectors $\partial_i = \bar\partial_i + f_i\bar\partial_0$ to
the graph slice, and hence must be the correct expression.

Following our convention for the definition of the second
fundamental form defined in equation \ref{eqn:sffdef}, we thus have
that the components of the second fundamental form $h$ are
\begin{align*}
h_{ij} &= h(\partial_i, \partial_j) \\
&= - \langle \tilde\nabla_{\partial_i} \partial_j ,
{n}\rangle \\
&= \langle \tilde\nabla_{\partial_i} {n} ,
\partial_j \rangle \\
&= \left\langle
\tilde\nabla_{(\bar{\partial}_{i}+f_{i}\bar{\partial}_{0})} \left[
\frac{\bar{\partial}_{0} + \phi^{2}f^{\bar{k}}
\bar{\partial}_{k}}{\phi \left( 1- \phi^{2}\left\vert df
\right\vert^{2}_{\bar{g}} \right)^{1/2}} \right], \
\bar{\partial}_{j} + f_{j}\bar{\partial}_{0} \right\rangle \\
&=  \frac{\left\langle
\tilde\nabla_{(\bar{\partial}_{i}+f_{i}\bar{\partial}_{0})} \left[
\bar{\partial}_{0} + \phi^{2}f^{\bar{k}} \bar{\partial}_{k} \right],
\ \bar{\partial}_{j} + f_{j}\bar{\partial}_{0} \right\rangle}
{\phi\left( 1- \phi^{2}\left\vert df \right\vert^{2}_{\bar{g}}
\right)^{1/2}}
\end{align*}
where $\tilde\nabla$ is the Levi-Civita connection on our spacetime
$\left(\Real \times M^3, \tilde{g} \right)$. The third and fifth
equalities follow from the fact that $\langle {n},
\partial_j \rangle = 0$ on $M$.

From the form of the above expression, we see that the Christoffel
symbols of the spacetime, defined and computed in appendix
\ref{app:curv}, are going to come into play.  From those
computations, it follows that
\begin{align*}
  \tilde\nabla_{\bar{\partial}_{0}} \bar{\partial}_{0} &= \tilde{\Gamma}_{00}^{\ \ 0} \bar{\partial}_{0} + \tilde{\Gamma}_{00}^{\ \ k} \bar{\partial}_{k} = \phi \phi^{\bar{k}}\bar{\partial}_{k} \\
  \tilde\nabla_{\bar{\partial}_{0}} \bar{\partial}_{i} &= \tilde{\Gamma}_{0i}^{\ \ 0} \bar{\partial}_{0} + \tilde{\Gamma}_{0i}^{\ \ k} \bar{\partial}_{k} = \frac{\phi_{i}}{\phi}\bar{\partial}_{0} \\
  \tilde\nabla_{\bar{\partial}_{i}} \bar{\partial}_{0} &= \tilde{\Gamma}_{i0}^{\ \ 0} \bar{\partial}_{0} + \tilde{\Gamma}_{i0}^{\ \ k} \bar{\partial}_{k} = \frac{\phi_{i}}{\phi} \bar{\partial}_{0}
\end{align*}
where we remind the reader that Latin indices, when summation is
implied by one raised and one lowered, only sum from 1 to 3, and a
bar over a raised index indicates that the index was raised with
$\bar{g}$ as opposed to $g$.

It now becomes convenient to use normal coordinates on
$(\bar{M},\bar{g})$.  Note that these are not normal coordinates on
the whole spacetime, just on the $t=0$ slice of the spacetime
$(\bar{M},\bar{g})$.  Since by symmetry this slice has zero second
fundamental form, $\tilde\nabla_{\bar\partial_i} \bar\partial_k
=\overline\nabla_{\bar\partial_i} \bar\partial_k = 0$ at a single
point of our choosing.  In normal coordinates, derivatives of the
metric components $\bar{g}_{ij}$ and $\bar{g}^{ij}$ are zero at the
chosen point, thereby making the Christoffel symbols zero at that
point as well.  In addition, $(f^{\bar{k}})_i = (f_m \bar{g}^{mk})_i
= f_{im} \bar{g}^{mk} = f_i^{\ \bar{k}}$ at this single arbitrary
point. Note that $f_{im} = \frac{\partial^2 f}{\partial x^i \partial
x^m}$ is simply a coordinate chart second derivative in our
notation.

Hence,
\begin{align*}
\tilde\nabla_{(\bar{\partial}_{i}+f_{i}\bar{\partial}_{0})} & \left[
\bar{\partial}_{0} + \phi^{2}f^{\bar{k}} \bar{\partial}_{k} \right]
\\
 &= \frac{\phi_i}{\phi}\bar\partial_0
   + (2 \phi \phi_i f^{\bar{k}} + \phi^2 f_i^{\
   \bar{k}})\bar\partial_k
 + f_i (\phi \phi^{\bar{k}}\bar{\partial}_{k}
       + \phi^2 f^{\bar{k}} \cdot \frac{\phi_k}{\phi} \bar\partial_0)\\
&= \left( \frac{\phi_i}{\phi} + \phi f_i f^{\bar{k}} \phi_k
\right)\bar\partial_0
 + \left(2 \phi \phi_i f^{\bar{k}} + \phi^2 f_i^{\
   \bar{k}} + \phi f_i \phi^{\bar{k}} \right) \bar\partial_k.
\end{align*}
Since $\langle \bar\partial_0, \bar\partial_0 \rangle = -\phi^2$ and
$\langle \bar\partial_k, \bar\partial_j \rangle = \bar{g}_{kj}$
(which then lowers indices on other terms), we have computed that
\begin{align}\notag
h_{ij} &=  \frac{ -\phi^2 f_j \left( \frac{\phi_i}{\phi} + \phi f_i
f^{\bar{k}} \phi_k \right)
 + \left(2 \phi \phi_i f^{\bar{k}} + \phi^2 f_i^{\
   \bar{k}} + \phi f_i \phi^{\bar{k}} \right) \bar{g}_{kj}
  } {\phi\left( 1- \phi^{2}\left\vert df
\right\vert^{2}_{\bar{g}} \right)^{1/2}} \\ \notag
 &=  \frac{  -\phi_i f_j - \phi^2 f_i f_j f^{\bar{k}} \phi_k
 + 2 \phi_i f_j + \phi f_{ij} + f_i \phi_j
  } {\left( 1- \phi^{2}\left\vert df
\right\vert^{2}_{\bar{g}} \right)^{1/2}} \\
&=  \frac{ \phi f_{ij}
  + f_i \phi_j +  \phi_i f_j - \phi^2 \langle df, d\phi\rangle_{\bar{g}} f_i f_j
  } {\left( 1- \phi^{2}\left\vert df
\right\vert^{2}_{\bar{g}} \right)^{1/2}}
\label{eqn:secfundformnormal}
\end{align}
at the chosen point in normal coordinates for $(\bar{M},\bar{g})$.
Since $\overline\Hess_{ij} f = f_{ij}$ as well at the chosen point
in these normal coordinates, we have that
\begin{align}\label{eqn:tensorcomp}
h_{ij} &=  \frac{ \phi \overline\Hess_{ij} f
  + f_i \phi_j +  \phi_i f_j - \phi^2 \langle df, d\phi\rangle_{\bar{g}} f_i f_j
  } {\left( 1- \phi^{2}\left\vert df
\right\vert^{2}_{\bar{g}} \right)^{1/2}}
\end{align}
at the chosen point.  However, the above equation represents the
components of the tensorial equation
\begin{align}\label{eqn:tensor}
h &=  \frac{ \phi \overline\Hess f
  + df \otimes d\phi +  d\phi \otimes df - \phi^2 \langle df, d\phi\rangle_{\bar{g}} df \otimes df
  } {\left( 1- \phi^{2}\left\vert df
\right\vert^{2}_{\bar{g}} \right)^{1/2}}
\end{align}
by which we mean both the left and right sides of the equation are
tensors.  Since tensorial equations may be verified in any
coordinate chart, we conclude that equation \ref{eqn:tensor} is true
at our chosen point.  Since our chosen point was arbitrary, it
follows that equation \ref{eqn:tensor} is true at every point. Thus,
equation \ref{eqn:tensorcomp} is true at every point as well, in any
coordinate chart.

The beginner differential geometer should take note on how using
normal coordinates simplified these computations substantially.
However, the computation can also be done straightforwardly without
using normal coordinates, just not as elegantly.

As we mentioned originally in section \ref{sec:gsyi}, we extend $h$
trivially in our constructed static spacetime so that $h(\partial_t,
\cdot) = 0$.  Note that this gives $h(\partial_i,\partial_j) =
h(\bar\partial_i, \bar\partial_j)$, so we can call this term
$h_{ij}$ without ambiguity.  Our next goal is to convert the above
formula for $h_{ij}$ expressed with respect to $\bar{g}$ to one
expressed with respect to $g$.

To convert the tensor $\overline\Hess f$, we must recall that it is
defined to be the covariant derivative $\overline\nabla$ of the
1-tensor $df$.  Note that $df$ does not involve any metric, since
$df(W) = W(f)$ by definition.  However, $\overline\nabla$ does
involve the metric $\bar{g}$ when applied to tensors. For example,
in coordinates,
\begin{align*}
\overline\Hess_{ij} f =&
\overline\nabla(df)(\bar\partial_i,\bar\partial_j) \\ =&
\bar\partial_i (df(\bar\partial_j)) -
df(\overline\nabla_{\bar\partial_i} \bar\partial_j) \\ =&
\bar\partial_i(\bar\partial_j(f)) -
(\overline\nabla_{\bar\partial_i} \bar\partial_j)(f) \\ =& f_{ij} -
\bar\Gamma_{ij}^{\ \ k} f_k
\end{align*}
which involves the metric $\bar{g}$ and its first derivatives by
equation \ref{eqn:Christoffel}.  Hence,
\begin{align}\label{HessDiff}
\overline\Hess_{ij} f - \Hess_{ij} f = & - (\bar\Gamma_{ij}^{\ \ k}
- \Gamma_{ij}^{\ \ k}) f_k.
\end{align}
Thus, we need to compute the difference of the Christoffel symbols
of $\bar{g}$ and $g$.  By equations \ref{eqn:metric2} and
\ref{eqn:inversemetric2},
\begin{align}
  \notag 2\Gamma_{ij}^{\ \ k} =&  g^{km}(g_{im,j}+g_{jm,i}-g_{ij,m}) \\
  \notag =&  \left(\bar{g}^{km} + \frac{\phi^{2}f^{\bar{k}}f^{\bar{m}}}{1-\phi^{2} \left\vert df \right\vert_{g}^{2}}
  \right)\cdot \\
  \notag &\quad \left(\bar{g}_{im,j} + \bar{g}_{jm,i} - \bar{g}_{ij,m} - (\phi^{2}f_{i}f_{m})_{,j} - (\phi^{2}f_{j}f_{m})_{,i} + (\phi^{2}f_{i}f_{j})_{,m} \right) \\
  \notag =& \left(\bar{g}^{km} + v^{k}v^{m} \right) \left\{\bar{g}_{im,j} + \bar{g}_{jm,i} - \bar{g}_{ij,m} - 2\phi\phi_{j}f_{i}f_{m} - \phi^{2}f_{m}f_{ij} - \phi^{2}f_{i}f_{mj} \right. \\
  \notag &\quad \left. -2 \phi \phi_{i}f_{j}f_{m} - \phi^{2}f_{m}f_{ji} - \phi^{2}f_{j}f_{mi} + 2\phi\phi_{m}f_{i}f_{j} + \phi^{2}f_{j}f_{im} + \phi^{2}f_{i}f_{jm} \right\} \\
  \notag =& \left(\bar{g}^{km} + v^{k}v^{m} \right)
   \left\{\bar{g}_{im,j} + \bar{g}_{jm,i} - \bar{g}_{ij,m} - 2\phi^{2}f_{m}f_{ij} \right. \\
   &\quad \left.
   - 2\phi\phi_{j}f_{i}f_{m}  - 2 \phi \phi_{i}f_{j}f_{m} + 2\phi\phi_{m}f_{i}f_{j} \right\}
\end{align}
so that in normal coordinates on $(\bar{M},\bar{g})$,
\begin{align}
  \overline\Gamma_{ij}^{\ \ k} - {\Gamma}_{ij}^{\ \ k} =& \phi^{2}f_{ij}f^{\bar{k}}
  + \phi f_{i} \phi_{j}f^{\bar{k}}+ \phi\phi_{i}f_{j}f^{\bar{k}}-\phi f_{i}f_{j}\phi^{\bar{k}}
   \notag \\
  &\quad + v^{k}\left[
  \phi^{2} v(f) f_{ij}+ \phi  v(f) \phi_{j}f_{i}
  + \phi  v(f) \phi_{i}f_{j}-\phi  v(\phi) f_{i}f_{j}
  \right]\notag \\
  \notag =& \phi f^{\bar{k}} \left(\phi f_{ij} + f_{i} \phi_{j} + \phi_{i}f_{j} \right) - \phi f_{i} f_{j}\phi^{\bar{k}} \notag \\
  &\quad + v^{k} \phi v(f) \left(\phi f_{ij} + f_{i}\phi_{j} + \phi_{i}f_{j} \right) - \phi  v(\phi) f_{i}f_{j} v^{k} \notag \\
  =& \left(\phi f^{\bar{k}} + \frac{\phi^{3}\left\vert df \right\vert_{\bar{g}}^{2} f^{\bar{k}}}{1-\phi^{2}\left\vert df \right\vert_{\bar{g}}^{2}} \right)
  \left(\phi f_{ij} + f_{i} \phi_{j} + \phi_{i} f_{j} \right) \notag \\
  &\quad - \frac{ \phi^{3} \langle df, d\phi \rangle_{\bar{g}} f_{i} f_{j} f^{\bar{k}}}{1-\phi^{2} \left\vert df \right\vert_{\bar{g}}^{2}} - \phi
  f_{i}f_{j}\phi^{\bar{k}}\notag \\
  =& \frac{\phi f^{\bar{k}}}{1-\phi^{2}\left\vert df \right\vert_{\bar{g}}^{2}} \left[ \phi f_{ij} + f_{i} \phi_{j} + \phi_{i}f_{j} - \phi^{2} \langle df, d\phi \rangle f_{i}f_{j} \right] - \phi
  f_{i}f_{j}\phi^{\bar{k}} \notag \\
  =& h_{ij}v^{k} - \phi f_{i}f_{j} \phi^{\bar{k}}.
  \label{eqn:proofofi3}
\end{align}
(Note that the above equation is a proof of identity \ref{i3} from
section \ref{sec:gsyi}.)  Plugging this equation into equation
\ref{HessDiff} gives us
\begin{align}\label{HessDiff2}
\overline\Hess_{ij} f - \phi f_{i}f_{j} \langle
d\phi,df\rangle_{\bar{g}}&= \Hess_{ij} f  - h_{ij}v(f) .
\end{align}

Hence, using the above formula and equation
\ref{eqn:usefulidentity2}, we can transform equation
\ref{eqn:tensorcomp} to
\begin{align*}
h_{ij} =&  \frac{ \phi \overline\Hess_{ij} f
  + f_i \phi_j +  \phi_i f_j - \phi^2 \langle df, d\phi\rangle_{\bar{g}} f_i f_j
  } {\left( 1- \phi^{2}\left\vert df
\right\vert^{2}_{\bar{g}} \right)^{1/2}}  \\
= &  {\left( 1 + \phi^{2}\left\vert df \right\vert^{2}_{{g}}
\right)^{1/2}} \left( \phi \Hess_{ij} f - \phi h_{ij}v(f)
  + f_i \phi_j +  \phi_i f_j
  \right)
\end{align*}
so that by the definition of $v$ in equation \ref{vdef2}
\begin{align*}
\left( 1 + \phi^{2}\left\vert df \right\vert^{2}_{{g}} \right)h_{ij}
=& {\left( 1 + \phi^{2}\left\vert df \right\vert^{2}_{{g}}
\right)^{1/2}} \left( \phi \Hess_{ij} f
  + f_i \phi_j +  \phi_i f_j
  \right).
\end{align*}
Thus, we see that the second fundamental form of the graph slice
expressed with respect to the metric $g$ is
\begin{eqnarray}\label{eqn:tensorcomponents2}
   h_{ij} &=& \frac{\phi \mbox{Hess}_{ij}f + (f_i \phi_j + \phi_i f_j)}
              { (1+\phi^2|df|^2_g)^{1/2}},
\end{eqnarray}
as claimed in section \ref{sec:gsyi}.  Technically, we've shown that
the above equation is true at an arbitrary point in normal
coordinates.  Again, the above equation represents the components of
the tensorial equation
\begin{eqnarray}
   h &=& \frac{\phi \mbox{Hess}f + (df \otimes d\phi + d\phi \otimes df)}
              { (1+\phi^2|df|^2_g)^{1/2}},
\end{eqnarray}
which is therefore true at the arbitrary point, and thus is true
everywhere.  Hence, equation \ref{eqn:tensorcomponents2} is true at
every point as well, in any coordinate chart.

\section{Derivation of Identities}\label{app:identities}

In this appendix we finish the proof of the generalized Schoen-Yau
identity sketched out in section \ref{sec:gsyi}.  We begin with
equation \ref{eqn:startofappendix} derived in that section, which we
repeat for clarity:
\begin{eqnarray*}
\bar{R} &=& 16\pi(\mu - J(v))
        + (\mbox{tr}_g h)^2 -  (\mbox{tr}_g k)^2
        - \|h\|_g^2 + \|k\|_g^2   \notag \\ &&
        + 2 \, v(\mbox{tr}_g h) - 2 \, v(\mbox{tr}_g k)
        - 2 \, \mbox{div}(h)(v) + 2 \, \mbox{div}(k)(v).
\end{eqnarray*}

We will convert this formula for $\bar{R}$ to an expression in terms
of the $\bar{g}$ metric.  To perform the conversion, we need several
identities (originally listed in section \ref{sec:gsyi}) for
arbitrary symmetric 2-tensors $k$ which we now prove.  We continue
with the notation established in the previous appendix, which might
be thought of as an introduction to this appendix.

\vspace{.2in}\noindent {\bf Identity \ref{i1}}
\[
   (\mbox{tr}_g(k))^2 - \|k\|_g^2 = (\mbox{tr}_{\bar{g}}k)^2 -
   \|k\|_{\bar{g}}^2 + 2 k(\bar{v},\bar{v}) \mbox{tr}_{\bar{g}}k -
   2|k(\bar{v},\cdot)|_{\bar{g}}^2
\]
\noindent {\it Proof:  }
\begin{align*}
  (\tr_{g} k)^{2} - \left\| k \right\|_{g}^{2} &= (g^{ij} k_{ij})^{2} - g^{ik}g^{jl}k_{ij}k_{kl} \\
  \notag &= \left[\left( \bar{g}^{ij} + v^i v^j \right)k_{ij} \right]^{2}
  - \left(\bar{g}^{ik} + v^i v^k \right) \left( \bar{g}^{jl} + v^j v^l \right) k_{ij}k_{kl} \\
  &= \left[\tr_{\bar{g}} k + k({\bar{v}},{\bar{v}}) \right]^{2} - \left\| k \right\|_{\bar{g}}^{2} - k({\bar{v}},{\bar{v}})^{2} - 2\left\vert k({\bar{v}},\cdot) \right\vert_{\bar{g}}^{2} \\
  &= (\tr_{\bar{g}} k)^{2} - \left\| k \right\|_{\bar{g}}^{2} + 2k({\bar{v}},{\bar{v}})\tr_{\bar{g}}k - 2\left\vert k({\bar{v}},\cdot) \right\vert_{\bar{g}}^{2}.
\end{align*}
The first equality is true by definition.  The second equality
follows from equation \ref{eqn:inversemetric2}.  For the third
equality, remember that $k$ is defined to have zero time-time
components and time-spatial components.  Hence, $k(v,w) =
k(\bar{v},w)$ for all $w$ since $v$ projects to $\bar{v}$ (and
consequently $v$ and $\bar{v}$ are equal except for their time
components).

\vspace{.2in}\noindent {\bf Identity \ref{i2}}
\[
   v(\mbox{tr}_g k) = \bar{v}(\mbox{tr}_{\bar{g}}k + k(\bar{v},\bar{v}))
\]
\noindent {\it Proof:  }
\begin{align*}
  v(\tr_{g}k) &= v\left(\tr_{\bar{g}}k + k({\bar{v}},{\bar{v}})\right) \\
  &= {\bar{v}}\left(\tr_{\bar{g}}k + k({\bar{v}},{\bar{v}}) \right).
\end{align*}
The first equality was shown in the proof of identity \ref{i1}.  The
second equality follows since $v$ and $\bar{v}$ only differ by their
time components, and the function being differentiated does not
depend on the time coordinate, by definition.

\vspace{.2in}\noindent {\bf Identity \ref{i3}}
\[
   \overline\Gamma_{ij}^k - \Gamma_{ij}^k = h_{ij}v^k - \phi f_i f_j
   \phi^{\bar{k}}
\]
\noindent {\it Proof:  }  See equation \ref{eqn:proofofi3}.

\vspace{.2in}\noindent {\bf Identity \ref{i4}}
\begin{eqnarray*}
   \mbox{div}(k)(v) &=& \overline{\mbox{div}}(k)(\bar{v})
   + (\overline\nabla_{\bar{v}}k)(\bar{v},\bar{v})
   - 2 |\bar{v}|_{\bar{g}}^2\; k\left(\bar{v},\frac{\overline\nabla \phi}{\phi}\right)
   \\ &&
   + \langle h(\bar{v},\cdot) , k(\bar{v},\cdot) \rangle_{\bar{g}}
   + 2 h(\bar{v},\bar{v}) k(\bar{v},\bar{v})
   + (\mbox{tr}_{\bar{g}}h) k(\bar{v},\bar{v})
\end{eqnarray*}
\noindent {\it Proof:  }
\begin{align*}
  \mbox{div}(k)(v) &=  g^{ij} (\nabla_{\partial_{i}}k)(\partial_{j},v) \\
  &=  g^{ij} \left[\partial_{i} (k(\partial_{j},v)) - k(\nabla_{\partial_{i}}\partial_{j},v) - k(\partial_{j},\nabla_{\partial_{i}} v) \right] \\
  &= g^{ij} \left[(k_{j\alpha} v^{\alpha})_{,i} - \Gamma_{ij}^{\ \ m} k_{m\alpha} v^{\alpha} - k\left(\partial_{j},\  (v^{\alpha})_{,i}\ \partial_{\alpha} + v^{\alpha} \nabla_{\partial_{i}}\partial_{\alpha}\right) \right] \\
  &= g^{ij}\left[ (k_{j\alpha} v^{\alpha})_{,i} - \Gamma_{ij}^{\ \ m} k_{m\alpha} v^{\alpha} - (v^{\alpha})_{,i}\ k_{j\alpha} - v^{\alpha} \Gamma_{i\alpha}^{\ \ \, m} k_{jm} \right] \\
  &= g^{ij} \left(k_{j\alpha,i} - \Gamma_{ij}^{\ \ m}k_{m\alpha} - \Gamma_{i\alpha}^{\ \ \, m} k_{jm} \right)
  v^{\alpha} \\
  &= \left(\bar{g}^{ij} + v^i v^j \right) \left[ k_{j\alpha,i} + k_{m\alpha}\left(-\overline{\Gamma}_{ij}^{\ \ m} + h_{ij}v^{m} - \phi f_{i}f_{j} \phi^{\bar{m}}\right) \right. \\
  &\quad\hspace{1.12in} \left. + k_{jm} \left( -\overline{\Gamma}_{i\alpha}^{\ \ \, m} + h_{i\alpha}v^{m} - \phi f_{i}f_{\alpha} \phi^{\bar{m}} \right) \right] v^{\alpha} \\
  \notag &= \overline{\mbox{div}}(k)({\bar{v}}) + (\tr_{\bar{g}} h)k({\bar{v}},{\bar{v}}) - \phi \left\vert df \right\vert_{\bar{g}}^{2} k(\overline{\nabla}\phi,\ {\bar{v}}) + \langle h({\bar{v}}, \cdot), k({\bar{v}},\cdot) \rangle_{\bar{g}} \\
  \notag &\quad - \phi\  \bar{v}(f) \ k(\overline{\nabla}f, \overline{\nabla}\phi) + (\overline{\nabla}_{{\bar{v}}} k)({\bar{v}},{\bar{v}}) +
  h({\bar{v}},{\bar{v}})k({\bar{v}},{\bar{v}})\\
  &\quad - \phi \ \bar{v}(f)^{2} \ k(\overline{\nabla}\phi,\ {\bar{v}})
  + h({\bar{v}},{\bar{v}})k({\bar{v}},{\bar{v}}) - \phi \ \bar{v}(f)^{2} \ k({\bar{v}}, \overline{\nabla}\phi) \\
  \notag &= \overline{\mbox{div}}(k)({\bar{v}}) + (\tr_{\bar{g}}h)k({\bar{v}},{\bar{v}}) + \langle h({\bar{v}},\cdot), k({\bar{v}},\cdot) \rangle_{\bar{g}} + (\overline{\nabla}_{{\bar{v}}} k)({\bar{v}},{\bar{v}}) \\
  &\quad + 2h({\bar{v}},{\bar{v}})k({\bar{v}},{\bar{v}}) - k\left({\bar{v}},\frac{\overline{\nabla}\phi}{\phi}\right)\left\{2 \phi^{2}\left\vert df \right\vert_{\bar{g}}^{2} \left\vert {\bar{v}} \right\vert_{\bar{g}}^{2} + 2 \phi^{2} \left\vert df \right\vert_{\bar{g}}^{2}
  \right\} \\
  \notag &= (\overline{\nabla}\cdot k)({\bar{v}}) + (\tr_{\bar{g}} h) k({\bar{v}},{\bar{v}}) + \langle h({\bar{v}},\cdot), k({\bar{v}}, \cdot) \rangle_{\bar{g}} + (\overline{\nabla}_{{\bar{v}}} k)({\bar{v}},{\bar{v}}) \\
  &\quad +2h({\bar{v}},{\bar{v}})k({\bar{v}},{\bar{v}}) - 2 \left\vert {\bar{v}} \right\vert_{\bar{g}}^{2}k\left({\bar{v}}, \frac{\overline{\nabla}\phi}{\phi}
  \right).
\end{align*}
The first equality is the definition of divergence.  The second
equality is the definition of the covariant derivative of a tensor.
The third and fourth equalities use Christoffel symbols as defined
in the previous appendix.  The sixth equality uses equation
\ref{eqn:inversemetric2} and then identity \ref{i3}.  The seventh
equality is most easily seen by using normal coordinates with
respect to $\bar{g}$.  The eighth equality combines terms using the
fact that $\bar{v}$ is parallel to $\overline\nabla f$ in
$(M^3,\bar{g})$ by the definition of $\bar{v}$ in equations
\ref{vdef2} and \ref{vvecdef2}.  The ninth equality is the
simplification
\begin{equation*}
2 \phi^{2}\left\vert df \right\vert_{\bar{g}}^{2} \left\vert
{\bar{v}} \right\vert_{\bar{g}}^{2} + 2 \phi^{2} \left\vert df
\right\vert_{\bar{g}}^{2} =
2 \phi^{2}\left\vert df \right\vert_{\bar{g}}^{2} \left( \left\vert
{\bar{v}} \right\vert_{\bar{g}}^{2} + 1 \right)
= \dfrac{2\phi^{2} \left\vert df \right\vert_{\bar{g}}^{2}}{1-
\phi^{2}\left\vert df \right\vert_{\bar{g}}^{2}} = 2 \left\vert
{\bar{v}} \right\vert^{2}.
\end{equation*}

\vspace{.2in}\noindent {\bf Identity \ref{i5}}
\[
   v_{\bar{\imath}\bar{;}j} = h_{ij} + v_{\bar{\imath}} h(\bar{v},\cdot)_j
                 - \frac{\phi_i v_{\bar{\jmath}}}{\phi}
\]
\noindent {\it Proof:  }  First, let us clarify our notation. Recall
that bars refer to the metric $\bar{g}$.  Hence, for example,
$v_{\bar{\imath}} = \bar{g}_{ik}v^k = \langle v, \bar\partial_i
\rangle_{\bar{g}}$, where $v^k$ is defined in equation \ref{vdef2}.
As is standard, semicolons refer to covariant differentiation
(whereas commas refer to coordinate chart derivatives).  Of course
in our case, we need to specify with respect to which metric are we
performing covariant differentiation.  Hence, we place a bar over
the semicolon to denote covariant differentiation with respect to
$\bar{g}$.  Hence, $v_{\bar{\imath}\bar{;}j} = \langle
\overline\nabla_{\bar\partial_j} v, \bar\partial_i
\rangle_{\bar{g}}$.

All of our computations in the proof of this identity and the two
that follow only involve the metric $\bar{g}$, so it is notationally
convenient (though not really necessary) to use normal coordinates
with respect to this metric. Then at that point,
\begin{align*}
  v_{\bar{\imath}\bar{;}j} &= v_{\bar{\imath},j}
  = \left(\frac{\phi f_{i}}{\left(1-\phi^{2}\left\vert df \right\vert_{\bar{g}}^{2}\right)^{1/2}} \right)_{,j} \\
  &= \frac{\phi f_{ij} + f_{i} \phi_{j}}{(1- \phi^{2}\left\vert df \right\vert_{\bar{g}}^{2})^{1/2}} + \frac{\left(\phi \phi_{j}\left\vert df \right\vert_{\bar{g}}^{2} + \phi^{2}f_{\alpha j} f^{\bar{\alpha}} \right)(\phi f_{i})}{(1-\phi^{2}\left\vert df \right\vert_{\bar{g}}^{2})^{3/2}} \\
  &= h_{ij} + \frac{\phi^{2}\langle df, d\phi \rangle_{\bar{g}}f_{i}f_{j} - \phi_{i}f_{j}}{(1- \phi^{2}\left\vert df \right\vert_{\bar{g}}^{2})^{1/2}} + \left[\frac{\phi_{j}}{\phi}\left\vert \bar{v} \right\vert_{\bar{g}}^{2} + \frac{\phi f_{\alpha j}v^{\alpha}}{\left( 1- \phi^{2} \left\vert df \right\vert_{\bar{g}}^{2}\right)^{1/2}}\right]\cdot v_{\bar{\imath}} \\
  &= h_{ij} + \left\vert \bar{v} \right\vert_{\bar{g}}^{2} \cdot \frac{v_{\bar{\imath}}\phi_{j}}{\phi}+ \frac{\phi^{2}\langle df, d\phi \rangle_{\bar{g}}f_{i}f_{j} - \phi_{i}f_{j} + \phi v_{\bar\imath} f_{\alpha j}v^{\alpha}}{(1- \phi^{2}\left\vert df \right\vert_{\bar{g}}^{2})^{1/2}} \\
  &= h_{ij} + \left\vert \bar{v} \right\vert_{\bar{g}}^{2} \cdot \frac{v_{\bar{\imath}}\phi_{j}}{\phi} + v_{\bar\imath}h(\bar{v},\cdot)_{j} + \left(1-\phi^{2} \left\vert df \right\vert_{\bar{g}}^{2}\right)^{-1/2} \cdot \\
  &\quad \left\{v_{\bar\imath}v^{\alpha}\left(\phi^{2} \langle df, d\phi \rangle_{\bar{g}} f_{\alpha}f_{j} - f_{\alpha} \phi_{j} - \phi_{\alpha} f_{j} \right) + \phi^{2} \langle df, d\phi \rangle_{\bar{g}} f_{i}f_{j} - \phi_{i}f_{j}
  \right\} \\
  &= h_{ij} + v_{\bar\imath} h(\bar{v},\cdot)_{j}- \frac{\phi_{i} v_{\bar\jmath}}{\phi}  \\
 &\quad + \left\vert \bar{v} \right\vert_{\bar{g}}^{2} \frac{v_{\bar\imath}\phi_{j}}{\phi}- v_{\bar\imath}v^{\alpha} \left(\frac{v_{\bar\alpha}\phi_{j}}{\phi} + \frac{\phi_{\alpha}v_{\bar\jmath}}{\phi} \right)  + \phi \bar{v}(\phi) \left(f_{i}f_{j} + \left\vert \bar{v} \right\vert_{\bar{g}}^{2} f_{i}f_{j} \right) \\
  &= h_{ij} + v_{\bar\imath} h(\bar{v},\cdot)_{j}- \frac{\phi_{i} v_{\bar\jmath}}{\phi}
  + \frac{\bar{v}(\phi)}{\phi} \left((1 + \left\vert \bar{v} \right\vert_{\bar{g}}^{2}) \phi^2 f_{i}f_{j} - {v_{\bar\imath}v_{\bar\jmath}}     \right) \\
 &= h_{ij} + v_{\bar\imath} h(\bar{v},\cdot)_{j}- \frac{\phi_{i}
 v_{\bar\jmath}}{\phi}.
\end{align*}
The above calculations follow from our formula for $h$ in equation
\ref{eqn:secfundformnormal}, our definition of $\bar{v}$ in
equations \ref{vdef2} and \ref{vvecdef2}, and the substitution
\[
   v_{\bar\imath} = \frac{\phi f_i}{\left(1-\phi^{2}\left\vert df \right\vert_{\bar{g}}^{2}\right)^{1/2}}
\]
which we use a number of times.

\vspace{.2in}\noindent {\bf Identity \ref{i6}}
\[
   \overline{\mbox{div}}(k)(\bar{v}) = \overline{\mbox{div}}(k(\bar{v},\cdot))
        - \langle h , k \rangle_{\bar{g}}
        - \langle h(\bar{v},\cdot) , k(\bar{v},\cdot) \rangle_{\bar{g}}
        + k\left(\bar{v},\frac{\overline\nabla \phi}{\phi}\right)
\]
\noindent {\it Proof:  } By identity \ref{i5},
\begin{align*}
  \overline{\mbox{div}}\left(k({\bar{v}},\cdot)\right) &= \overline{\mbox{div}}(k)({\bar{v}}) + \langle k_{ij}, v_{\bar\imath\bar{;}j} \rangle_{\bar{g}} \\
  &= \overline{\mbox{div}}(k)({\bar{v}}) + \langle k, h \rangle_{\bar{g}} + \langle k({\bar{v}},\cdot), h({\bar{v}},\cdot) \rangle_{\bar{g}} - k\left(\frac{\overline\nabla \phi}{\phi}, {\bar{v}}\right).
\end{align*}

\vspace{.2in}\noindent {\bf Identity \ref{i7}}
\[
   (\overline\nabla_{\bar{v}} k)(\bar{v},\bar{v}) = \bar{v}(k(\bar{v},\bar{v}))
      - 2 \langle h(\bar{v},\cdot) , k(\bar{v},\cdot) \rangle_{\bar{g}}
      - 2 h(\bar{v},\bar{v}) k(\bar{v},\bar{v})
      + 2 |\bar{v}|_{\bar{g}}^2 \;k\left(\bar{v},\frac{\overline\nabla\phi}{\phi}\right)
\]
\noindent {\it Proof:  } By identity \ref{i5},
\begin{equation*}
   \left( \overline\nabla_{\bar{v}} \bar{v} \right)_{\bar\imath} =
   v^j v_{\bar\imath \bar{;} j}
   = h(\bar{v},\cdot)_{i} + h(\bar{v},\bar{v}) v_{\bar\imath} -
   |\bar{v}|_{\bar{g}}^2 \frac{\phi_i}{\phi}
\end{equation*}
so that by the definition of covariant differentiation of a
symmetric 2-tensor,
\begin{align*}
  (\overline{\nabla}_{{\bar{v}}}k)({\bar{v}},{\bar{v}}) &= {\bar{v}}(k({\bar{v}},{\bar{v}})) - 2k({\bar{v}},\overline{\nabla}_{{\bar{v}}}{\bar{v}}) \\
  &= {\bar{v}}(k({\bar{v}},{\bar{v}})) - 2\langle k({\bar{v}},\cdot),h({\bar{v}},\cdot) \rangle_{\bar{g}} - 2k({\bar{v}},{\bar{v}})h({\bar{v}},{\bar{v}}) + 2 \left\vert {\bar{v}}\right\vert_{\bar{g}}^{2} k\left( {\bar{v}}, \frac{\overline{\nabla} \phi}{\phi}
  \right)
\end{align*}
proving the identity.

\vspace{.2in}\noindent {\bf Identity \ref{i8}}
\begin{eqnarray*}
   \mbox{div}(k)(v) &=& \overline{\mbox{div}}(k(\bar{v},\cdot))
      + \bar{v}(k(\bar{v},\bar{v}))
      + k\left(\bar{v}, \frac{\overline\nabla\phi}{\phi}\right) \\&&
      - \langle h , k \rangle_{\bar{g}}
      - 2\langle h(\bar{v},\cdot) , k(\bar{v},\cdot) \rangle_{\bar{g}}
      + (\mbox{tr}_{\bar{g}}h) k(\bar{v},\bar{v})
\end{eqnarray*}
\noindent {\it Proof:  }  Plugging identities \ref{i6} and \ref{i7}
into identity \ref{i4} and simplifying proves the identity.

\vspace{.2in}\noindent {\bf Identity \ref{mainidentity}}
(\textbf{The Generalized Schoen-Yau Identity})
\begin{eqnarray*}
\bar{R} &=& 16\pi(\mu - J(v)) + \|h-k\|_{\bar{g}}^2+2|q|_{\bar{g}}^2
        - \frac{2}{\phi}\overline{\mbox{div}}(\phi q) \\&&
        + (\mbox{tr}_{\bar{g}}h)^2 - (\mbox{tr}_{\bar{g}}k)^2
        + 2 \bar{v}(\mbox{tr}_{\bar{g}}h - \mbox{tr}_{\bar{g}}k)
        + 2 k(\bar{v},\bar{v})(\mbox{tr}_{\bar{g}}h - \mbox{tr}_{\bar{g}}k)
\end{eqnarray*}
where
\[
   q = h(\bar{v}, \cdot) - k(\bar{v}, \cdot) = h(v, \cdot) - k(v, \cdot) \; .
\]
\noindent {\it Proof:  }First we recall equation
\ref{eqn:startofappendix} derived in section \ref{sec:gsyi},
\begin{eqnarray*}
\bar{R} &=& 16\pi(\mu - J(v))
        + (\mbox{tr}_g h)^2 -  (\mbox{tr}_g k)^2
        - \|h\|_g^2 + \|k\|_g^2   \notag \\ &&
        + 2 \, v(\mbox{tr}_g h) - 2 \, v(\mbox{tr}_g k)
        - 2 \, \mbox{div}(h)(v) + 2 \, \mbox{div}(k)(v).
\end{eqnarray*}
Next, we plug in identities \ref{i1}, \ref{i2}, and \ref{i8}.  Note
that these identities are true for arbitrary symmetric 2-tensors $k$
and hence are true for $h$ as well.  Thus,
\begin{eqnarray*}
\bar{R} &=& 16\pi(\mu - J(v)) \\
   &&+(\mbox{tr}_{\bar{g}}h)^2 -
   \|h\|_{\bar{g}}^2 + 2 h(\bar{v},\bar{v}) \mbox{tr}_{\bar{g}}h -
   2|h(\bar{v},\cdot)|_{\bar{g}}^2 \\
   &&-(\mbox{tr}_{\bar{g}}k)^2 +
   \|k\|_{\bar{g}}^2 - 2 k(\bar{v},\bar{v}) \mbox{tr}_{\bar{g}}k +
   2|k(\bar{v},\cdot)|_{\bar{g}}^2 \\
   && + 2\bar{v}(\mbox{tr}_{\bar{g}}h) + 2\bar{v}(h(\bar{v},\bar{v})) \\
   && - 2\bar{v}(\mbox{tr}_{\bar{g}}k) + 2\bar{v}(k(\bar{v},\bar{v}))\\
    && - 2\overline{\mbox{div}}(h(\bar{v},\cdot))
      - 2\bar{v}(h(\bar{v},\bar{v}))
      - 2h\left(\bar{v}, \frac{\overline\nabla\phi}{\phi}\right) \\&&
      + 2\|h\|_{\bar{g}}^2
      + 4|h(\bar{v},\cdot)|_{\bar{g}}^2
      - 2(\mbox{tr}_{\bar{g}}h) h(\bar{v},\bar{v}) \\
    && + 2\overline{\mbox{div}}(k(\bar{v},\cdot))
      + 2\bar{v}(k(\bar{v},\bar{v}))
      + 2k\left(\bar{v}, \frac{\overline\nabla\phi}{\phi}\right) \\&&
      - 2\langle h , k \rangle_{\bar{g}}
      - 4\langle h(\bar{v},\cdot) , k(\bar{v},\cdot) \rangle_{\bar{g}}
      + 2(\mbox{tr}_{\bar{g}}h) k(\bar{v},\bar{v}).
\end{eqnarray*}
Simplifying and combining terms then gives us that
\begin{eqnarray*}
\bar{R} &=& 16\pi(\mu - J(v)) \\
&& + \|h\|_{\bar{g}}^2 - 2\langle h , k \rangle_{\bar{g}} + \|k\|_{\bar{g}}^2 \\
&& + 2|h(\bar{v},\cdot)|_{\bar{g}}^2
   - 4\langle h(\bar{v},\cdot) , k(\bar{v},\cdot) \rangle_{\bar{g}}
   + 2|k(\bar{v},\cdot)|_{\bar{g}}^2\\
&& - 2\overline{\mbox{div}}(h(\bar{v},\cdot))
   + 2\overline{\mbox{div}} (k(\bar{v},\cdot) )
   - 2h\left(\bar{v}, \frac{\overline\nabla\phi}{\phi}\right)
   + 2k\left(\bar{v}, \frac{\overline\nabla\phi}{\phi}\right)   \\
        && + (\mbox{tr}_{\bar{g}}h)^2 - (\mbox{tr}_{\bar{g}}k)^2
           + 2 \bar{v}(\mbox{tr}_{\bar{g}}h - \mbox{tr}_{\bar{g}}k)
           + 2 k(\bar{v},\bar{v})(\mbox{tr}_{\bar{g}}h - \mbox{tr}_{\bar{g}}k)
\end{eqnarray*}
so that
\begin{eqnarray*}
\bar{R} &=& 16\pi(\mu - J(v)) + \|h-k\|_{\bar{g}}^2+2|q|_{\bar{g}}^2
        - 2\overline{\mbox{div}}(q) - 2q\left(\frac{\overline\nabla \phi}{\phi}\right) \\&&
        + (\mbox{tr}_{\bar{g}}h)^2 - (\mbox{tr}_{\bar{g}}k)^2
        + 2 \bar{v}(\mbox{tr}_{\bar{g}}h - \mbox{tr}_{\bar{g}}k)
        + 2 k(\bar{v},\bar{v})(\mbox{tr}_{\bar{g}}h - \mbox{tr}_{\bar{g}}k)
\end{eqnarray*}
where
\[
   q = h(\bar{v}, \cdot) - k(\bar{v}, \cdot) = h(v, \cdot) - k(v, \cdot) \; .
\]
Note that these two definitions of $q$ exist on the entire
constructed static spacetime and are equal since both $h$ and $k$
are extended trivially in the time direction of the constructed
static spacetime and $v$ and $\bar{v}$ differ only in their time
components.  By the product rule the above equation proves the
identity.

\newpage
\section{The Transformation of the Mean Curvature of the
Boundary}\label{sec:transformationofthemeancurvature}

In this appendix we derive the transformation formula for the mean
curvature of the apparent horizon boundary $\Sigma$ of $(M,g)$ as
approximated by the level sets of $\phi$.  Hence, to be useful, we
need to assume that $\phi = 0$ on $\Sigma$, is strictly positive
elsewhere, and has level sets converging smoothly to $\Sigma$.  The
discussion here supplements the discussion in section
\ref{sec:boundaryconditions}.

More precisely, given
\[
\bar{g} = g + \phi^2 df^2,
\]
we compute the mean curvature $\bar{H}$ of the level sets of $\phi$
in $(M, \bar{g})$ in terms of the mean curvature $H$ of those level
sets in $(M,g)$ and $f$ and $\phi$.

\begin{identity}\label{meancurvident}
\textbf{Transformation of Mean Curvature Identity}

The mean curvature with respect to $\bar{g} = g + \phi^2 df^2$ of a
level set $\Sigma$ of $\phi$ is given by
\begin{eqnarray*}
   \bar{H} &=& \left(\frac{1 + \phi^2|\nabla f|_g^2}{1 + \phi^2|\nabla_\Sigma
   f|_g^2}\right)^{1/2} \left[ (H - \II(T,T)) - (\tr_\Sigma^g h - h(T,T))\langle \nu,v \rangle_g
   \right] \\&+&
   \left(\frac{1 + \phi^2|\nabla_\Sigma f|_g^2}{1 + \phi^2|\nabla f|_g^2}\right)^{1/2}
   \frac12  \nu(\phi^2) |\nabla_\Sigma f|_g^2 (1 - |T|_g^2)
\end{eqnarray*}
where $\nabla$ is the gradient with respect to $g$,
$\nabla_\Sigma$ is the gradient with respect to $g$ restricted to
$\Sigma$, $\II$ is the second fundamental form of $\Sigma$ in
$(M,g)$ (so that $H = \tr_\Sigma(\II)$), $\nu$ is the outward unit
normal vector to $\Sigma$ in $(M,g)$, and
\begin{eqnarray*}
   v &=& \frac{\phi \nabla f}{(1+\phi^2|\nabla f|_g^2)^{1/2}} ,\\
   T &=& \frac{\phi \nabla_\Sigma f}{(1+\phi^2|\nabla f|_g^2)^{1/2}} =
   \mbox{tan}_\Sigma(v) ,\\
   h &=& \frac{\phi \mbox{Hess}f + (df \otimes d\phi + d\phi \otimes df)}
              { (1+\phi^2|df|^2_g)^{1/2}}.
\end{eqnarray*}
\end{identity}

\noindent {\it Proof:} To derive this identity, it is convenient to
let $\gamma = \phi df$ so that $\bar{g} = g + \gamma \otimes
\gamma$.  We also define $\gamma_{\tan} = \tan_{\Sigma}^{g}\gamma$.
Then it is a short exercise to verify that
\begin{equation}
  \bar{{\nu}} = \left(\frac{1+\left\vert \gamma \right\vert^{2}_{g}}{1+ \left\vert \gamma_{\tan} \right\vert^{2}_{g}} \right)^{1/2} \left( {\nu} - \frac{\gamma({\nu})\gamma^{\ast g}}{1+\left\vert \gamma \right\vert^{2}_{g}} \right)
\end{equation}
is the unit normal to $\Sigma$ in $(M,\bar{g})$, where $\gamma^{\ast
g}$ is defined to be the vector dual to the covector $\gamma$ with
respect to $g$ (which of course is $\phi \nabla f$).  The above
formula for $\bar{\nu}$ has the property that $|\bar{\nu}|_{\bar{g}}
= 1$ and $\langle \bar{\nu}, T \rangle_{\bar{g}} = 0$ for all
vectors $T$ tangent to $\Sigma$.

To compute the mean curvature of $\Sigma$ with respect to
$(M^n,\bar{g})$ at a point p, choose a coordinate chart in a
neighborhood of $p$ so that the first $n-1$ coordinate directions
are tangent to $\Sigma$.  Our convention is that for a sphere in
$R^n$ with the standard flat metric, we choose the normal vector to
be the one pointing outwards and the mean curvature to be positive.
Hence, by identity \ref{i3},

\begin{align*}
  \bar{H} &= - \sum_{i,j = 1}^{n-1} \bar{g}^{ij} \left\langle \bar{{\nu}}, \overline{\nabla}_{\bar{\partial}_{i}} \bar{\partial}_{j} \right\rangle_{\bar{g}}
  = -\sum_{i,j=1}^{n-1} \sum_{\theta = 1}^{n} \bar{g}^{ij} \left\langle \bar{{\nu}}, \overline{\Gamma}_{ij}^{\ \ \theta} {\partial}_{\theta} \right\rangle_{\bar{g}} \\
  &= - \sum \bar{g}^{ij} \left\langle \bar{{\nu}}, \ \left(\Gamma_{ij}^{\ \ \theta} + h_{ij}v^{\theta} - \phi f_{i}f_{j} \phi^{\bar{\theta}} \right){\partial}_{\theta} \right\rangle_{\bar{g}} \\
  &= -\sum \bar{g}^{ij} \left(\Gamma_{ij}^{\ \ \theta} + h_{ij}v^{\theta} - \phi f_{i}f_{j} \phi^{\bar{\theta}} \right) \left[ \langle \bar{{\nu}}, \partial_{\theta} \rangle_{g} + \gamma(\bar{{\nu}}) \gamma(\partial_{\theta}) \right] \\
  \notag &= -\sum \bar{g}^{ij} \left(\Gamma_{ij}^{\ \ \theta} + h_{ij}v^{\theta} - \phi f_{i}f_{j} \phi^{\bar{\theta}} \right) \left( \frac{ 1+ \left\vert \gamma \right\vert^{2}_{g} }{1+\left\vert \gamma_{\tan} \right\vert^{2}_{g} } \right)^{1/2} \\
  &\quad \cdot \left[ \left\langle {\nu} - \frac{\gamma({\nu})\gamma^{\ast g}}{1 + \left\vert \gamma \right\vert^{2}_{g}}, \ \partial_{\theta} \right\rangle_{g} + \gamma\left( {\nu} - \frac{\gamma({\nu})\gamma^{\ast g}}{ 1+ \left\vert \gamma \right\vert^{2}_{g}} \right)\gamma(\partial_{\theta}) \right] \\
  \notag &= -\sum \bar{g}^{ij} \left(\Gamma_{ij}^{\ \ \theta} + h_{ij}v^{\theta} - \phi f_{i}f_{j} \phi^{\bar{\theta}} \right) \left( \frac{ 1+ \left\vert \gamma \right\vert^{2}_{g} }{1+\left\vert \gamma_{\tan} \right\vert^{2}_{g} } \right)^{1/2} \\
  &\quad \cdot \left[\langle {\nu}, \partial_{\theta} \rangle_{g} - \frac{\gamma({\nu})\gamma(\partial_{\theta})}{1+ \left\vert \gamma \right\vert^{2}_{g}} + \gamma({\nu})\gamma(\partial_{\theta})\left(1 - \frac{\left\vert \gamma \right\vert^{2}_{g}}{1+ \left\vert \gamma \right\vert^{2}_{g}} \right) \right]
\end{align*}
so that
\begin{equation*}
  \bar{H} = - \sum \bar{g}^{ij} \left(\Gamma_{ij}^{\ \ \theta} + h_{ij}v^{\theta} - \phi f_{i}f_{j} \phi^{\bar{\theta}} \right) \langle {\nu}, \partial_{\theta} \rangle_{g} \left( \frac{ 1+ \left\vert \gamma \right\vert^{2}_{g} }{1+\left\vert \gamma_{\tan} \right\vert^{2}_{g} }
  \right)^{1/2}.
\end{equation*}
Substituting
\begin{align*}
  \phi^{\bar{\theta}} = \bar{g}^{\theta k}\phi_{k} = \left(g^{\theta k} - \frac{\gamma^{\theta}\gamma^{k}}{1+\left\vert \gamma \right\vert^{2}_{g} } \right)\phi_{k}
  =\phi^{\theta} - \frac{\langle \gamma, d\phi \rangle_{g}}{1+ \left\vert \gamma \right\vert^{2}_{g}} \gamma^{\theta}
\end{align*}
we get that
\begin{align*}
  \notag \bar{H} &= - \left( \frac{ 1+ \left\vert \gamma \right\vert^{2}_{g} }{1+\left\vert \gamma_{\tan} \right\vert^{2}_{g} } \right)^{1/2} \sum \left( g^{ij} - \frac{\gamma^{i} \gamma^{j}}{1+ \left\vert \gamma \right\vert^{2}_{g}} \right) \\
  &\quad \cdot \left( \Gamma_{ij}^{\ \ \theta} + h_{ij}v^{\theta} - \phi f_{i}f_{j} \left( \phi^{\theta} - \frac{\langle \gamma, d\phi \rangle_{g}}{1+ \left\vert \gamma \right\vert^{2}_{g}} \gamma^{\theta} \right) \right) \langle {\nu}, \partial_{\theta} \rangle_{g} \\
  \notag &= - \left( \frac{ 1+ \left\vert \gamma \right\vert^{2}_{g} }{1+\left\vert \gamma_{\tan} \right\vert^{2}_{g} } \right)^{1/2} \sum \left( g^{ij} - \frac{\gamma^{i} \gamma^{j}}{1+ \left\vert \gamma \right\vert^{2}_{g}} \right) \\
  &\quad \cdot \left\langle {\nu}, \ \nabla_{\partial_{i}} \partial_{j} + h_{ij}v - \phi f_{i}f_{j} \left( \nabla \phi - \frac{\langle \gamma, d\phi \rangle_{g}}{1+ \left\vert \gamma \right\vert^{2}_{g}} \gamma^{\ast g} \right) \right\rangle_{g} \\
  \notag &= \left( \frac{ 1+ \left\vert \gamma \right\vert^{2}_{g} }{1+\left\vert \gamma_{\tan} \right\vert^{2}_{g} } \right)^{1/2} \left[ \left(H - \frac{\II(\gamma_{\tan}, \gamma_{\tan})}{1+\left\vert \gamma \right\vert^{2}_{g}} \right) - \left(\left(\tr_{\Sigma}^{g} h \right) - \frac{ h(\gamma_{\tan},\gamma_{\tan})}{1 + \left\vert \gamma \right\vert^{2}_{g}} \right) \langle {\nu}, v \rangle_{g} \right. \\
  &\quad \left. + \phi \left(\left\vert \nabla_{g}^{\Sigma} f \right\vert^{2}_{g} - \frac{\langle df, \gamma_{\tan} \rangle^{2}_{g}}{1 + \left\vert \gamma \right\vert^{2}_{g}} \right) \left({\nu}(\phi) - \frac{\gamma({\nu}) \langle \gamma, d\phi \rangle_{g}}{1 + \left\vert \gamma \right\vert^{2}_{g} } \right)
  \right].
\end{align*}
Substituting $\gamma = \phi \ df$, we get
\begin{align}
  \notag \bar{H} &= \left( \frac{1+ \phi^{2} \left\vert \nabla f \right\vert^{2}_{g} }{1 + \phi^{2} \left\vert \nabla_{\Sigma} f \right\vert^{2}_{g}} \right)^{1/2}
  \left\{ \left( H - \frac{\phi^{2} \II(\nabla_{\Sigma} f, \nabla_{\Sigma} f ) }{ 1+ \phi^{2} \left\vert \nabla f \right\vert^{2}_{g} }
  \right)\right. \\ \notag
  &\quad - \left(\tr_{\Sigma}^{g} h - \frac{\phi^{2} h(\nabla_{\Sigma} f, \nabla_{\Sigma} f)}{1+ \phi^{2} \left\vert \nabla f \right\vert^{2}_{g}} \right)\langle {\nu}, v \rangle_{g}  \\
  &\quad \left. + \phi \left\vert \nabla_{\Sigma} f \right\vert^{2}_{g} \left(1 - \frac{\phi^{2} \left\vert \nabla_{\Sigma} f \right\vert^{2}_{g}}{ 1  + \phi^{2} \left\vert \nabla f\right\vert^{2}_{g}} \right) \left( {\nu}(\phi) - \frac{\phi^{2} {\nu}(f) \langle df, d\phi \rangle_{g}}{1 + \phi^{2} \left\vert \nabla f \right\vert^{2}_{g} }\right)\right\}
\end{align}

So far we have not used the assumption that $\Sigma$ is a level set
of $\phi$, so the above formula would be of interest if one wanted
to analyze the mean curvatures of a family of surfaces converging to
the boundary other than the level sets of $\phi$.  Since in our case
$\Sigma$ is a level set of $\phi$, it follows that $\nabla\phi \ ||
\ {\nu}$ and $\nabla \phi = {\nu}(\phi) {\nu}$.  Hence,
\begin{align*}
{\nu}(\phi) - \frac{\phi^{2} {\nu}(f) \langle df, d\phi
\rangle_{g}}{1 + \phi^{2} \left\vert \nabla f \right\vert^{2}_{g} }
=& {\nu}(\phi) \left(1 - \frac{\phi^2 {\nu}(f)^2}{1 + \phi^{2}
\left\vert \nabla f
\right\vert^{2}_{g}} \right) \\
=&  {\nu}(\phi) \left( \frac{1 + \phi^2 |\nabla_\Sigma f|_g^2}{1 +
\phi^{2} \left\vert \nabla f \right\vert^{2}_{g}}  \right).
\end{align*}
Thus, if we recall that
\begin{equation*}
 v = \frac{\phi \nabla f}{\left(1 + \phi^{2} \left\vert
\nabla f \right\vert^{2}_{g} \right)^{1/2}}
\end{equation*}
and we let
\begin{equation*}
\displaystyle T = \frac{\phi \nabla_{\Sigma} f}{\left(1 + \phi^{2}
\left\vert \nabla f \right\vert^{2}_{g} \right)^{1/2}} \;\;\;\left(
= \tan_{\Sigma}(v) \right)
\end{equation*}
we get that
\begin{align}
  \notag \bar{H} &= \left(\frac{1+\phi^{2} \left\vert \nabla f \right\vert^{2}_{g}}{1 + \phi^{2} \left\vert \nabla_{\Sigma} f \right\vert^{2}_{g}} \right)^{1/2} \left[ (H- \II(T,T)) - \left(\tr_{\Sigma}^{g} h - h(T,T) \right) \langle {\nu}, v \rangle_{g} \right] \\
  &\quad + \left(\frac{1+\phi^{2} \left\vert \nabla_{\Sigma} f \right\vert^{2}_{g}}{1 + \phi^{2} \left\vert \nabla f \right\vert^{2}_{g}} \right)^{1/2} \cdot \frac{1}{2} {\nu}(\phi^{2}) \cdot \left\vert \nabla_{\Sigma} f \right\vert^{2}_{g} \cdot \left(1 - \left\vert T \right\vert^{2}_{g} \right)
\end{align}
as claimed in the identity.

The purpose of including this identity in this paper is to help
those who want to study the existence theories of one of the systems
of equations described in this paper, such as the Jang - zero
divergence equations.  The above identity may be useful for
understanding boundary behavior.

In particular, to reduce the Penrose conjecture to the Riemannian
Penrose case, it is necessary for $\bar{H} \le 0$ (or something
very close to this).  Typically, one would even expect $\bar{H} =
0$. The cases of blowup, blowdown, or bounded behavior everywhere
are discussed in section \ref{sec:boundaryconditions} from a
different point of view.  In the cases of blowup or blowdown the
level sets of $f$ were used instead of those of $\phi$.

However, if we are going to allow mixed blowup and blowdown behavior
on generalized apparent horizons, then we can no longer use the
level sets of $f$ everywhere.  However, since $\phi$ is always
assumed to go to zero on the boundary, it is natural to study the
mean curvatures of the level sets of $\phi$. Since an existence
theory is necessary before we can make very many conclusions about
boundary behavior, we restrict this final discussion to a few
observations.

At points on the boundary in the interior of a blowup region, a
blowdown region, or a bounded behavior region, it is plausible that
$v = -\nu$, $v = \nu$, or $v=0$, respectively, in the limit as the
level sets of $\phi$ approach the boundary. In those three cases, it
follows that $T=0$. If we further assume that the second term in the
formula for $\bar{H}$ can be shown to be zero, we then get
\begin{align}
  \notag \bar{H} &= \left(\frac{1+\phi^{2} \left\vert \nabla f \right\vert^{2}_{g}}{1 + \phi^{2} \left\vert \nabla_{\Sigma} f \right\vert^{2}_{g}} \right)^{1/2}
  \left[ H - \left(\tr_{\Sigma}^{g} h  \right) \langle {\nu}, v \rangle_{g}
  \right].
\end{align}

The term in brackets is then $[H + \left(\tr_{\Sigma}^{g}
h\right)]$, $[H - \left(\tr_{\Sigma}^{g} h  \right)]$, and $[H]$ in
those three respective cases.  Modulo possible issues with taking
limits, one could then use the generalized Jang equation to conclude
that
\begin{equation}
 0 = \tr_{\bar{g}}(h-k) = \tr_g(h-k) - (h-k)(v,v) = \tr^g_\Sigma(h-k)
\end{equation}
in the case of either blowup or blowdown, since $\bar{g}^{ij} =
g^{ij} - v^i v^j$.  Hence, the term in brackets equals zero in the
three respective cases of a local future apparent horizon, a local
past apparent horizon, or a local future and past apparent horizon,
as desired.  One would then need to show that the term in front of
the brackets still allows one to conclude $\bar{H} = 0$ in the
limit, even when it diverges as the level sets of $\phi$ approach
the boundary.

Of course the really tricky part is understanding points on the
boundary where every open set around the point contains two or
more of blowup, blowdown, and bounded behavior.  We offer the
above formula for the mean curvature $\bar{H}$ of the level sets
of $\phi$ in case it is helpful to others who approach this
problem.



\newpage
\begin{thebibliography}{10}

\bibitem{ANDERSSONMETZGER}
L. Andersson, and J. Metzger, \textit{The area of horizons and the
trapped region}, preprint, arXiv:0708.4252, 2007.


\bibitem{BRAY}
H.L. Bray, \textit{Proof of the Riemannian Penrose inequality using
the positive mass theorem,} J. Differential Geom. $\mathbf{59}$
(2001), 177-267.

\bibitem{BRAYsurvey}
H.L. Bray, \textit{Black holes, geometric flows, and the Penrose
inequality in general relativity}, Notices of the AMS {\bf 49}
(2002), 1372-1381.

\bibitem{BRAYCHRUSCIEL}
H.L. Bray and P.T. Chrusciel, \textit{The Penrose inequality} in
\textit{The Einstein equations and the large scale behavior of
gravitational fields (50 years of the Cauchy problem in general
relativity)}, edited by P.T. Chrusciel, H. Friedrich, Birkhaeuser,
Basel, 2004, pp. 39-70.  ESI preprint 1390,
http://arxiv.org/abs/gr-qc/0312047

\bibitem{BRAYLEE}
H.L. Bray, and D. Lee, \textit{On the Riemannian Penrose inequality
in dimensions less than 8}, preprint, arXiv:0705.1128, 2008,
(accepted by the Duke Mathematical Journal).

\bibitem{BUNTINGALAM}
G. Bunting, and A. Masood-ul-Alam, \textit{Nonexistence of
multiple black holes in asymptotically Euclidean static vacuum
space-time}, Gen. Relativity Gravitation $\mathbf{19}$ (1987), no.
2, 147-154.

\bibitem{CDGH}
P.T. Chrusciel, E. Delay, G. Galloway, and R. Howard, {\it
Regularity of horizons and the area theorem}, Annales Henri Poincare
{\bf 2} (2001), 109-178, gr-qc/0001003.

\bibitem{EICHMAIR}
M. Eichmair, \textit{Existence, regularity, and properties of
generalized apparent horizons}, preprint, arXiv:0805.4454, 2008.

\bibitem{GEROCH}
R. Geroch, \textit{Energy extraction}, Ann. New York Acad. Sci.
$\mathbf{224}$ (1973), 108-117.

\bibitem{HAWKINGELLIS}
S. Hawking, and G. Ellis, \textit{The large scale structure of
space-time}, Cambridge Monographs on Mathematical Physics,
Cambridge, 1973.

\bibitem{HERZLICH}
M. Herzlich, \textit{A Penrose-like inequality for the mass of
Riemannian asymptotically flat manifolds}, Commun. Math. Phys.
$\mathbf{188}$ (1997), 121-133.

\bibitem{HEUSLER}
M. Heusler, \textit{Black hole uniqueness theorems}, Cambridge
Lecture Notes in Physics (Cambridge univ. Press) Vol. 6, Cambridge,
1996.

\bibitem{HOROWITZ}
G.T. Horowitz, \textit{The positive energy theorem and its
extensions} in \textit{Asymptotic behavior of mass and spacetime
geometry}, Springer Lecture Notes in Physics {\bf 202} Ed. F.
Flaherty (Springer, New York, 1984), 1-20.

\bibitem{HUISKENILMANEN}
G. Huisken, and T. Ilmanen, \textit{The inverse mean curvature
flow and the Riemannian Penrose inequality}, J. Differential Geom.
$\mathbf{59}$ (2001), 353-437.

\bibitem{JANG}
P.-S. Jang, \textit{On the positivity of energy in general
relativity}, J. Math. Phys. $\mathbf{19}$ (1978), 1152-1155.

\bibitem{JANGWALD}
P.-S. Jang, and R. Wald, \textit{The positive energy conjecture
and the cosmic censorship hypothesis}, J. Math. Phys.
$\mathbf{18}$ (1977), 41-44.

\bibitem{JAUREGUI}
J. Jauregui, thesis, Duke University (in preparation).

\bibitem{KHURI}
M. Khuri, \textit{Nonexistence of generalized apparent horizons in
Minkowski space}, Class. Q. Grav. $\mathbf{26}$ (2009), 078001.

\bibitem{KHURI2}
M. Khuri, \textit{A Penrose-like inequality for general initial
data sets}, Commun. Math. Phys., to appear, 2009.

\bibitem{MARS}
M. Mars, \textit{The Penrose inequality in general relativity} (in
preparation).

\bibitem{MARSSENOVILLA}  M. Mars, and J. M. M. Senovilla,
\textit{Trapped surfaces and symmetries}, Class. Q. Grav.
$\mathbf{20}$ (2003), L293-L300.

\bibitem{MEEKSSIMONYAU}
W. Meeks III, L. Simon, and S.-T. Yau, \textit{Embedded minimal
surfaces, exotic spheres, and manifolds with positive Ricci
curvature}, Ann. of Math. $\mathbf{116}$ (1982), 621-659.

\bibitem{ONEILL}
B. O'Neill, \textit{Semi-Riemannian geometry with applications to
relativity}, Pure and Applied Mathematics (Academic Press) Vol. 103,
New York, 1983.

\bibitem{PENROSE1}
R. Penrose, \textit{Gravitational collapse: the role of general
relativity}, Rivista del Nuovo Cimento Numero Speciale $\mathbf{1}$
(1969), 252-276.

\bibitem{PENROSE2}
R. Penrose, \textit{Naked singularities}, Ann. N. Y. Acad. Sci.
$\mathbf{224}$ (1973), 125-134.

\bibitem{SCHOENYAU1}
R. Schoen, and S.-T. Yau, \textit{On the proof of the positive mass
conjecture in general relativity}, Commun. Math. Phys. $\mathbf{65}$
(1979), 45-76.

\bibitem{SCHOENYAU2}
R. Schoen, and S.-T. Yau, \textit{Proof of the positive mass
theorem II}, Commun. Math. Phys. $\mathbf{79}$ (1981), no. 2,
231-260.

\bibitem{SCHOENYAU3}
R. Schoen, and S.-T. Yau, \textit{Positivity of the total mass of
a general space-time}, Phys. Rev. Lett. $\mathbf{43}$ (1979),
1457-1459.

\bibitem{SCHOENYAU4}
R. Schoen, and S.-T. Yau, \textit{Existence of incompressible
minimal surfaces and the topology of three dimensional manifolds
with nonnegative scalar curvature}, Ann. of Math. $\mathbf{110}$
(1979), 127-142.

\bibitem{SENOVILLA} J. M. M. Senovilla, \textit{Classification of
spacelike surfaces in spacetime}, Class. Q. Grav. $\mathbf{24}$
(2007), 3091-3124.

\bibitem{STREETS}
J. Streets, \textit{Applications of inverse mean curvature flow to
negative mass singularities} (in preparation).

\bibitem{WALD1}
R. Wald, \textit{General relativity}, Univ. Chicago Press, Chicago,
1984.

\bibitem{WITTEN}
E. Witten, \textit{A new proof of the positive energy theorem},
Commun. Math. Phys. $\mathbf{80}$ (1981), 381-402.

\end{thebibliography}
\end{document}